\documentclass[11pt,paper=a4]{article}
 \usepackage{indentfirst, latexsym,bm}
 \usepackage{amsmath}
 \usepackage{pifont}
 \usepackage{amsfonts}
 \usepackage{mathrsfs}
 \usepackage{array}
 \usepackage{multirow}
 \usepackage{graphicx}
 \usepackage{subfigure}
 \usepackage{picinpar}
 \RequirePackage[colorlinks,citecolor=blue,urlcolor=blue,linkcolor=blue]{hyperref}
 \RequirePackage[numbers]{natbib}
 \usepackage{enumitem} 

 \usepackage{booktabs}
 \usepackage{threeparttable}
 \usepackage{mathrsfs,amsfonts,amsmath}
 \usepackage{color}

 \makeatletter
 \def\namedlabel#1#2{\begingroup
	#2%
	\def\@currentlabel{#2}%
	\phantomsection\label{#1}\endgroup
}
\makeatother

 \newcommand\email[1]{\href{mailto:#1}{ \nolinkurl{#1}}}

 \renewcommand{\theequation}
 {\arabic{section}.\arabic{equation}}

 \newtheorem{theorem}{Theorem}[section]
 \newtheorem{definition}[theorem]{Definition}
 \newtheorem{lemma}[theorem]{Lemma}
 \newtheorem{corollary}[theorem]{Corollary}
 \newtheorem{proposition}[theorem]{Proposition}
 \newtheorem{remark}[theorem]{Remark}
 \newtheorem{condition}[theorem]{Condition}
 \newtheorem{example}{Example}[section]

 \def\blemma{\begin{lemma}}\def\elemma{\end{lemma}}
 \def\bproposition{\begin{proposition}}\def\eproposition{\end{proposition}}
 \def\ttheorem{\begin{theorem}}\def\etheorem{\end{theorem}}
 \def\bcorollary{\begin{corollary}}\def\ecorollary{\end{corollary}}
 \def\bremark{\begin{remark}}\def\eremark{\end{remark}}
 \def\bcondition{\begin{condition}}\def\econdition{\end{condition}}

 \def\benumerate{\begin{enumerate}}\def\eenumerate{\end{enumerate}}
 \def\bitemize{\begin{itemize}}\def\eitemize{\end{itemize}}

 \def\beqlb{\begin{eqnarray}}\def\eeqlb{\end{eqnarray}}
 \def\beqnn{\begin{eqnarray*}}\def\eeqnn{\end{eqnarray*}}
 \def\ar{\!\!\!&}

 \def\proof{\noindent{\it Proof.~~}}\def\qed{\hfill$\Box$\medskip}

\begin{document} 
 	
 \title{\bf \LARGE Functional Limit Theorems for Hawkes Processes}
 
 \author{Ulrich Horst\footnote{Department of Mathematics and School of Business and Economics,  Humboldt-Universit\"at zu Berlin, Unter den Linden 6, 10099 Berlin; email: horst@math.hu-berlin.de. Horst gratefully acknowledges support from DFG CRC/TRR 388 ``Rough Analysis, Stochastic Dynamics and Related Fields", Project B02. }
 \quad\ and\quad 
 Wei Xu\footnote{School of Mathematics and Statistics, Beijing Institute of Technology, Beijing 100081, China; email: xuwei.math@gmail.com} 
 } 
 	  
 \maketitle

 \begin{abstract}
 We prove that the long-run behavior of Hawkes processes is fully determined by the average number {\sl and} the dispersion of child events. 
 For subcritical processes we provide FLLNs and FCLTs under minimal conditions on the kernel of the process with the precise form of the limit theorems depending strongly on the dispersion of child events. 
 For a critical Hawkes process with weakly dispersed child events, functional central limit theorems do not hold. 
 Instead, we prove that the rescaled intensity processes and rescaled Hawkes processes behave like CIR-processes without mean-reversion, respectively integrated CIR-processes. 
 We provide the rate of convergence by establishing an upper bound on the Wasserstein distance between the  distributions of rescaled Hawkes process and the corresponding limit process. 
 By contrast, critical Hawkes process with heavily dispersed child events share many properties of subcritical ones. 
 In particular, functional limit theorems hold. However, unlike  subcritical processes critical ones with heavily dispersed child events display long-range dependencies. 
 
 \bigskip
 
 \noindent {\it MSC 2020 subject classifications:} Primary 60F17, 60G55; secondary 60J80
  
 \smallskip
  
 \noindent  {\it Keywords and phrases:} Hawkes process, functional limit theorem, regular variation, convergence rate.

 \end{abstract}
  	

   \section{Introduction}

 A Hawkes process $N:=\{N(t):t\geq 0\}$ is a random point process that models self-exciting arrivals of random events. Its \textit{intensity} $\Lambda:=\{\Lambda(t):t\geq 0\}$ is usually of the form
 \beqlb\label{HawkesDensity}
 \Lambda(t):= \mu(t)+ \sum_{0<\tau_i<t} \phi(t-\tau_i) = \mu(t) + \int_{(0,t)} \phi(t-s)N(ds), 
 \eeqlb
 for some \textit{immigration density} $\mu \in    L^1_{\rm loc} (\mathbb{R}_+;\mathbb{R}_+)$ that captures the immigration of exogenous events, and some \textit{kernel}
 $\phi\in    L^1 (\mathbb{R}_+;\mathbb{R}_+)$ that captures the self-exciting impact of past events on the arrivals of future events. 
 The random variable $\tau_i$ denotes the arrival time of the $i$-th event, for each $i \in \mathbb N$. 
 
 A Hawkes process can be reconstructed as a Poisson cluster process associated to an age-dependent branching process with an average offspring number $m:=\|\phi\|_{L^1}$. 
 The criticality of branching processes gives rise to three regimes.
 Depending on the average number of offsprings, a Hawkes process is called subcritical/stationary $(m<1)$, {critical/quasi-stationary} $(m=1)$ or {super-critical/non-stationary} $(m>1)$.   
 
 In this paper we analyse the asymptotic behavior of subcritical and critical Hawkes processes by establishing FLLNs and FCLTs (classical and non-clasical) for such processes.\footnote{Supercritical Hawkes processes grow to infinity exponentially fast in which case FLLNs and FCLTs cannot be established. However, under the light-tailed condition on the kernel $\phi$, it is possible to obtain weak LLNs and CLTs by using the asymptotic results for the corresponding supercritical Crump-Mode-Jagers branching processes given in \cite{Jagers1975,Sagitov1995,Xu2021}.} 

\subsection{Literature review}  
  
 First introduced by Hawkes in \cite{Hawkes1971a,Hawkes1971b} to understand cross-dependencies between earthquakes and their aftershocks, Hawkes processes have been generalized in many directions in recent years including marked Hawkes processes and point measures \cite{BremaudMassoulie2002,Ogata1988}, nonlinear Hawkes processes \cite{BremaudMassoulie2001,Reynaud-BouretSophieSchbath2010}, infinite-dimensional processes \cite{BieleckiJakubowskiNieweglowski2022,HorstXu2019}, and nearly unstable processes \cite{JaissonRosenbaum2015,JaissonRosenbaum2016,Xu2021}. 
 
 Hawkes processes and their generalizations have become a powerful tool to model a variety of phenomena in biology and neuroscience \cite{HorstXu2021,Johnson1996,Pernice2011,Reynaud-BouretSophieSchbath2010}, sociology and criminology \cite{Blundell2012,Crane2008,Mohler2011}, seismology \cite{Bremaud2010,Ogata1988} and, in particular, finance. Applications in finance range from intraday transaction dynamics \cite{Bauwens2009,Bowsher2007} to asset price processes \cite{BacryDelattreHoffmannMuzy2013a,BacryDelattreHoffmannMuzy2013}and  rough volatility modeling \cite{ElEuchFukasawaRosenbaum2018,HorstXu2022,JaissonRosenbaum2015,JaissonRosenbaum2016}, and from limit order book modelling  \cite{HorstXu2019,Large2007} to financial contagion \cite{SahaliaCacho-DiazLaeven2015,Giesecke2011,JorionZhang2009} and exchange rate dynamics \cite{Hewlett2006}. We refer to  \cite{BacryMastromatteoMuzy2015} for a review of Hawkes processes and their applications.

The more applied literature on Hawkes processes is accompanied by an increasing theoretical literature that studies {\sl microscopic} and {\sl macroscopic} properties of Hawkes processes and their generalizations. Three well-known microscopic properties, including statistical characterization, cluster representation and probability generating function were first established in \cite{Hawkes1971a,Hawkes1971b,HawkesOakes1974} for subcritical Hawkes processes and extended recently in \cite{BacryMuzy2014,GaoZhu2018c} to general processes. The genealogy and event cascades of a Hawkes process with exponential kernel were studied in \cite{HorstXu2022} by exploring the intrinsic branching structure.  
  
Macroscopic properties are usually established by proving scaling limit or laws of large numbers and central limit theorems. Hawkes and Oakes  \cite{HawkesOakes1974} were the first to establish a central limit theorem (CLT) for stationary Hawkes processes whose kernels satisfy the integrability condition $\int_0^\infty t\, \phi(t)dt<\infty$. Their result was later generalized by Bacry et al. \cite{BacryDelattreHoffmannMuzy2013} who established a functional law of large numbers (FLLN) and a functional central limit theorem (FCLT) under the weaker integrability condition 
 \beqlb\label{Con.Bacry2013}
 \int_0^\infty \sqrt{t}\, \phi(t)dt<\infty. 
 \eeqlb
 A FCLT and a large deviation principle (LPD) for Hawkes processes with exponential kernel and large exogenous density have been established by Gao and Zhu \cite{GaoZhu2018b,GaoZhu2018a}. 
 LDPs for nonlinear and mean-field limits of Hawkes processes have been established in \cite{Zhu2014,Zhu2015} and \cite{GaoGaoZhu2023,GaoZhu2023} respectively.  
     
Limit theorems and LDPs  for {\sl marked} Hawkes processes have also been studied by many authors. LDPs for marked Hawkes processes were first established in \cite{StabileTorrisi2010,Zhu2013a}; their results have recently been extended to path-wise LDPs in \cite{DenkertXu2023}. CLTs for marked Hawkes processes were first established in \cite{KarabashZhu2015,YaoXiao2018}. Under a light-tailed condition on the kernel, FLLNs and FCLTs for Hawkes random measures (of which marked Hawkes processes are a special case) have been established in \cite{HorstXu2021}. The results were used to analyze the population dynamics of budding microbes in a host. Scaling limits for marked Hawkes processes with exponential, respectively general light-tailed kernels have been analyzed in  \cite{HorstXu2022}  in  \cite{Xu2021}. 
For asymptotically critical processes Xu \cite{Xu2021b} proved the weak convergence of their rescaled densities to a multi-type continuous-state branching process with immigration. Horst et al.~\cite{HorstXuZhang2023} considered a class of marked Hawkes processes whose kernels are step functions. Under a heavy-tailed condition on the kernels they proved that the rescaled intensities converge weakly to the unique solution of a stochastic Volterra equation driven by a Brownian motion and a Poisson random measure.

%
%
%
%
 
Nearly unstable Hawkes processes\footnote{A nearly unstable Hawkes process is obtain by taking the limit of a subcritical Hawkes processes as $m \uparrow 1$. } with light-tailed kernels were first analyzed by Jaisson and Rosenbaum \cite{JaissonRosenbaum2015}.  They proved the weak convergence of the rescaled intensity to a Feller diffusion - also known as CIR-model in finance - and the convergence of the rescaled point process to the integrated diffusion; their result was extended to multi-variate processes in \cite{Xu2021}. Under a heavy-tailed condition on the kernel, the same authors later proved the weak convergence of the rescaled  point process to the integral of a rough fractional diffusion; see \cite{JaissonRosenbaum2016}. Analogous scaling limits in the multivariate case were established in  \cite{ElEuchFukasawaRosenbaum2018,RosenbaumTomas2021}.  A more refined convergence results has recently been established by Horst et al. \cite{HorstXuZhang2023a}; they proved the weak convergence of the rescaled intensities, instead of their integrals, to a rough fractional diffusion. Our analysis shows that critical and nearly unstable Hawkes processes display very different asymptotic behavior. 

\subsection{Our contribution}

This paper establishes functional and scaling limit theorems for Hawkes processes under minimal conditions on their kernels. The following key questions will be addressed. First, in \cite{BacryDelattreHoffmannMuzy2013} the moment condition (\ref{Con.Bacry2013}) plays an important role in establishing the classical FCLT for subcritical Hawkes processes. Our first key question is hence:
  \begin{enumerate} 	
 	\item[\namedlabel{Q1}{(Q1)}]  \textit{Is the condition (\ref{Con.Bacry2013}) sharp or can FCLTs still be established for subcritical Hawkes processes if this condition fails to hold?}
 \end{enumerate}

The FCLT for subcritical processes uncovers that despite the cross-dependencies on the levels of event arrivals, subcritical Hawkes processes not only do {\sl not} display long-range dependencies, but actually enjoy similar asymptotic properties as Poisson point processes where the event arrivals are entirely uncorrelated. This motivates our second key question:  
   \begin{enumerate} 	
   \item[\namedlabel{Q2}{(Q2)}]  \textit{Which classes of Hawkes processes display long-range dependencies?}
 \end{enumerate}

 It is well known that subcritical Hawkes processes with constant immigration rate converge to a stationary regime in the long run. 
 To the best of our knowledge, no general result on the asymptotic stationarity or non-stationary of critical Hawkes processes has so far been established in the literature.\footnote{The existence of non-trivial stationary Hawkes processes on $\mathbb{R}$ with zero exogenous arrival intensity was proved in \cite{BremaudMassoulie2001} under regular variation conditions on the kernel $\phi$.} Our third key question is thus:  
  \begin{enumerate} 	
 	\item[\namedlabel{Q3}{(Q3)}]  \textit{Are all critical Hawkes processes with constant immigration asymptotically non-stationary?  Do FLLNs and FCLTs hold for non-stationary processes?}
 \end{enumerate}
 
 We provide detailed answers to all three questions, thereby providing a comprehensive analysis of the long-run behavior of both subcritical and critical Hawkes processes. 
 To answer questions \ref{Q1}-\ref{Q3} we first extend the duality method developed in \cite{Jaber2021,JaberLarssonPulido2019,Xu2021b} and establish an explicit exponential-affine representation of the Fourier-Laplace functional  of the process of $(N,\Lambda)$ in term of the unique solution to a certain nonlinear Volterra equation. 
 The functional captures most microscopic properties of  $(N,\Lambda)$ including the finite-dimensional distributions, event cascades, and  the distribution and genealogy of all child events triggered by a mother event. 
 The representation of the Lapalce functional is the first main contribution of this paper.
 
 \subsubsection{Subcritical Hawkes processes}
 
 Our second main contribution is twofold and fully answers question \ref{Q1}. First, we show that the moment condition~(\ref{Con.Bacry2013}) is not sharp; 
 second we show that FCLTs for subcritical Hawkes processes with constant immigration rate can be established under much weaker conditions with the precise form depending on the dispersion of child events. 
 
 More precisely, by exploring the release of excitation of each event on future events, we prove that the long-term behavior of a subcritical Hawkes process is fully determined by the function 
 \beqnn
 {\it \Psi}_1(t):= \int_0^t s \phi(s)ds,\quad t\geq 0
 \eeqnn
 that describes the dispersion of child events. We prove that a classical FCLT holds if and only if $\frac{{\it \Psi}_1(t)}{\sqrt{t}} \to 0$ as $t\to\infty$, which is much weaker than  the condition~(\ref{Con.Bacry2013}). If $\frac{{\it \Psi}_1(t)}{\sqrt{t}} \to {\it \Psi}_*\in (0,\infty)$, a non-classical FCLT holds. In this case, rescaled Hawkes process behaves asymtotically like a Brownian motion with linear drift plus an additional square root function, i.e.,
 \beqnn
 \frac{N(nt)}{n} \sim C_0\cdot t + \frac{C_2}{\sqrt{n}} \cdot B(t) + \frac{C_1}{\sqrt{n}}   \cdot\sqrt{t}. 
 \eeqnn
Finally, if $\frac{{\it \Psi}_1(t)}{\sqrt{t}} \to \infty$ and the integrated tail-function
 \beqnn
 {\it\Phi}(t) := \int_t^\infty \phi(s)ds,\quad t\geq 0,
 \eeqnn
 of the kernel $\phi$ is regularly varying with index $-\alpha \in[-1/2,0]$, a degenerate FCLT holds. In this case, 
 \beqnn
  \frac{N(nt)}{n} \sim C_0\cdot t+ C_3\cdot {\it\Phi}(n) \cdot t^{1-\alpha}. 
 \eeqnn
 
\subsubsection{Critical Hawkes processes}
 
 Whereas subcritical Hawkes processes do not display long-range dependencies we prove that critical Hawkes processes do. This answers question \ref{Q2}. 
 Furthermore, 
 we prove that the long run behavior of critical Hawkes processes  and hence the answer to question \ref{Q3} critically depends on the {\sl dispersion coefficient} 
 \[
 	\sigma:= {\it\Psi}_1(\infty). 
 \] 
 We call a Hawkes process {\it weakly critical} if $\sigma<\infty$ and {\sl strongly critical} if $\sigma=\infty$. It turns out that weakly critical processes are not asymptotically stationary so that neither a classical FLLN nor a classical FCLT can be established for such processes. By contrast, classical FLLNs and FCLTs can be established for strongly critical processes. This answers question \ref{Q3}.

 Our limit results for weakly critical processes complement the work of Jaisson and Rosenbaum \cite{JaissonRosenbaum2015}. They considered scaling limits for sequences of weakly critical processes with light-tailed kernels when the number of average offsprings increases to one. In this case, the limiting model can be described by a CIR model with mean-reversion. By contrast, we consider limit theorems for critical Hawkes processes (with a fixed number of average offsprings). It turns out that in the weakly critical case our limit model can also be described by a CIR model, but {\sl without} mean-reversion. The long run dynamics of our model is hence very different from the one in \cite{JaissonRosenbaum2015}. 

 Our results for strongly critical processes complement the scaling limits for nearly unstable processes with heavy-tailed kernels established in \cite{JaissonRosenbaum2016} and refined in Horst et al \cite{HorstXuZhang2023}. Again, the limiting results are very different. In particular, it does not seem to be possible to obtain the rough fractional diffusion without mean-reversion obtained in  \cite{HorstXuZhang2023,JaissonRosenbaum2016} by scaling a single critical Hawkes processes. 

 In what follows we summarize our main results for critical Hawkes processes in greater detail.


 \subsubsection{Weakly critical Hawkes processes}

 Motivated by scaling limits established in \cite{HorstXu2022}, we prove that after a suitable scaling in time and space, the intensity of a weakly critical Hawkes process can be approximated in law by a CIR model without mean-reversion $\Lambda^*$, and that the Hawkes process $N$ behaves asymptotically as the integrated process of $\Lambda^*$. That is, we prove that
 \beqnn
 \Lambda^{(n)}(t):= \frac{\Lambda(nt)}{n} \sim \Lambda^* (t) 
 \quad \mbox{and}\quad 
 N^{(n)}(t):= \frac{ N(nt)}{n^2} \sim \int_0^t  \Lambda^* (s)ds =: \mathcal{I}_{\Lambda^*}(t), \quad t\geq 0.  
 \eeqnn
 The CIR model without mean-reversion satisfies $\Lambda^*(t) \to \infty$ a.s.~as $t\to\infty$ and   
 \[
 	{\rm Corr}(\Lambda^*(t),\Lambda^*(t+T))\sim C/T \quad  \mbox{as} \quad T\to\infty. 
\]	
The latter suggests that the correlation between increments of the point process $N$ decreases at a rate $\frac{1}{T}$, i.e.~that
 	\beqnn
 	{\rm Corr}\big(N(t,t+1),N(t+T,t+T+1)\big) \sim 
 	{\rm Corr}\big(\Lambda(t),\Lambda(t+T)\big) \sim \frac{C}{T}.
 	\eeqnn
As a result, we cannot expect weakly stationary processes to be asymptotically stationary but we do expect them to display long-range dependencies. 

We also investigate the rate of convergence by establishing  two upper bounds on the distance between the path-distributions of $(\Lambda^{(n)}, N^{(n)})$ and $(\Lambda^*, \mathcal{I}_{\Lambda^*})$ on compact time intervals with respect to the  Fourier-Laplace functional and the Wasserstein distance. 
It turns out that the rate of convergence is fully determined by the decay rate of $\|\mathcal{I}_{R^{(n)}} -\mathcal{I}_{R^*} \|_{L^\infty_T}$ to $0$, where $\mathcal{I}_{R^{(n)}}$ is the integral function of the time-scaled resolvent of $\phi$ and $\mathcal{I}_{R^*}$ is the corresponding limit function. Under a regular variation condition on the integrated tail-function ${\it\Phi}$,  an upper estimate in term of ${\it\Phi}$ on $\|\mathcal{I}_{R^{(n)}} -\mathcal{I}_{R^*} \|_{L^\infty_T}$ is established.
  
 The lack of Markovianity, martingale and regularity properties of kernel $\phi$ renders the analysis of the convergence rate challenging as many standard methods such as Stein's method, Lyapunov function and Markov coupling developed to study convergence rates and ergodicity of Markov processes and martingales cannot be applied in our setting. 
 
 To overcome these difficulties, we use the explicit exponential-affine representations of the Fourier-Laplace functionals of $(\Lambda^{(n)}, N^{(n)})$ and $(\Lambda^*, \mathcal{I}_{\Lambda^*})$ to bound the distance of their distribution with respect to the  Fourier-Laplace functional through the difference of their Fourier-Laplace exponents. Bounds for the latter difference are successfully established by considering the corresponding nonlinear Volterra equations. 
For the distributional distance with respect to the Wasserstein distance, we first use the Plancherel theorem to rewrite the distance in terms of  Fourier-Laplace functionals, for which upper estimates can be established through a series of tail estimates for the Laplace-Stieltes transforms of the integrated kernel function and integrals of functions that describe the dispersion  of child events.

  \subsubsection{Strongly critical Hakes processes}
 
 
 While FLLNs and FCLTs do not hold for weakly critical processes the situation is very different for strongly critical processes. 
 If the arrivals of child events are highly dispersed, then it is almost impossible to trace a child event back to its mother event and we expect the distribution of the Hawkes process to display some form of asymptotic loss of memory. 
 
 While some form of stationarity is expected in the long run and hence FLLNs and FCLTs are expected to hold, we will see that strongly critical Hawkes processes display long-range dependencies.  
 %
%
%
 More precisely, we prove that when the integrated tail-function ${\it\Phi}$ is regularly varying with index $-\alpha\in[-1,0]$,  strongly critical Hawkes processes asymptotically behave like Gaussian processes with long-range dependencies when suitably normalized and hence display long-range dependencies.
 Specifically, we prove that
 \beqnn
 \frac{N(nt)-\mathbf{E}[N(nt)]}{\mathcal{I}^2(n)} \sim C_4\cdot \frac{\sqrt{\mathcal{I}^2_R(n)}}{n} \cdot \int_0^t (t-s)^{\alpha}s^{\alpha/2} dB(s)
 \eeqnn
 where the function $\mathcal{I}^2_R$ denotes the repeated integral of order two of the resolvent of $\phi$. 
 The function $\mathcal{I}^2_R$ is regularly varying with index $1+\alpha\in[1,2]$ and $\mathcal{I}^2_R(n)=o(n^2)$ as $n\to\infty$. Furthermore, if $\mathcal{I}^2_R$ is second-order regularly varying with auxiliary function $A$ that is  regularly varying with index $\rho\leq 0$, then we  prove that 
\beqlb\label{eqn.00}
 \frac{N(nt) }{\mathcal{I}^2(n)}  \sim C_5\cdot  t^{\alpha+1} +C_6\cdot \frac{\sqrt{\mathcal{I}^2_R(n)}}{n} \cdot \int_0^t (t-s)^{\alpha}s^{\alpha/2} dB(s) + C_7\cdot t +  C_8 \cdot t^{\alpha+1}\cdot \int_1^t u^{\rho-1}du.
 \eeqlb

 \medskip
 The remainder of this paper is organized as follows. In Section \ref{Sec.MR}, we first introduce some elementary properties of Hawkes processes including the exponential-affine representation of their Fourier-Laplace functionals. Subsequently we state FLTs for Hawkes processes in different regimes. Detailed proofs for the elementary properties are given in Section~\ref{Sec.EleProp}. The FLLNs and FCLTs for subcritical, weakly critical and strongly critical Hawkes processes are proved separately in Sections~\ref{Sec.SubHP}, \ref{Sec.WeaklyHP} and \ref{Sec.StronglyHP}, respectively. Selected auxiliary results for distances between finite measures and regularly varying functions are presented in Appendix~\ref{Appendix.DM} and \ref{Appendix.RV}, respectively.

\bigskip

 {\bf Notation.}
 Let $\mathbb{R}_+=[0,\infty)$, $\mathbb{N}=\{0,1,2,\cdots  \}$ and $\mathbb{B}$ be a Banach space endowed with a norm $\|\cdot\|_\mathbb{B}$. 
 For $\mathcal{T} \subset \mathbb{R}$ and $k\in\mathbb{N}$, we denote by $C^k(\mathcal{T};\mathbb{B})$ the space of $\mathbb{B}$-valued functions on $\mathcal{T}$ with continuous $k$-derivative. Let $C_b (\mathcal{T};\mathbb{B})$ be the space of bounded, continuous $\mathbb{B}$-valued functions  on $\mathcal{T}$.  
 For each $p\in[1,\infty]$, denote by $L^p(\mathcal{T};\mathbb{B})$ the space of $\mathbb{B}$-valued functions on $\mathcal{T}$ that satisfy $\|f\|_{L^p_\mathcal{T}}^p:=\int_\mathcal{T} \|f(s)\|_\mathbb{B}^p ds<\infty$. We set $L^p_{\rm loc}(\mathbb{R}_+;\mathbb{B}):=\bigcap_{T\geq 0} L^p([0,T];\mathbb{B})$.
 For simplicity, we write $\|f\|_{L^p}$ for $\|f\|_{L^p_\infty}$ and $\|f\|_{L^p_T}$ for $\|f\|_{L^p_{[0,T]}}$ for any $T>0$.
  Let $M(\mathcal{T};\mathbb{B})$ be the space of  all countably-additive, $\mathbb{B}$-valued finite vector measures on $\mathcal{T}$ and $M(\mathbb{R}_+;\mathbb{B}):=\cup_{T\geq 0} M([0,T];\mathbb{B})$ be the space of $\sigma$-finite measures on $\mathbb{R}_+$. Specially, $M(\mathcal{T};\mathbb{B})$ turns to be the space of all  finite measures, signed measures and complex measures on $\mathcal{T}$ when $\mathbb{B}=\mathbb{R}_+$, $\mathbb{B}=\mathbb{R}$ and  $\mathbb{B}=\mathbb{C}$ respectively.  
  We make the convention that for any $t_1\leq t_2$,
 \beqnn
 \int_{t_1}^{t_2} = \int_{(t_1,t_2]}
 \quad \mbox{and}\quad
 \int_{t_1}^\infty = \int_{(t_1,\infty)}
 \eeqnn
 and for each $k\geq 1$, we denote by $\mathcal{I}^k_f$ the $k$-th repeated integral of $f$ with base point $0$ defined by
 \beqnn
 \mathcal{I}^k_f(t):= \int_0^t\int_0^{t_1}\cdots\int_0^{t_{k-1}}f(t_k)dt_k\cdots dt_2 dt_1, \quad t\geq 0.
 \eeqnn
 We set $\mathcal{I}_f=\mathcal{I}^1_f$ for convention. 
 For $w \in M(\mathbb{R}_+;\mathbb{B})$, we also define $\mathcal{I}_w(t):= w([0,t])$ for any $t\geq 0$. 
 We denote by $\hat{f}$ the Laplace-Stieltjes transform of function  $f \in L^1_{\rm loc}(\mathbb{R}_+;\mathbb{R})$ 
 \beqnn
 \hat{f}(\lambda):= \int_0^\infty \lambda e^{-\lambda t}f(t)dt,\quad \lambda >0.
 \eeqnn
  For any $f,g \in \big( L^1_{\rm loc}(\mathbb{R}_+;\mathbb{R}) \big)^n$ and $\nu:=(\nu_1,\cdots,\nu_n)\in \big(M(\mathbb{R}_+;\mathbb{B})\big)^n$, we denote by $f*g$ the convolution of $f$ and $g$ and by $f*d\nu$ the convolution of $f$ and $\nu$, i.e.
 \beqnn
 f*g(t):=\sum_{i=1}^n\int_0^t f_i(t-s)g_i(s)ds  
 \quad \mbox{and}\quad 
 f*d\nu(t):= \sum_{i=1}^n\int_0^t f_i (t-s)\nu_i(ds),\quad t\geq 0.
 \eeqnn
Throughout this paper, we assume the generic constant $C$ may vary from line to line.

   \section{Preliminaries and main results}\label{Sec.MR}
 \setcounter{equation}{0}
 
 In this section, we first state key properties of Hawkes processes with general kernels including functional Fourier-Laplace transforms and a series asymptotic results on their intensities and resolvents that will be important to our analysis. 
 These results are required to answer our questions \ref{Q1}-\ref{Q3}.
 Subsequently, we state FLTs for Hawkes process in various settings. In particular, we show that classical limit theorems fail to hold for critical Hawkes processes unless the arrivals of child events is sufficiently dispersed. 
 Our limit theorems provide full answers to the above questions. 
 
 We assume throughout that our Hawkes process $N$ with kernel $\phi$ and immigration density $\mu$ is defined on a complete probability space $(\Omega,\mathscr{F},\mathbf{P})$ endowed with a filtration $\{\mathscr{F}_t:t\geq 0\}$ that satisfies the usual hypotheses. Moreover, we extend the kernel $\phi$ to the whole real line $\mathbb{R}$ by setting  $\phi(t)=0$ if $t<0$.  

 \subsection{Hawkes processes and their Fourier-Laplace transforms}
  
 Associated to the kernel $\phi$ is the {\it resolvent} $R$ that is defined by the unique solution to the linear Volterra equation 
  \beqlb\label{Resolvent}
 R= \phi+ R*\phi.
 \eeqlb
This equaiton is also known as \textit{resolvent equation} or \textit{renewal equation}. The resolvent describes the combined direct and indirect impact of an individual event on the arrival of future events. It admits the Neumann series expansion 
 \beqnn
 R(t)=\phi(t)+ \sum_{k=2}^\infty \phi^{(*k)}(t),\quad t\geq 0,
 \eeqnn
where $\phi^{(*k)}$ denotes the $k$-th convolution of $\phi$. The first term $\phi$  represents the direct impact of an immigrant event, and the term $\phi^{(*k)}$ can be interpreted as describing the indirect impact on future arrivals through the direct impact of the descendant events in the $(k-1)$-th generation. We shall also need the following three important functions:
 \beqnn
  \mathcal{I}_\Lambda(t):= \int_0^t \Lambda(s)ds, \quad 
  \mathcal{I}_R(t):= \int_0^t R(s)ds
 \quad\mbox{and}\quad
 \mathcal{I}_R^2(t):= \int_0^t \mathcal{I}_R(s) ds, \quad t\geq 0.
 \eeqnn
 The first function describes the {\it cumulative event rate}. 
 The second and third function specify the {\it cumulative impact} and {\it integrated cumulative impact} of each immigrant event and its descendants on the arrivals of future events, respectively. Let 
 \[
 	m:= \|\phi\|_{L^1}
 \]
 denote the average number of descendants of each event. 
 The \textit{total impact} $ \mathcal{I}_R(\infty)$ of an immigrant event is finite if and only if the Hawkes process is subcritical $(m<1)$, in which case $$\mathcal{I}_R(\infty)= \frac{m}{1-m}.$$ 
 
 \begin{remark}
 For supercritical processes $(m > 1)$ both $ \mathcal{I}_R$ and $\mathcal{I}_\Lambda$ grow exponentially in time and each event has an infinite number descendants with positive probability. In the critical case $(m=1)$ the long-term behavior of  $ \mathcal{I}_R$ strongly depends on the choice of $\phi$. The function often grows to infinity at a power rate, but sometimes it increases only slowly; see Proposition~\ref{Thm.AsymR} for details.
 \end{remark}

 The Hawkes process $N$ has the compensator $\mathcal{I}_\Lambda$ and the compensated point process 
 \begin{equation}
 	\widetilde{N}:=N-\mathcal{I}_\Lambda 
\end{equation}	
is an $(\mathscr{F}_t)$-martingale. By using the variation of constants formula for linear Volterra integral equations as given in \cite[p.36, Equation (1.2)]{GripenbergLondenStaffans1990}, Bacry et al. \cite{BacryDelattreHoffmannMuzy2013} established the following martingale representation theorem for the intensity process in terms of the compensated point process.

 \begin{lemma}[Martingale representation]\label{MartRep}
 	The intensity process $\Lambda$ admits the representation
 \beqlb\label{SVR}
 \Lambda(t) \ar=\ar \mu(t) + \mu*R(t) + \int_0^t R(t-s)\widetilde{N}(ds) ,\quad t\geq 0.
  \eeqlb 
 Moreover, the expected intensity and the expected number of events are given by, respectively, 
 \[
 	\mathbf{E}[\Lambda(t)]= H_\mu(t) \quad \mbox{and} \quad 
	\mathbf{E}[N(t)]= \mathbf{E}[\mathcal{I}_\Lambda(t)]= \mathcal{I}_\mu(t) + \mu *\mathcal{I}_R(t), 
\]
where 
\[
	H_{\mu}(t) := \mu(t) + \mu * R(t)  \quad \mbox{and} \quad 
	 \mathcal{I}_\mu(t) := \int_0^t \mu(s) ds.
\]
 \end{lemma}

The first two terms on the right side of (\ref{SVR}) specify the average direct and indirect impact of external events on the arrivals of future events. 
 The third term specifies the random variations in the number of child events on the arrivals of future events. Integrating both sides of (\ref{SVR}) and then using the stochastic Fubini theorem we see that the integrated intensity can be represented as  
 \beqlb\label{IntegralLambda}
 \mathcal{I}_\Lambda(t)= \mathcal{I}_\mu(t) + \mu* \mathcal{I}_{R}(t) +   \int_0^{t}   R(t-s) \widetilde{N}(s)ds,\quad t\geq 0.
 \eeqlb
 
 \begin{remark}\label{rem:mu}
For the special case $\mu(t) \equiv \mu_0>0$, it follows from the above lemma that
\[
	\mathbf{E}[N(t)] = \mu_0 \cdot t + \mu_0 \cdot \mathcal{I}^2_{R}(t).
\]
This shows that the long-run dynamics of $N$ strongly depends on the asymptotics of the function $\mathcal{I}^2_{R}$. We study its asymptotics in the next section. 
 \end{remark}
 
Based on the martingale representation theorem we can now state our first main result. It extends Theorem~4.3 in \cite{JaberLarssonPulido2019}  and Theorem~1.6 in \cite{Xu2021b} by providing an explicit exponential-affine representation of the Fourier-Laplace transform of Hawkes processes with general kernels in terms of the unique solution to the equivalent nonlinear Volterra integral equations\footnote{The second equation above can be obtained from the first using (1.1)-(1.2) in \cite[p.35-36]{GripenbergLondenStaffans1990}. To obtain the first equation from the second we rewrite the second equation as $ V = \phi*\nu+ \phi*(e^{V+f}-1-V)+\phi*V$ and then apply (\ref{Resolvent}) and Lemma~3 in \cite{BacryDelattreHoffmannMuzy2013} with $h=\phi*\nu+ \phi*(e^{V+f}-1-V)$. }
 	\beqlb\label{VolRiccati}
 	V = R*d\nu + R*( e^{V+ f}-1-V)
 	\quad \mbox{and}\quad
 	V = \phi*d\nu + \phi*(e^{V+ f}-1).
 	\eeqlb
%

 A  pair $(V,T_\infty)$ with $T_\infty\in (0,\infty]$ and $V\in L^\infty_{\rm loc}([0,T_\infty),\mathbb{C})$ is called a \textit{non-continuable solution} of (\ref{VolRiccati}) if $V$
 satisfies (\ref{VolRiccati}) on $[0,T_\infty)$ and $\|V\|_{L^\infty_{T_\infty}}=\infty$ if $T_\infty<\infty$; we say that $V$ is a {\it global solution} if $T_\infty=\infty$.  The existence of a unique non-continuable solution can be obtained  by using Banach's fixed point theorem; see e.g.  \cite[p.341-343]{GripenbergLondenStaffans1990}.  We establish the existence of a unique global solution. To state our first main result, we set 
 $$
 f*dN(t):=\int_0^t f(t-s)N(ds)\quad \mbox{and}\quad f*d\widetilde{N}(t):=\int_0^t f(t-s)d\widetilde{N}(s),\quad t\geq 0.
 $$

 \begin{theorem}[Fourier-Laplace functionals]\label{Thm201}
 	For any $\nu \in M(\mathbb{R}_+;\mathbb{C}_-)$ and any $f\in L^\infty_{\rm loc}(\mathbb{R}_+;\mathbb{C}_-) $ the equivalent nonlinear Volterra integral equations \eqref{VolRiccati} admit a unique global solution
\[
	V\in L^\infty_{\rm loc}(\mathbb{R}_+;\mathbb{C}_-),
\]
	and the Fourier-Laplace functional of the process $(N,\Lambda)$ admits the representation 
 	\beqlb\label{FourLapFunHP}
 	\mathbf{E}\big[\exp\big\{ \Lambda*d\nu(T)+ f*dN(T)\big\}\big]= \exp \big\{H_\mu*d\nu(T)+ H_\mu * W(T)  \big\},  \quad T\geq 0,
 	\eeqlb
 	where $$W:= e^{V+ f}-1-V.$$ 
 	Moreover, ${\tt Re} (\phi*d\nu) -2\mathcal{I}_\phi \leq  {\tt Re} V \leq {\tt Re}(\phi*d\nu) \leq0$ and $  |{\tt Im} V| \leq |{\tt Im}(\phi*d\nu)|+\mathcal{I}_\phi.$
 	\end{theorem}

The proof of the above theorem is given in Section \ref{Sec.EleProp}. By choosing an appropriate measure $\nu \in M(\mathbb{R}_+;\mathbb{C}_-)$ and a suitable function $f\in L^\infty_{\rm loc}(\mathbb{R}_+;\mathbb{C}_-) $ in $(\ref{FourLapFunHP})$ we obtain an explicit representation of the Fourier-Laplace functional of the triple $$(\Lambda, N,\widetilde{N})$$ as well as the characteristic function of their finite-dimensional distributions. For instance, for $g_1,g_2\in L^\infty_{\rm loc}(\mathbb{R}_+;\mathbb{C}_-)$ and $g_3\in L^1_{\rm loc}(\mathbb{R}_+;\mathtt{i}\mathbb{R}) $, the Fourier-Laplace functional
 \beqnn
 \mathbf{E}\big[\exp\big\{ g_1*\Lambda(T)+  g_2*dN(T) +g_3*d\widetilde{N}(T)\big\}\big]
 \eeqnn
 can be expressed in terms of (\ref{FourLapFunHP}) by choosing $\nu(dt):=(g_1(t)- g_3(t))dt$ and $f=g_2+g_3$. Moreover, for $ t_1,\cdots,t_k \in[0,T]$, $ z_1,\cdots,z_k \in \mathbb{C}_-$ and $u_1,\cdots,u_k \in \mathbb{C}_-$, the characteristic function
 \beqnn
 \mathbf{E}\Big[\exp\Big\{ \sum_{i=1}^k\big(z_i\Lambda(t_i)+u_iN(t_i) \big)\Big\}\Big]
  \eeqnn
 can be obtained from $(\ref{FourLapFunHP})$ by choosing $\nu(dt)= \sum_{i=1}^k z_i \cdot \delta_{T-t_i}(dt)$ and $f(t)= \sum_{i=1}^k u_i \cdot \mathbf{1}_{\{T-t\leq  t_i\}}$, where $\delta_x(dt)$ is the Dirac measure at point $x\in\mathbb{R}$. 
 For $k=1$, an alternat ive and equivalent representation of $\mathbf{E}[\exp \{ uN(t)  \}]$ has been established in  \cite{GaoZhu2018b,HawkesOakes1974} by using the probability generating function of the corresponding birth-immigration process.

 \subsection{Asymptotic results for resolvents}

 We proceed to establish various asymptotic results for integrated resolvents that will be key to our analysis of the long-run behavior of Hawkes processes. As we will see the long-run behavior is fully determined by the following quantities:
 \beqnn
 {\it\Phi}(t):=\int_t^\infty \phi(s)ds
 \quad \mbox{and}\quad
  {\it \Psi}_k(t):= \int_0^t s^k \phi(s)ds, \quad
  t\geq 0,\ k\in\mathbb{Z}_+.
 \eeqnn
 Since $m={\it \Phi}(0)$ specifies the total impact of an individual event on the arrival of future events measured by its expected number of descendants, the function $ {\it\Phi} $ can be interpreted as describing the {\it temporal release} of the total impact, and the function $ {\it \Psi}_k$ can be interpreted as a measure for the {\it dispersion of child events}. 
 These functions measure how close child events arrive to their respective mother events. The not necessarily finite quantity 
 \beqnn
 \sigma :=  {\it \Psi}_1(\infty) \in(0,\infty]
 \eeqnn
 will be of particular relevance to us. As we will see, the long-run behavior of critical Hawkes processes strongly depends on whether this quantity is finite, or not. Using integration by parts, we see that the Laplace-Stieltjes transforms of ${\it \Phi}$ and ${\it\Psi}_k$ are given by
 \beqlb
 {\it\hat\Phi}(\lambda)
 \ar:=\ar  \int_0^\infty \lambda e^{-\lambda t} {\it\Phi}(t)dt = \int_0^\infty \big(1-e^{-\lambda t}\big)\phi(t)dt,\label{LapPhi}\\
 {\it\hat\Psi}_k(\lambda)
 \ar:=\ar \int_0^\infty \lambda e^{-\lambda t} {\it\Psi}_k(t)dt = \int_0^\infty e^{-\lambda t}t^k\phi(t)dt.\label{LapPsi}
 \eeqlb
 Using integration by parts again and then taking Laplace-Stieltjes transforms on both sides of (\ref{Resolvent}), the Laplace-Stieltjes transforms of the integral processes $\mathcal{I}_R$ and $\mathcal{I}_R^2$ are given by
 \beqlb\label{eqn.LapIR}
 \hat{\mathcal{I}}_R(\lambda)= \frac{m - {\it\hat\Phi}(\lambda)}{1-m+{\it\hat\Phi}(\lambda)} 
 \quad \mbox{and}\quad
 \hat{\mathcal{I}}_R^2(\lambda)= \hat{\mathcal{I}}_R(\lambda) /\lambda ,\quad \lambda >0.
 \eeqlb
 
Before stating our asymptotic results on the functions  $\mathcal{I}_R$ and $\mathcal{I}_R^2$, we recall some elements of regular variation, which is the appropriate mathematical tool for describing heavy-tail phenomena, long-range dependencies and limit theorems.\footnote{We refer to \cite{BinghamGoldieTeugels1987,DeHaanFerreira2006} for a general theory of regular variation.}
 For $\alpha\in\mathbb{R}$, we denote by ${\rm RV}_\alpha^\infty$ the space of all functions on $\mathbb{R}$ that are \textit{regularly varying at infinity} with index $\alpha$, i.e., the class of all functions $f$ that are eventually positive or negative and that satisfy
 \beqlb\label{Def.RV}
  \lim_{t\to\infty} \frac{f(tx)}{f(t)} = x^\alpha, \quad x>0.
 \eeqlb
The functions in ${\rm RV}_0^\infty$ are also said to be \textit{slowly varying}. For each $f\in {\rm RV}_\alpha^\infty$, we have that $t^{-\alpha}f(t)\in {\rm RV}_0^\infty$. 
 Applications of Karamata's theorem and the Tauberian theorem given in  Appendix~\ref{Appendix.RV} immediately yield the equivalence between the regular variation of ${\it \Phi}$, ${\it \Psi}_k$ and their Laplace-Stieltjes transforms. 

  \begin{proposition}\label{LapTheta}
 For $\alpha\in(0,1)$ and $k\geq 1$, the following four statements are equivalent:
 \medskip \smallskip \\ \medskip \smallskip
  \centerline{ {\rm(1)} ${\it \Phi} \in{\rm RV}_{-\alpha}^\infty$;\quad  {\rm(2)} ${\it \Psi}_k\in{\rm RV}_{k-\alpha}^\infty$; \quad {\rm(3)}  ${\it\hat\Phi}(1/\cdot) \in {\rm RV}_{-\alpha}^\infty$;\quad  {\rm(4)} ${\it\hat\Psi}_k(1/\cdot) \in {\rm RV}_{k-\alpha}^\infty$.}
 \noindent In either case, as $t\to\infty$,
 \beqnn
 {\it\hat \Psi}_k(1/t)\sim \Gamma(1+k-\alpha) {\it \Psi}_k(t)  
 \quad \mbox{and}\quad 
 {\it\hat \Phi}(1/t)\sim 
 \left\{
 \begin{array}{ll}
 	 \Gamma(1-\alpha) {\it \Phi}(t),  & \mbox{if }\alpha\in[0,1);\\
 	 {\it\Psi}_1(t), & \mbox{if } \alpha=1.
 	\end{array}
 \right.
 \eeqnn 
 \end{proposition}

 We are now ready to specify the tail behavior of (iterated) integrals of the resolvent $R$. The proof of the following result is given in Section \ref{Sec.EleProp}.

  \begin{proposition}\label{Thm.AsymR}
 Three regimes need to be distinguished when analyzing the long-term behavior of the integral functions $\mathcal{I}_R$ and $\mathcal{I}_R^2$.
   \begin{enumerate}
   \item[(1)] When $m<1$, we have as $t\to\infty$, 
    \beqnn
    \mathcal{I}_R(t)\to \frac{m}{1-m}  
    \quad \mbox{and}\quad 
    \mathcal{I}_R^2(t)\sim \frac{m\cdot t}{1-m}.
    \eeqnn

   \item[(2)] When $m=1$ and $\sigma <\infty$,  we have as $t\to\infty$, 
    \beqnn
    \mathcal{I}_R(t)\to \frac{t}{\sigma}  
    \quad \mbox{and}\quad 
    \mathcal{I}_R^2(t)\sim \frac{t^2}{2\sigma}.
    \eeqnn
  		
   \item[(3)] When $m=1$ and $\sigma =\infty$, if ${\it \Phi}\in{\rm RV}_{-\alpha}^\infty$ for some $\alpha\in[0,1]$,  we have as $t\to\infty$, 
    \beqnn
 	 \mathcal{I}_R(t)   \sim  
 	 \begin{cases}
      \displaystyle{ \frac{1/{\it\Phi}(t)}{\Gamma(1-\alpha)\Gamma(1+\alpha)}  \in {\rm RV}_{\alpha }^\infty,} & \mbox{if }\alpha\in[0,1);\vspace{5pt} \\
 	  \displaystyle{\frac{t}{{\it\Psi}_1(t)}   \in {\rm RV}_{ 1}^\infty,} & \mbox{if }\alpha=1,
     \end{cases} 
 	 \quad \mbox{and}\quad
 	 \mathcal{I}_R^2(t) 
 	 \sim \frac{t\cdot 	\mathcal{I}_R(t)}{1+\alpha}   \in {\rm RV}_{\alpha +1}^\infty.
 	\eeqnn 
	In particular, 
	 \beqlb\label{eqn.3.21}
		 \frac{\mathcal{I}^2_R(n)}{n} \sim \frac{\mathcal{I}_R(n)}{1+\alpha} \to \infty \quad \mbox{as} \quad n \to \infty.
 	\eeqlb

 	\end{enumerate} 
 \end{proposition}

 For weakly critical Hawkes processes, integration by parts along with an application of Proposition~\ref{Thm.AsymR}(2) shows that 
 \beqnn
 \int_0^t R(ns)ds =\frac{\mathcal{I}_{R}(nt)}{n} \to \frac{t}{\sigma},\quad t\geq 0,
 \eeqnn
 as $n\to\infty$. This yields the weak convergence of  the $\sigma$-finite measure with density function $$R^{(n)}:=\{R^{(n)}(t)=R(nt):t\geq 0\}$$ to the $\sigma$-finite measure with density $R^*:=\{R^*(t)=1/\sigma:t\geq 0\}$ on any compact set.  In particular, we have the following result.

   \begin{corollary}\label{ConvergenceR}
  	When $m=1$ and $\sigma<\infty$,  we have $\mathcal{I}_{R^{(n)}}\to \mathcal{I}_{R^*}$ locally uniformly as $n\to\infty$.
  \end{corollary}

The above corollary establishes the convergence of the rescaled integrated resolvent function.  The speed of convergence will be studied detail below; as we will see it fully specifies the rate of convergence of rescaled weakly critical Hawkes processes to their limits.

 \subsection{Main results}

 In this section we state our functional limit theorems (FLTs) for subcritical and critical Hawkes processes with constant exogenous density, i.e., we assume that
 \beqnn
 \mu(t)\equiv \mu_0 \in (0,\infty),\quad t\geq 0.
 \eeqnn 
 Analogous results also can be established for time-inhomogeneous or nondeterministic exogenous intensities $\mu$ under mild conditions on $\mu+R*\mu$. For instance, if $\mu(t)=\mu_0 + \hat{\mu}(t)$ with $\hat{\mu} \in L^1(\mathbb{R}_+;\mathbb{R}_+)$ representing the total impacts of all past events prior to time $0$ on the arrivals of future events, then all results in the sequel remain valid if $\hat\mu+R*\hat\mu$ is bounded on $\mathbb{R}_+$. 
 
 Under the assumption of a constant immigration rate it follows from Remark \ref{rem:mu} that the expected value of the rescaled point process satisfies 
 \[
 	\frac{\mathbf{E}[N(nt)]}{n} = \mu_0 t + \mu_0 \frac{\mathcal{I}^2_R(nt)}{n}.
 \]
 In what follows, we  denote by $B:=\{B(t):t\geq 0\}$ a one-dimensional standard Brownian motion.

   \subsubsection{Subcritical processes: $m<1$} \label{ThmStationaryHP}
 
 In this section we answer question \ref{Q1}. 
 We start by recalling the following functional limit theorems for subcritical Hawkes processes, due to Bacry et al.~\cite{BacryDelattreHoffmannMuzy2013}.  
 
 \begin{theorem}[Bacry et al. \cite{BacryDelattreHoffmannMuzy2013}]\label{MainThm.01}
 If $m<1$, then the following holds as $n\to\infty$. 
 \begin{enumerate}
   \item[(1)] $ \sup_{t\in[0,T]}\big| \frac{N(nt)}n- \frac{t \mu_0}{ 1-m} \big|   \to 0$ in a.s.~and in $L^2 (\mathbf{P})$, for any $T\geq0$;

   \item[(2)] $  n^{-1/2} \big( N(nt)- \mathbf{E}[N(nt)] \big) \to \sqrt{\mu_0
       (1-m)^{-3}} \cdot B(t)$ weakly in $\mathbf{D}(\mathbb{R}_+; \mathbb{R})$.
 \end{enumerate}
 \end{theorem}
 
 The FCLT uncovers that even if the kernel $\phi$ is heavy-tailed and no matter how dispersed child events are, the increments of subcritical Hawkes processes are always asymptotically independent and long-range dependencies can never be observed in the subcritical case. 
 This calls for the use of critical Hawkes processes to capture  the phenomena of long-range dependence.  
 We study critical Hawkes processes in the next subsection. 

In the subcritical case it follows from Proposition \ref{LapTheta} that $\frac{\mathbf{E}[N(nt)]}{n} \sim  \frac{t \mu_0}{1-m}$. 
This suggests that Theorem~\ref{MainThm.01}(2) also holds if $\mathbf{E}[N(nt)]$ is replaced by $\frac{t \mu_0}{1-m}$. Bacry et al. \cite{BacryDelattreHoffmannMuzy2013} proved that this is  in fact true, if
 \beqlb\label{BacryCon}
 \int_0^\infty \sqrt{t} \cdot \phi(t)dt<\infty.
 \eeqlb
The next theorem, whose proof is given in Section \ref{Sec.SubHP}, generalizes their result under a much weaker condition and shows that the limiting Brownian motion may have to be replaced by a Brownian motion with  polynomial drift if the condition (\ref{BacryCon}) fails.  Moreover, the FCLT may degenerate if the tail of the kernel $\phi$ is too heavy and child events are too dispersed.

 \begin{theorem} \label{MainThm.02}
If $m<1$ and $ {\it \Psi}_1(t)/\sqrt{t}\to {\it \Psi}_*\in[0,\infty]$ as $t\to\infty$, then
the following holds as $n\to\infty$:
 \begin{enumerate}
 \item[(1)] if ${\it \Psi}_*<\infty $, then weakly in $\mathbf{D}([0,\infty);\mathbb{R})$,
   \beqnn
   \sqrt{n} \Big(\frac{N(nt)}{n}-\frac{\mu_0\cdot t}{1-m} \Big) \to \sqrt{\mu_0(1-m)^{-3}}\cdot  B(t)  - \frac{\mu_0 {\it \Psi}_*}{ (1-m)^2}\cdot \sqrt{t} ;
   \eeqnn
 	
 \item[(2)] if ${\it \Psi}_*=\infty$ and ${\it \Phi}\in {\rm RV}_{-\alpha}^\infty$ with $\alpha\in[0, 1/2]$, then for any $T>0$,
   \beqnn
   \sup_{t\in[0,T]}\Big| \frac{1}{{\it\Phi}(n)} \Big( \frac{N(nt)}{n} -\frac{\mu_0\cdot t}{1-m}  \Big)+\frac{\mu_0\cdot t^{1-\alpha}}{(1-\alpha)(1-m)^{2}} \Big|\to 0,\quad \mbox{in }\mathbf{P}.
 \eeqnn		
 \end{enumerate}
 \end{theorem}

Under  the hypothesis (\ref{BacryCon}), it always holds that ${\it \Psi}_* =0$. In fact for each $K>0$ and $t>K$,
  \beqnn
 \frac{{\it \Psi}_1(t)}{\sqrt{t}} = t^{-1/2}\int_0^K s \phi(s)ds+ \int_K^t \frac{s}{t^{1/2}} \phi(s)ds
  \leq  t^{-1/2}\int_0^K s \phi(s)ds + \int_K^\infty \sqrt{s} \phi(s)ds ,
  \eeqnn
  which tends to $0$ as $t\to\infty$ and then $K\to\infty$. The converse implication, however, does not hold in general. In particular, the conditions are not equivalent; our condition is in fact much weaker. To see this, let  
  \[
  	\phi(s)= (1+s)^{-3/2}\cdot(2+\log(1+s))^{-1}. 
\]	
Then, $\|\phi\|_{L^1}<1$ so that the process is indeed subcritical and condition (\ref{BacryCon}) fails because
  \beqnn
  \int_0^\infty \sqrt{s}\phi(s)ds \geq   \int_1^\infty \frac{s^{-1}}{2+\log s}  ds=\infty. 
  \eeqnn
At the same time, for any $\epsilon > 0$ by first choosing $K$ and then $t$ large enough,  
 \beqnn
 \frac{{\it \Psi}_1(t)}{\sqrt{t}} \leq \int_0^K \frac{ t^{-1/2}(1+s)^{-1/2}}{2+\log(1+s)}ds  + \int_K^t\frac{t^{-1/2}  (1+s)^{-1/2} }{2+\log K} ds
  \leq \sqrt{\frac{K}{t}} + \frac{2}{2+\log K} < \epsilon.
  \eeqnn

\medskip

%

   \subsubsection{Weakly critical processes:  $m=1$ and $\sigma<\infty$}
 \label{ThmWeaklyHP}
 
 We proceed to weakly critical Hawkes processes, providing partial answers to questions \ref{Q2} and \ref{Q3}. 
 Our goal is to prove a non-standard convergence result for rescaled weakly critical Hawkes processes and to provide explicit bounds for the rate of convergence. 

  \smallskip
 
{\bf Scaling limits.} For weakly critical Hawkes processes neither a standard FLLN nor a standard FCLT can be expected. This can already bee seen for the benchmark case of weakly critical processes with exponential kernel $\phi(t):=\beta e^{-\beta t}$, for some $\beta >0$. 

In this case the sequence of rescaled intensity process $ \{ \Lambda(n \cdot)/n :n \in \mathbb{N}\}$ satisfies Condition 3.4 in \cite{HorstXu2022}. By Theorem 3.12 in \cite{HorstXu2022} the rescaled intensity process converges weakly to a critical branching diffusion with immigration $\Lambda^*$ and the point process $N$ behaves as the integrated branching process process as for large $n \in \mathbb N$ when suitably rescaled. Specifically, 
 \beqnn
 \Lambda^{(n)}(t) \sim \Lambda^* (t) 
 \quad \mbox{and}\quad 
\frac{ N(nt)}{n^2} \sim  \mathcal{I}_{\Lambda^*}(t), \quad t\geq 0.  
 \eeqnn

We prove that the above result holds for any kernel, light-tailed or not, as long as the child events are not too dispersed. To this end, we recall that $\widetilde N = N - \mathcal{I}_\Lambda$ denotes the compensated Hawkes process and introduce the family of rescaled processes $$X^{(n)}:=\left( \Lambda^{(n)}, \mathcal{I}^{(n)}_\Lambda, N^{(n)},
 \widetilde{N}^{(n)} \right)$$ where
 \begin{equation} \label{tildeN}
  \Lambda^{(n)}(t):= \frac{\Lambda(nt)}{n},\quad \mathcal{I}^{(n)}_\Lambda(t):=
  \frac{\mathcal{I}_{\Lambda}(nt)}{n^2},\quad N^{(n)}(t):= \frac{N(nt)}{n^2}
  \quad \mbox{and}\quad
  \widetilde{N}^{(n)}(t):= \frac{\widetilde{N}(nt)}{n},\quad t\geq 0.
 \end{equation}

The assumption that $\phi\in L^1 (\mathbb{R}_+,\mathbb{R}_+)$ implies that $X^{(n)} \in L^1_{\rm loc} (\mathbb{R}_+; \mathbb{R}_+)\times\mathbf{D}(\mathbb{R}_+;\mathbb{R})^3$.
Due to the lack of regularity of the kernel, we consider the rescaled intensity $\Lambda^{(n)}$ as an $ L^1_{\rm loc} (\mathbb{R}_+; \mathbb{R}_+)$-valued random variable. 
Since the space $L^1_{\rm loc} (\mathbb{R}_+; \mathbb{R}_+)$ is not relatively compact, establishing the tightness of the sequence $\{\Lambda^{(n)} : n \geq 1\}$ turns to be challenging. To bypass this difficulty, we consider the convergence in the $\sigma$-finite measure space $M(\mathbb{R}_+;\mathbb{R}_+)$. 

We say a sequence of $L^1_{\rm loc}(\mathbb{R}_+;\mathbb{R}_+)$-valued random variables $\{Y_n : n \geq 1\}$ converges to $Y_*\in L^1_{\rm loc} (\mathbb{R}_+;\mathbb{R}_+)$ weakly in $M(\mathbb{R}_+;\mathbb{R}_+)$ if the  $\sigma$-finite random measure with density function $Y_n$ converges weakly to the $\sigma$-finite random measure with density function $Y_*$ on any compact set\footnote{The convergence is equivalent to the weak convergence of $\mathcal{I}_{Y_n}$ to $\mathcal{I}_{Y_*}$ in $\mathbf{C}(\mathbb{R}_+;\mathbb{R}_+)$. Indeed, for each $g\in L^\infty_{\rm loc}(\mathbb{R}_+;\mathbb{R}_+)$ and $t\geq 0$, we have $g* \mathcal{I}_{Y_n}(t)= \mathcal{I}_g*Y_n(t)$ and $g*\mathcal{I}_{Y_n}(t)= \mathcal{I}_g*Y_*(t)$. Hence $Y_n\to Y_*$ weakly in $M(\mathbb{R}_+;\mathbb{R}_+)$ if and only if $\mathbf{E}[\exp\{g*\mathcal{I}_{Y_n}(t)\}]\to \mathbf{E}[\exp\{g*\mathcal{I}_{Y_*}(t)\}]$, which is equivalent to $\mathcal{I}_{Y_n}\to \mathcal{I}_{Y_*}$ weakly in $\mathbf{C}(\mathbb{R}_+;\mathbb{R}_+)$; see Theorem~3.37 in \cite[p.354]{JacodShiryaev2003}. \label{Footnote.03}}.
  
  We are now ready to state our second main result. The proof of the following two theorems is given in Section \ref{sec:MainThm04}. 
  
  \begin{theorem}\label{MainThm.04}
  If $m =1$ and $\sigma<\infty$,  then $X^{(n)}\to X^*$  weakly in
  $M(\mathbb{R}_+;\mathbb{R}_+)\times\mathbf{D}(\mathbb{R}_+;\mathbb{R})^3 $ as
  $n\to\infty$, where 
  \[
  	X^*:=\left\{\big(\Lambda^*(t), \mathcal{I}_{ \Lambda^*}(t), \mathcal{I}_{\Lambda^*}(t), \sigma\Lambda^*(t)-\mu_0\cdot t\big):t\geq 0 \right\}
\]	
and $\Lambda^*$ is the unique strong solution to the SDE
   \beqlb\label{CBI}
   \Lambda^*(t)= \frac{\mu_0}{\sigma} \cdot  t + \int_0^t \frac{1}{\sigma}\sqrt{
   \Lambda^*(s)}d B(s).
   \eeqlb
  \end{theorem}

 The limit process $\Lambda^*$ is a critical continuous-state branching process with a deterministic immigration rate, which has been extensively studied in the literature; we refer to \cite{DuffieFilipovicSchachermayer2003,Li2019}  for reviews. The Laplace transform of the process $ \big\{\int_0^t \Lambda^*(s)w(ds):t\geq 0 \big\}$ with $w\in M(\mathbb{R}_+;\mathbb{R}_-)$ can be found in Theorem~4.2 in \cite{Li2019}. The next theorem extends these results; it  provides an explicit exponential-affine representation of the Fourier-Laplace functional of the vector $X^*$ in term of solutions to Riccati equations. 

 \begin{theorem}\label{Thm.FourLapCIR}
  For each $w\in M(\mathbb{R}_+;\mathbb{C}_-)$ and $g\in L^\infty_{\rm
  loc}(\mathbb{R}_+;\mathtt{i}\mathbb{R})$, we have that
  \beqlb \label{FourLapCIR}
 \mathbf{E}\Big[\exp\Big\{  \Lambda^* * dw(T) + \int_0^T g(T-s)\sqrt{ \Lambda^*(s)}d
 B(s)\Big\}\Big]= \exp\big\{ \mu_0 \cdot \mathcal{I}_{V^*}(T)  \big\}, \quad T\geq 0,
 \eeqlb
 where $ V^*\in L^\infty_{\rm loc}(\mathbb{R}_+;\mathbb{C}_-)$ is the unique solution to
 the Riccati equation\,\footnote{ For fixed $w \in M(\mathbb{R}_+;\mathbb{C}_-)$ and $g \in L^\infty_{\rm loc}(\mathbb{R}_+;\mathtt{i}\mathbb{R})$, the existence, uniqueness, regularity and local boundedness of solution $V^*$ to (\ref{CBI.Riccati}) can be founded in many textbooks on ordinary differential equations. Moreover, $V^*$ is locally of bounded variation.
 }
 \beqlb\label{CBI.Riccati}
 V^*(t) = \frac{\mathcal{I}_w (t)}\sigma     +  \frac{\mathcal{I}_{(V^*+g)^2}(t)}{2\sigma},\quad t\geq
 0.
 \eeqlb
 Here  $\mathcal{I}_w (t)= w([0,t])$. Moreover, the following holds:
 \begin{enumerate}
 	\item[(1)]  If $w\in M(\mathbb{R}_+;\mathbb{R}_-)$ and $g \equiv 0$, then $ V^*\in
 	L^\infty_{\rm loc}(\mathbb{R}_+;\mathbb{R}_-)$ and $\mathcal{I}_w/\sigma \leq V^*\leq 0$.
 	
 	\item[(2)] If $$w(dt)=f(t)dt \quad \mbox{with some density} \quad  f\in
 	L^\infty_{\rm loc}(\mathbb{R}_+;\mathbb{C}_-),$$ then the function $V^*$ is differentiable on $\mathbb{R}_+$. Moreover,
for each $T\geq0$, there exists a constant $C>0$ such that
 	for any $f\in L^\infty_{\rm loc}(\mathbb{R}_+;\mathbb{C}_-)$ and $g\in L^\infty_{\rm loc}(\mathbb{R}_+;\mathtt{i}\mathbb{R})$,
 	\beqlb\label{UpperBound}
 	\|V^*\|_{L^\infty_T}\leq C \cdot \big(\|f\|_{L^\infty_T}+\|g\|_{L^\infty_T}^2\big).
 	\eeqlb
 	\end{enumerate}

 \end{theorem}

 
 The uniform upper bound in (\ref{UpperBound}) is key to our analysis. To the best of our knowledge the bound is new and seems to be quite difficult to achieve directly from (\ref{CBI.Riccati}). 
 To obtain the bound, we first give exact representations for Fourier-Laplace functionals of the vectors $\{X^{(n)}: n \geq 1\}$ by using Theorem~\ref{Thm201} and then prove that their exponents converge locally uniformly to those of $X^*$. The desired upper bound is then obtained by establishing an analogous uniform upper bound  for the Fourier-Laplace exponents of $X^{(n)}$. 
 
 \medskip
 
 {\bf Speed of convergence.} We now analyze the speed of convergence of $X^{(n)}$ to $X^*$ by establishing upper bounds on the distance between their distributions with respect to the Fourier-Laplace functional and the Wasserstein distance $d_{\rm W}$ repectively; see Appendix~\ref{Appendix.DM}. 
To this end, we introduce for any $K\geq 0$ the function spaces
  \beqnn
 \mathcal{A}_{K} := \big\{ f \in
 C^1(\mathbb{R}_+;\mathbb{C}^4) :  \|f\|_{L^\infty}, \|f'\|_{L^\infty}\leq K \big\}
 \quad\mbox{and}\quad
 \mathcal{B}_{K} := \big\{ f\in L^\infty(\mathbb{R}_+;\mathbb{C}^4):  \|f\|_{L^\infty}\leq K \big\},
 \eeqnn
denote by $\mathcal{A}_{K}^-$ [resp.~$\mathcal{B}_{K}^-$] the subset of functions $f\in
 \mathcal{A}_{K}$ [resp.~$f\in \mathcal{B}_{K}$] with ${\tt Re}f_1,{\tt Re}f_2\leq 0$ and ${\tt Re}f_3,{\tt Re}f_4=0$, and recall the constant function $R^* \equiv 1/\sigma$ defined above Corollary~\ref{ConvergenceR}.  The following result is proved in Section \ref{sec:MainThm.07}.

   \begin{theorem}\label{MainThm.07}
 For each $T> 0$, $K\geq 0$ and $\kappa\in(0,1/2)$, there exists a constant  $C>0$ such that
 for any $n\geq 1$ the following holds:
 \beqlb
  \sup_{f\in \mathcal{A}_{K}^-}\sup_{ t\in[0,T]} \big|\mathbf{E}\big[e^{f*X^{(n)}(t)
 }\big]-\mathbf{E}\big[e^{f*X^*(t) }\big]\big|
 \ar\leq\ar C\cdot \big( n^{-1} +  \|\mathcal{I}_{R^{(n)}}- \mathcal{I}_{R^*}
 \|_{L^\infty_T} \big),\label{UpperBLapTran}\\
  \sup_{f\in \mathcal{A}_{K}}\sup_{ t\in[0,T]}d_{\rm W} \big(
 f*X^{(n)}(t),f*X^*(t)\big)
 \ar\leq\ar C\cdot \big|\log   \big( n^{-1} +  \| \mathcal{I}_{R^{(n)}}- \mathcal{I}_{R^*} \|_{L^\infty_T} \big) \big|^{-\kappa}. \label{UpperBWass}
 \eeqlb
If $\phi$ has bounded variation, then $\|R^{(n)}-R^* \|_{L^2_T}\to 0$ as $n\to\infty$
 and the preceding upper bounds also hold with $\big(\mathcal{A}_{K}, \mathcal{A}_{K}^-,
 \|\mathcal{I}_{R^{(n)}}- \mathcal{I}_{R^*} \|_{L^\infty_T}\big)$ replaced by $\big(\mathcal{B}_{K}, \mathcal{B}_{K}^-, \|R^{(n)}-R^* \|_{L^2_T}\big)$.
 \end{theorem}
 
Choosing $f=(0,0,1,0)$, respectively $f=(0,0,0,1)$, the Wasserstein distance in the inequality (\ref{UpperBWass}) reduces to 
 \beqnn
 d_{\rm W}\Big(N^{(n)}(T), \mathcal{I}_{\Lambda^*}(T)\Big), \quad \mbox{respectively,}\quad d_{\rm W}\Big(\widetilde{N}^{(n)}(T), \sigma\Lambda^*(T)-\mu_0\cdot T\Big).
 \eeqnn
Thus, Theorem~\ref{MainThm.07} establishes the convergence rates of rescaled or normalized Hawkes process to their corresponding limit processes. Moreover, the theorem reveals that 
the convergence rate of $X^{(n)}$ to $X^*$ is determined by the speed of convergence of $\mathcal{I}_{R^{(n)}}$ to  $\mathcal{I}_{R^*}$; see also Corollary~\ref{ConvergenceR}. 

\begin{example}
If the kernel $\phi$ is exponential, then  $R^{(n)}=R=R^*$ and hence $$\|\mathcal{I}_{R^{(n)}}-\mathcal{I}_{R^*} \|_{L^\infty_T} =\|\mathcal{I}_{R^{(n)}}-
 \mathcal{I}_{R^*} \|_{L^2_T}=0.$$ 
\end{example}

Beyond the benchmark case of exponential kernels closed-form representations of the quantities $\|\mathcal{I}_{R^{(n)}}-\mathcal{I}_{R^*} \|_{L^\infty_T} $ and $ \|\mathcal{I}_{R^{(n)}}- \mathcal{I}_{R^*} \|_{L^2_T}$ are usually not available. However, under additional regular variation assumptions, the next proposition establishes upper bounds in term of the integrated kernel ${\it \Phi}$. The proof is given in Section \ref{sec:MainThm.08}.  

 \begin{proposition}\label{Thm.ConvRateRLinfty}
 \begin{itemize}
\item[(i)] For each $T> 0$, there exists a constant  $C>0$ such that for any $n\geq 1$,
 \beqnn
 \displaystyle \big\|\mathcal{I}_{R^{(n)}}- \mathcal{I}_{R^*}\big\|_{L^\infty_T} \leq \left\{
 \begin{array}{ll}
 C\cdot n^{-1/3}, & \mbox{if ${\it\Psi}_2(\infty)<\infty$};\vspace{3pt}\\
 C \cdot |n{\it\Phi}(n)|^{1/3} , & \mbox{if ${\it\Phi}\in {\rm RV}_{-\alpha-1}^\infty$ with
 $\alpha \in (0,1)$}.
 \end{array}
 \right.
 \eeqnn 
%
\item[(ii)] If $\phi$ has bounded variation, for each $T> 0$ and $\ell\in(0,1/2)$, there exists a
 constant $C>0$ such that for any $n\geq 1$,
 \beqnn
 \displaystyle \big\|\mathcal{I}_{R^{(n)}}-\mathcal{I}_{R^*} \big\|_{L^2_T}\leq \left\{
 \begin{array}{ll}
 C\cdot n^{-1/2}, & \mbox{if ${\it \Psi}_2(\infty)<\infty$};\vspace{3pt} \\
 C\cdot n{\it\Phi}(n), & \mbox{if ${\it \Phi}\in {\rm RV}_{-\alpha-1}^\infty$ with
 $\alpha\in(0,1/2)$};\vspace{3pt} \\
 C\cdot n^{\alpha+\ell +1/2} {\it\Phi}(n)  , & \mbox{if ${\it \Phi}\in {\rm
 RV}_{-\alpha-1}^\infty$ with $\alpha\in[1/2,1)$}.
 \end{array}
 \right.
 \eeqnn
\end{itemize}
\end{proposition}
 

 \medskip

   \subsubsection{Strongly critical processes: $m=1$ and $\sigma=\infty$}  
 \label{ThmStonglyHP}
 
 We conclude our analysis of Hawkes processes by formulating FLTs for strongly critical Hawkes processes. 
 Along with the previously established limit results for weakly critical processes, they provide a full answer to questions \ref{Q2} and \ref{Q3}. 
 All proofs are given in Section \ref{Sec.StronglyHP}. 
 
 Contrary to the weakly critical case, strongly dispersed child events render Hawkes processes more stable. In contrast to weakly critical Hawkes processes where the excitation of each event on the future events will be fully released in a short time, for strongly critical processes the release of excitation is much slow and may last much longer. We refer to this phenomena as the loss of self-excitation. 
 
 As a result, it becomes increasingly difficult to trace a descendant event back to its mother event. This suggests that the distribution of child events and hence of all events should enjoy some mixing property in the long run in which case both a FLLN and a FCLT should prevail. The following result confirms our intuition; the proof is given in Section 6. 

 \begin{theorem}\label{MainThm.05}
 If $m=1$, $\sigma=\infty$ and ${\it \Phi}\in{\rm RV}_{-\alpha}^\infty$ with
 $\alpha\in[0,1]$, then the following holds as $n\to\infty$.
 \begin{enumerate}
  \item[(1)] $\displaystyle \sup_{t\in[0,T]}\Big| \frac{N(nt)}{\mathcal{I}^2_R(n)} -  \mu_0\cdot t^{\alpha+1}
      \Big|\to  0$ in $\mathbf{P}$ and $L^2(\mathbf{P})$ for any $T\geq 0$;

  \item[(2)] $\displaystyle n|\mathcal{I}^2_R(n)|^{-3/2}\cdot \big(  N(nt)-  \mathbf{E}[N(nt)] \big)\to \sqrt{\mu_0(\alpha+1)}\int_0^t (t-s)^\alpha s^{\alpha/2}dB(s)$ weakly in
  $\mathbf{D}(\mathbb{R}_+; \mathbb{R})$.
 \end{enumerate}
 \end{theorem}


 \begin{remark}\label{Remark.02}
 It follows from Theorem~\ref{MainThm.05}(2) that long-range dependencies in the limiting Gaussian process decrease as $\alpha$ decreases to zero; in particular, the limiting process is a standard Brownian motion when $\alpha=0$. This is very intuitive. By Proposition~\ref{LapTheta}, we have  ${\it \Psi}_1\in {\rm RV}_{1-\alpha}^\infty$, which implies that the smaller $\alpha$ is, the more dispersed and hence the more ``stationary'' child events are distributed in the long run. 
\end{remark}

For {subcritical} Hawkes processes we have shown above that the expected number of events in the FCLT can be replaced by the average long run number of events under mild conditions on the model dynamics. In what follows we establish a similar result for strongly critical processes using results from the theory second-order regular variation. 

Second-order regular variation refines the concept of first-order regular variation by specifying convergence rates in first-order regular variation. The concept of second-order variation has proven itself to be a powerful tool in extreme value theory. We refer the reader to \cite{DeHaanFerreira2006, DeHaanPeng1997, DeHaanResnick1996, DeHaanResnick1996, HorstXu2023a, Smith1982} for a detailed discussion of second-order regular variation and its applications to  extreme value theory.    

A function $f \in {\rm RV}_\alpha^\infty$ with  $\alpha\in\mathbb{R}$ is said to be of \textit{second-order regular variation at infinity} with second order parameter $\rho\leq 0$ if there is a function $A\in {\rm RV}_\rho^\infty$ vanishing at infinity such that
 \beqlb\label{2RV}
   \lim_{t\to\infty}\frac{f(tx)/f(t)-x^\alpha}{A(t)}=  x^\alpha \int_1^x u^{\rho-1}du,
   \quad x>0.
 \eeqlb
In this case we write $f\in {\rm 2RV}_{\alpha,\rho}^\infty(A)$ and call $A$ as the \textit{auxiliary function}. 
We write $f\in {\rm 2RV}_{\alpha,-\infty}^\infty(0)$ if $f(t)=C\cdot t^\alpha$ for any large $t$ and some $C \neq 0$. 
 For each $\rho\leq 0$, we denote by  $\mathscr{A}_\rho^\infty $  the space of functions  $f\in{\rm RV}_\rho^\infty $ that satisfy $f(t)\to 0$ as $t\to\infty$. 
 We also make the convention that $\mathscr{A}_{-\infty}^\infty:=\{0\}$ and $0/0=0$. 
 
 
 \begin{theorem}\label{MainThm.06}
 Assume that $m =1$, $\sigma=\infty$ and $\mathcal{I}^2_R\in {\rm 2RV}_{ \alpha+1,\rho}^\infty(A)$ for some $\alpha\in[0,1]$, $\rho\in[-\infty,0]$ and
 $A\in\mathscr{A}_\rho^\infty$. 
 If there exists a nonnegative function $\gamma$ on $\mathbb{R}_+$ such that as $n\to\infty$,
 \beqlb\label{eqn.20021}
 \gamma(n)\cdot\frac{n}{\mathcal{I}^2_R(n)} \to \gamma^*_1\in \mathbb{R}_+, \quad
 \gamma(n) A(n) \to \gamma^*_2\in \mathbb{R}
 \quad \mbox{and}\quad
 \gamma(n)\cdot\frac{|\mathcal{I}^2_R(n)|^{1/2}}{n} \to \gamma^*_3\in \mathbb{R}_+,
 \eeqlb
 with at least one of the limit coefficients being non-zero, then we have as $n\to\infty$, weakly in $\mathbf{D}(\mathbb{R}_+,\mathbb{R})$ that
 \beqnn
 \begin{split}
&  \gamma(n)\Big( \frac{N(nt)}{\mathcal{I}^2_R(n)}-  \mu_0\cdot t^{\alpha+1}   \Big) \\
& \quad \longrightarrow
 \mu_0 \Big(\gamma_1^*t + \gamma_2^* t^{\alpha+1} \int_1^t u^{\rho -1}du \Big)   +\gamma_3^* \sqrt{\mu_0(\alpha+1)} \int_0^t  (t-s)^{\alpha}s^{\alpha/2}   dB(s) .
 \end{split}
 \eeqnn
 \end{theorem}

\begin{remark}
The scaling constant in the above FCLT is highly non-standard. Under additional assumption it can be replaced by a more standard one. In fact, if $\rho <0$ then it follows from Corollary~2.4 in \cite{HorstXu2023a} that there exists a constant $C_{\mathcal{I}^2_R}>0$ such that 
\[
	\mathcal{I}^2_R(t) = C_{\mathcal{I}^2_R} \cdot t^{\alpha+1}\cdot \left(1+ \frac{A(t)}{ \rho} + o(A(t))  \right) \quad \mbox{as} \quad t\to\infty.
\] 
In this case we thus have as $n\to\infty$, weakly in $\mathbf{D}(\mathbb{R}_+,\mathbb{R})$ that
  \beqnn
  \begin{split}
&  \frac{\gamma(n)}{ C_{\mathcal{I}^2_R}}\Big( \frac{N(nt)}{n^{\alpha+1} }-  \mu_0\cdot C_{\mathcal{I}^2_R} \cdot t^{\alpha+1}   \Big) \\
& \quad  \longrightarrow  \mu_0  \Big(\gamma_1^* \cdot t + \frac{\gamma_2^* }{\rho} \cdot t^{\alpha+\rho+1}  \Big)   +\gamma_3^* \cdot  \sqrt{\mu_0(\alpha+1)} \int_0^t  (t-s)^{\alpha}s^{\alpha/2}   dB(s) . 
 \end{split}
  \eeqnn
  \end{remark}

The key assumption in establishing the above FCLT is the second-order regular variation of the iterated resolvent integral $\mathcal{I}^2_R$. A class of Hawkes processes where this assumption can be verified in analyzed in Section 2.4.1 below. Sufficient conditions in terms of the resolvent and the kernel of the Hawkes process are given {in} our accompanying work \cite{HorstXu2023a}. In particular, the following result holds; it will be applied to a class of stable Hawkes processes in Section 2.4.2 below.  


%

 \begin{proposition}[\cite{HorstXu2023a}, Lemma 3.1]
 \label{Prop.2RVR2.000}
 For $\alpha\in [0,1]$, $\rho\leq 0$ and $A\in\mathscr{A}_\rho^\infty$, 
 if $\alpha+\rho>0$ and $$R \in {\rm 2RV}^\infty_{\alpha-1,\rho}(\frac{(\alpha+\rho)(\alpha+\rho+1)}{\alpha(\alpha+1)}\cdot A)$$ or if $\alpha+\rho +1>0$ and $\mathcal{I}_R \in {\rm
 2RV}^\infty_{\alpha,\rho}(\frac{\alpha+\rho+1}{\alpha+1}\cdot A)$,  then $$\mathcal{I}^2_R \in {\rm 2RV}^\infty_{\alpha+1,\rho}(A).$$
 \end{proposition}

    \subsection{Examples}
   
 The preceding FLTs for subcritical and weakly critical Hawkes processes were established under some conditions on the function ${\it\Phi}$ that can be easily verified for a large class of examples. By contrary, the conditions in Theorem~\ref{MainThm.06} are quite complex. In this subsection, we provide several specific examples that help us understanding the nature of the conditions. 

 \subsubsection{Mixed Mittag-Leffler type}

 For $\alpha\in(0,1)$ and $\kappa>0$, let us denote by $E_{\alpha,\kappa} $ the \textit{Mittag-Leffler function} defined by
 \beqlb\label{MLF}
 E_{\alpha,\kappa}(t) := \sum_{k=0}^\infty \frac{t^k}{\Gamma(\kappa+k\alpha)},\quad t\in\mathbb{R}.
 \eeqlb
 It is customary to write $E_{\alpha}$ for $E_{\alpha,1}$.
 For a constant $\beta>0$, we denote by $F^{\alpha,\beta}$ and $f^{\alpha,\beta}$ the \textit{Mittag-Leffler distribution} and \textit{density function}  on $\mathbb{R}_+$; they are given by
 \beqnn
 F^{\alpha,\beta}(t):= 1- E_\alpha(-\beta t^\alpha)
\quad \mbox{and}\quad
 f^{\alpha,\beta}(t): = \beta t^{\alpha-1}  E_{\alpha,\alpha}(-\beta t^{\alpha}),
 \quad t >0 .
 \eeqnn
 The Laplace transform of Mittag-Leffler distribution admits the representation
 \beqlb\label{ML.LP}
 \hat{F}^{\alpha,\beta}(\lambda) = \int_0^\infty \lambda e^{-\lambda t} F^{\alpha,\beta}(t)dt = \int_0^\infty  e^{-\lambda t} f^{\alpha,\beta}(t)dt=  \frac{\beta}{\beta+\lambda^\alpha},\quad \lambda > 0;
 \eeqlb
 the reader may refer to  \cite{HauboldMathaiSaxena2011,MathaiHaubold2008} for its abundant results. 
 
 \bigskip
 
{\bf Processes of Mittag-Leffler type.} A Hawkes process $(N,\Lambda)$ is said to be of \textit{Mittag-Leffler type} with index $(\alpha,\beta)$  if  $\phi=f^{\alpha,\beta}$. 
The process has also referred to as being of \textit{fractional type} in \cite{ChenHawkesScalas2021}; its  analytical and numerical properties have been extensively studied in \cite{HabyarimanaAdudaScalasChenHawkes2023}.
 In view of  (\ref{ML.LP}), we have that  
 \beqnn
{\it\Phi} = 1-F^{\alpha,\beta} ,\quad 
m=1,\quad
\sigma=\infty
\quad \mbox{and} \quad 
 {\it\hat\Phi}(\lambda)
 = 1-\hat{F}^{\alpha,\beta}(\lambda)
 =\frac{\lambda^\alpha}{\beta+\lambda^\alpha},\quad  \lambda \geq 0.
 \eeqnn
 From this and (\ref{eqn.LapIR}), we deduce that  
 \[
 	\hat{\mathcal{I}}_R(\lambda)= \beta/\lambda^\alpha \quad \mbox{and} \quad 
	\hat{\mathcal{I}}_R^2(\lambda)=\beta/\lambda^{\alpha+1}. 
\] 
The one-to-one correspondence between functions and their Laplace transforms yields that
 \beqnn
 \mathcal{I}_R(t) = \frac{\beta \cdot  t^{\alpha} }{\Gamma(\alpha+1)}
 \quad \mbox{and}\quad
 \mathcal{I}_R^2(t) = \frac{\beta \cdot  t^{\alpha+1} }{\Gamma(\alpha+2)},\quad t\geq 0.
 \eeqnn

 \begin{proposition}
 We have $\mathcal{I}_R \in {\rm 2RV}^\infty_{\alpha,-\infty}(0)$ and $\mathcal{I}_R^2 \in {\rm 2RV}^\infty_{\alpha+1,-\infty}(0)$. In particular, for Mittag-Leffler Hawkes process we have as $n\to\infty$, weakly in $\mathbf{D}(\mathbb{R}_+;\mathbb{R})$ that  
 \beqnn
 \begin{split}
 & n^{\frac{1-\alpha}{2}\wedge \alpha}\Big(  \frac{N(nt)}{n^{\alpha+1}} -  \frac{\beta\cdot \mu_0 }{\Gamma(\alpha+2)}\cdot t^{\alpha+1} \Big) \\
 & ~~ \longrightarrow  ~ \mathrm{1}_{\{\alpha\leq 1/3\}}\cdot \mu_0 t + \mathrm{1}_{\{\alpha\geq 1/3\}}\cdot   \sqrt{\frac{\beta^3\cdot \mu_0 \cdot(\alpha+1) }{|\Gamma(\alpha+2)|^3}}  \int_0^t (t-s)^{\alpha}s^{\alpha/2}dB(s) \quad \mbox{as $n\to\infty$.}
 \end{split}
 \eeqnn
 \end{proposition}
  

%
%
 %
 

\bigskip

{\bf Processes of mixed Mittag-Leffler type.} Let us now generalize the class of processes to processes of mixed Mittag-Leffler type. Specifically, for any $0< \alpha_1\leq \alpha_2< 1$ and $\beta_1,\beta_2>0$ we say that the Hawkes process $(N,\Lambda)$ is of \textit{mixed Mittag-Leffler type} with index $(\alpha_i,\beta_i)_{i=1,2} $  if  its kernel is given by  
\[
	\phi= f^{\alpha_1,\beta_1}* f^{\alpha_2,\beta_2}.
\]  
 By (4.55) and (4.56) in \cite{HorstXu2023a}, we have as $t\to\infty$, 
 \beqnn
 {\it\Phi}(t) \sim\frac{\mathcal{C}_\beta\cdot t^{-\alpha_1 }}{ \Gamma(1-\alpha_1)} \in {\rm RV}^\infty_{-\alpha_1 },
 \quad  
 \mathcal{I}_R(t)  \sim C_{\mathcal{I}_R} \cdot t^{\alpha_1} \in {\rm RV}_{\alpha_1}^\infty
  \quad \mbox{and}\quad
  \mathcal{I}^2_R(t) \sim C_{\mathcal{I}_R^2} \cdot t^{\alpha_1+1}  \in {\rm RV}_{\alpha_1+1}^\infty,
 \eeqnn
  with $\mathcal{C}_\beta:=1/\beta_1+1/\beta_2\cdot \mathbf{1}_{\{\alpha_1=\alpha_2\}}$, $C_{\mathcal{I}_R}= \big(\mathcal{C}_\beta  \Gamma(\alpha_1+1)\big)^{-1}$ and 
  $C_{\mathcal{I}_R^2}= \big(\mathcal{C}_\beta  \Gamma(\alpha_1+2)\big)^{-1}$.
 The next result on the second-order regular variation of $\mathcal{I}^2_R$ comes from Proposition~4.12 in \cite{HorstXu2023a}. 
 
  \begin{proposition}\label{Prop.HP}
 We have   $\mathcal{I}^2_R \in {\rm 2RV}^\infty_{\alpha_1+1, \rho_\alpha } (  A_*)$ with 
 \beqnn 
  \rho_\alpha:= \begin{cases}
 -\alpha_1, & \mbox{if } \alpha_1=\alpha_2;  \vspace{5pt}\\
  \alpha_1-\alpha_2, & \mbox{if } \alpha_1<\alpha_2 
  \end{cases}
  \quad\mbox{and}\quad
  A_*(t):=
  \begin{cases}
  \displaystyle{ \frac{\alpha_1\Gamma(2+\alpha_1)}{\beta_1+\beta_2} \cdot t^{-\alpha_1},} & \mbox{if } \alpha_1=\alpha_2;  \vspace{5pt}\\
  \displaystyle{\frac{\beta_1(\alpha_2-\alpha_1)\Gamma(2+\alpha_1)}{\beta_2\Gamma(2+2\alpha_1-\alpha_2)} \cdot t^{\alpha_1-\alpha_2},} & \mbox{if } \alpha_1<\alpha_2 .
  \end{cases}
 \eeqnn
 \end{proposition}

 The preceding proposition allows us to easily establish the following FCLT for this mixed Mittag-Leffler Hawkes processes. 

 \begin{corollary}
 For the mixed Mittag-Leffler Hawkes process $N$, the following hold.
 \begin{enumerate}
   \item[(1)] If $\alpha_1=\alpha_2$, we have  as $n\to\infty$, weakly in $\mathbf{D}(\mathbb{R}_+;\mathbb{R})$ that
   \beqnn
    \lefteqn{n^{\alpha_1\wedge \frac{1-\alpha_1}{2}}\Big(   \frac{N(nt)}{n^{\alpha_1+1}} -  C_{\mathcal{I}_R^2}\cdot \mu_0  \cdot t^{\alpha_1+1} \Big) }\ar\ar\cr
     \ar\to\ar \mathbf{1}_{\{\alpha_1\leq 1/3\}} \cdot\mu_0\Big( 1- \frac{\beta_1\beta_2}{(\beta_1+\beta_2)^2}  \Big)\cdot t   + \mathbf{1}_{\{\alpha_1\geq 1/3\}} \cdot \sqrt{|C_{\mathcal{I}_R^2}|^3 \cdot \mu_0 \cdot(\alpha_1+1) }  \int_0^t (t-s)^{\alpha_1}s^{\alpha_1/2}dB(s). 
   \eeqnn
   
   \item[(2)] If $\alpha_1<\alpha_2$, we have  as $n\to\infty$, weakly in $\mathbf{D}(\mathbb{R}_+;\mathbb{R})$  that
    \beqnn
  \lefteqn{n^{\frac{1-\alpha_1}{2} \wedge \alpha_1 \wedge (\alpha_2-\alpha_1)} \Big(   \frac{N(nt)}{n^{\alpha_1+1}} -  C_{\mathcal{I}_R^2}\cdot \mu_0  \cdot t^{\alpha_1+1} \Big)}\ar\ar\cr
  \ar\to\ar  \mathbf{1}_{\{\alpha_1\leq 1/3;\, 2\alpha_1 \leq \alpha_2 \}} \cdot  \mu_0 t - \mathbf{1}_{\{ \alpha_2 \leq (\frac{\alpha_1+1}{2}\wedge (2\alpha_1))\}}\cdot\frac{\mu_0\cdot\beta_1^2 }{\beta_2 \Gamma(2+2\alpha_1-\alpha_2)} \cdot t^{1+2\alpha_1-\alpha_2}  \cr
  \ar\ar+  \mathbf{1}_{\{\alpha_1\geq 1/3; \alpha_2\geq \frac{\alpha_1+1}{2}  \}}\cdot   \sqrt{|C_{\mathcal{I}_R^2}|^3 \cdot \mu_0 \cdot(\alpha_1+1) }  \int_0^t (t-s)^{\alpha_1}s^{\alpha_1/2}dB(s).
 \eeqnn 
 \end{enumerate}
 \end{corollary}

 \subsubsection{Stable type}
 
Let us  consider a one-sided stable random variable $Z^{\alpha}$ on $\mathbb{R}_+$ with stability parameter $\alpha \in (0,1)$ and unit scale parameter, whose probability density function and distribution function is denoted by $g^{\alpha}$ and $G^{\alpha}$, respectively. They are uniquely determined by the Laplace transform
  \beqnn
  \mathbf{E}\big[\exp\big\{-\lambda Z^{\alpha}\big\}\big] = \int_0^\infty  e^{-\lambda t}g^{\alpha}(t)dt =\int_0^\infty \lambda e^{-\lambda t}G^{\alpha}(t)dt= e^{-|\lambda|^\alpha},\quad \lambda >0.
  \eeqnn
  Let $\xi$ be a non-negative random variable with finite mean $m_\xi$.
  The Hawkes process $(N,\Lambda)$ is said to be of \textit{scaled-stable type} with \textit{index} $\alpha$ and \textit{scaled factor} $\xi$ if its kernel $\phi$ is equal to the probability density function\footnote{Since $g^\alpha$ is smooth on $(0,\infty)$; see (14.31) in \cite[p.88]{Sato1999}, the probability density function of $Z^{\alpha}\cdot |\xi|^{1/\alpha}$ exists.} of $Z^{\alpha}\cdot |\xi|^{1/\alpha}$.
  In this case, we have as $\lambda \to 0+$, 
  \beqnn
  {\it\hat\Phi}(\lambda) = 1- \int_0^\infty e^{-\lambda t}\phi(t)dt = 1- \mathbf{E}\big[\exp\big\{-\lambda |\xi|^{1/\alpha}\cdot Z^{\alpha} \big\}\big]
  = 1- \mathbf{E}\big[\exp\big\{-\lambda^\alpha \xi \big\}\big]  \sim  m_\xi\cdot \lambda^{\alpha}. 
  \eeqnn 
  By Proposition~\ref{Thm.KaramataTauberian} and  \ref{Thm.AsymR}(3), this shows that as $t\to\infty$,
  \beqnn
   {\it \Phi}(t)  \sim \frac{m_\xi}{\Gamma(1-\alpha)}\cdot t^{-\alpha} ,
   \quad
   \mathcal{I}_R(t) \sim \frac{t^{\alpha} }{m_\xi\cdot \Gamma(1+\alpha)}
   \quad \mbox{and}\quad  \mathcal{I}_R^2(t) \sim \frac{t^{\alpha+1} }{ m_\xi\cdot \Gamma(2+\alpha)}.
  \eeqnn

  \begin{lemma}\label{Lemma06}
  For each $\delta>0$, we have that $\Gamma(\alpha)\Gamma(2-\alpha)\cdot \mathcal{I}_ {\it \Phi}(t)\cdot\big(\mathcal{I}_R(t+\delta)-\mathcal{I}_R(t)\big)\to \delta $ as $t\to\infty$.
  \end{lemma}
  \proof  For  $\lambda >0$, we have that
  \beqnn
  \int_0^\infty \frac{1}{\lambda} e^{-t/\lambda} \big(\mathcal{I}_R(t+\delta)-\mathcal{I}_R(t)\big) dt
  \ar=\ar e^{\delta/\lambda} \int_\delta^\infty \frac{1}{\lambda} e^{-t/\lambda}\mathcal{I}_R(t)dt - \int_0^\infty \frac{1}{\lambda} e^{-t/\lambda }\mathcal{I}_R(t) dt \cr
  \ar=\ar  e^{\delta/\lambda } \int_0^\delta \frac{1}{\lambda} e^{-t/\lambda}\mathcal{I}_R(t)dt + (e^{\delta/\lambda}-1) \hat{\mathcal{I}}_R(1/\lambda).
  \eeqnn
  By Proposition~\ref{Thm.KaramataTauberian}, it holds that $\hat{\mathcal{I}}_R(1/\lambda) \sim \lambda^\alpha/m_\xi$ as $\lambda \to\infty$ and
  \beqnn
  \int_0^\infty \frac{1}{\lambda} e^{-t/\lambda} \big(\mathcal{I}_R(t+\delta)-\mathcal{I}_R(t)\big) dt
  \sim \frac{\delta}{  m_\xi  }\cdot \lambda^{ \alpha-1}.
  \eeqnn
  An application of Proposition~\ref{Thm.KaramataTauberian} and the non-negativity of $\mathcal{I}_R(t+\delta)-\mathcal{I}_R(t)$ shows that 
\[
 	\mathcal{I}_R(t+\delta)-\mathcal{I}_R(t) \sim \delta/(m_\xi \cdot \Gamma(\alpha))\cdot t^{ \alpha-1} \quad \mbox{as} \quad t\to\infty.
\] 
 By Proposition~\ref{Thm.Karamata}, we have $\mathcal{I}_ {\it \Phi}(t) \sim \frac{m_\xi}{\Gamma(2-\alpha) \cdot t^{1-\alpha}}$ from which the desired result follows. 
 \qed

 The preceding lemma shows that the renewal function $1+\mathcal{I}_R$  satisfies Lemma~1(iii) in \cite{AndersonAthreya1987}. It has been shown in 
 \cite[p.393]{AndersonAthreya1987} that the Theorem~2, 3 and 4 therein also hold for  $\alpha\in(0,1)$ under the results established in our Lemma \ref{Lemma06}.
 Moreover, by Lemma~2.1 in \cite{HorstXu2023a} and Proposition~\ref{Prop.2RVR2.000}, we can get the second-order regular variation of $\mathcal{I}_R$ and $\mathcal{I}^2_R$.
 
 \begin{proposition}\footnote{By the domain of attraction theorem \cite[Theorem~8.3.1]{BinghamGoldieTeugels1987} the assumptions (2) and (3) are equivalent to those in \cite[Theorem~3 and 4]{AndersonAthreya1987}.}
 Three regimes arise for the long-term behaviors of $\mathcal{I}_R$ and $\mathcal{I}_R^2$.
 \begin{enumerate}
  \item[(1)] If $\xi$ has finite variance $\sigma_\xi^2$, then  as $t\to\infty$,
  \beqnn
  \mathcal{I}_R (t ) - \frac{t^{\alpha} }{m_\xi\cdot \Gamma(1+\alpha)} \sim \frac{\sigma_\xi^2- m_\xi^2}{2 m_\xi^2}.
  \eeqnn
  Moreover, if $\sigma_\xi^2 \neq m_\xi^2$, then $\mathcal{I}_R \in {\rm 2RV}^\infty_{\alpha, -\alpha}(A_{\xi,1})$ and $\mathcal{I}_R^2 \in {\rm 2RV}^\infty_{\alpha+1, -\alpha}\big((1+\alpha )\cdot A_{\xi,1}\big)$ with
      \beqnn
      A_{\xi,1}(t) := - \frac{\alpha\Gamma(1+\alpha)}{2 m_\xi}(\sigma_\xi^2- m_\xi^2) \cdot t^{-\alpha}, \quad t>0.
      \eeqnn

  \item[(2)] If $\xi$ has infinite variance and $\int_0^t s^2 \mathbf{P}(\xi \in ds)\in  {\rm RV}^\infty_0$, then  as $t\to\infty$,
  \beqnn
  \mathcal{I}_R(t) - \frac{t^{\alpha} }{m_\xi\cdot \Gamma(1+\alpha)} \sim \frac{1}{m_\xi^2} \int_0^{t^\alpha} dx \int_x^\infty \mathbf{P}(\xi>u)du-1
  \in {\rm RV}^\infty_0, 
  \eeqnn
 and  $\mathcal{I}_R \in {\rm 2RV}^\infty_{\alpha, -\alpha} \big(m_\xi\cdot \alpha \Gamma(1+\alpha)\cdot A_{\xi,2} \big)$ and $\mathcal{I}_R^2 \in {\rm 2RV}^\infty_{\alpha+1, -\alpha} \big(m_\xi\cdot \alpha \Gamma(2+\alpha) \cdot A_{\xi,2} \big)$ with
      \beqnn
      A_{\xi,2}(t) :=  \Big( 1- \frac{1}{m_\xi^2} \int_0^{t^\alpha} dx \int_x^\infty \mathbf{P}(\xi>u)du\Big) \cdot t^{-\alpha},\quad t>0.
      \eeqnn

  \item[(3)] If $\mathbf{P}(\xi\geq t) \in {\rm RV}^\infty_{-\beta}$ with $\beta\in(1,2)$, then as $t\to\infty$,
  \beqnn
  \mathcal{I}_R (t) - \frac{t^{\alpha} }{m_\xi\cdot \Gamma(1+\alpha)} \sim -1 -\frac{\Gamma(1-\beta)\cdot t^{2\alpha}\mathbf{P}(\xi\geq t^\alpha)}{m_\xi^2\cdot\Gamma(1+2\alpha-\alpha\beta)}
  \in {\rm RV}^\infty_{\alpha(2-\beta)},
  \eeqnn
  $\mathcal{I}_R \in {\rm 2RV}^\infty_{\alpha, \alpha(1-\beta)} \big( m_\xi\cdot \alpha(\beta-1) \Gamma(1+\alpha) \cdot A_{\xi,3}\big)$ and $\mathcal{I}_R^2 \in {\rm 2RV}^\infty_{\alpha+1, \alpha(1-\beta)} \big(\frac{ m_\xi\cdot \alpha(\beta-1) \Gamma(2+\alpha)   }{1+\alpha(2-\beta)}\cdot A_{\xi,3}\big)$ with
      \beqnn
      A_{\xi,3}(t) := \Big( 1+ \frac{\Gamma(1-\beta)\cdot t^{2\alpha}\mathbf{P}(\xi\geq t^\alpha)}{m_\xi^2\cdot\Gamma(1+2\alpha-\alpha\beta)} \Big) \cdot t^{-\alpha},\quad t>0.
      \eeqnn

 \end{enumerate}
 \end{proposition}

 \begin{corollary}
 For the scaled-stable Hawkes process $N$, 
 the following limits hold weakly in $\mathbf{D}(\mathbb{R}_+;\mathbb{R})$ as $n\to\infty$. 
   \begin{enumerate}
  \item[(1)] If $\xi$ has finite variance $\sigma_\xi^2 \neq m_\xi^2$, then
   \beqnn
   n^{\alpha} \Big(   \frac{N(nt)}{n^{\alpha+1}} -  \frac{  \mu_0 \cdot t^{\alpha+1}}{m_\xi\cdot\Gamma(\alpha+2)} \Big)
   \ar\to\ar  \mathbf{1}_{\{\alpha\leq 1/3\}}\cdot\mu_0\cdot \Big( t +  \frac{\sigma_\xi^2- m_\xi^2}{2  m_\xi^2 } \cdot (t- t^{\alpha+1}) \Big)\cr
   \ar\ar + \mathbf{1}_{\{\alpha\geq 1/3\}}\cdot \sqrt{\frac{\mu_0(\alpha+1)}{m_\xi^3|\Gamma(2+\alpha)|^3}} \int_0^t (t-s)^\alpha s^{\alpha/2} dB(s).
   \eeqnn

   \item[(2)] If  $\xi$ has infinite variance and $\int_0^t s^2 \mathbf{P}(\xi \in ds)\in  {\rm RV}^\infty_0$, then if $\alpha\in(0,1/3]$,
     \beqnn
      \frac{1}{A_{\xi,2}(n)}  \Big(   \frac{N(nt)}{n^{\alpha+1}} -  \frac{  \mu_0 \cdot t^{\alpha+1}}{m_\xi\cdot\Gamma(\alpha+2)} \Big)
     \ar\to\ar  \mu_0\cdot (t^{\alpha+1} - t) ; 
     \eeqnn
     and if $\alpha\in(1/3,1)$,
     \beqnn
      n^{(1-\alpha)/2} \Big(   \frac{N(nt)}{n^{\alpha+1}} -  \frac{  \mu_0 \cdot t^{\alpha+1}}{m_\xi\cdot\Gamma(\alpha+2)} \Big)
     \ar\to\ar  \sqrt{\frac{\mu_0(\alpha+1)}{m_\xi^3|\Gamma(2+\alpha)|^3}} \int_0^t (t-s)^\alpha s^{\alpha/2} dB(s).
     \eeqnn

   \item[(3)] If $\mathbf{P}(\xi\geq t) \in {\rm RV}^\infty_{-\beta}$ with $\beta\in(1,2)$, then if $\beta<\frac{\alpha+1}{2\alpha}$,
       \beqnn
      \frac{1}{A_{\xi,3}(n)}  \Big(   \frac{N(nt)}{n^{\alpha+1}} -  \frac{  \mu_0 \cdot t^{\alpha+1}}{m_\xi\cdot\Gamma(\alpha+2)} \Big)
     \ar\to\ar  \frac{ \mu_0}{1+\alpha(2-\beta)} \cdot \big( t^{\alpha+1}- t^{1+\alpha(2-\beta)} \big) ; 
     \eeqnn
     and if $\beta>\frac{\alpha+1}{2\alpha}$,
     \beqnn
      n^{(1-\alpha)/2} \Big(   \frac{N(nt)}{n^{\alpha+1}} -  \frac{  \mu_0 \cdot t^{\alpha+1}}{m_\xi\cdot\Gamma(\alpha+2)} \Big)
     \ar\to\ar  \sqrt{\frac{\mu_0(\alpha+1)}{m_\xi^3|\Gamma(2+\alpha)|^3}} \int_0^t (t-s)^\alpha s^{\alpha/2} dB(s).
     \eeqnn
  \end{enumerate}
 \end{corollary}

   \section{Elementary properties of Hawkes processes} \label{Sec.EleProp}
 \setcounter{equation}{0} 
 
 In this section we prove Theorem~\ref{Thm201} and Proposition~\ref{LapTheta}, \ref{Thm.AsymR}. We first establish the 
%
%
%
representation (\ref{FourLapFunHP}) under the assumption of the existence of non-continuable solutions to (\ref{VolRiccati}). 


  \begin{proposition}\label{Lemma.FourLapGHP}
 	For $\nu \in M(\mathbb{R}_+;\mathbb{C}_-)$ and $f\in L^\infty_{\rm loc}(\mathbb{R}_+;\mathbb{C}_-) $, if $(V,T_\infty)\in L^\infty_{\rm loc}([0,T_\infty);\mathbb{C})\times(0,\infty]$ is a  non-continuable solution of (\ref{VolRiccati}), then (\ref{FourLapFunHP}) holds for any $T\in[0,T_\infty)$. 
 	
 \end{proposition}
 \proof 
 Recall the function $W=e^{V+f}-1-V$.  
 For each $T\in[0,T_\infty)$, let $\Lambda_T:= \Lambda*d\nu(T)+W *\Lambda (T)$. Taking expectations on both sides of (\ref{SVR}), we see that
 \beqnn
  \mathbf{E}[\Lambda(t)]=H_\mu(t) 
  \quad \mbox{and}\quad
  \mathbf{E}[\Lambda_T]=H_\mu * d\nu(T) + W * H_\mu (T).
 \eeqnn 

 To prove that the Fourier-Laplace functional equals $e^{\mathbf{E}[\Lambda_T]}$ as claimed, we show that this quantity can be expressed as the initial value of an exponential martingale $\{e^{Z_T(t)} : t \in [0,T]\}$. To this end, we first introduce the Doob martingale  $$ \{\Lambda_T(t):=\mathbf{E}[\Lambda_T|\mathscr{F}_t]:t\in[0,T]\}.$$
 Using the martingale representation of the intensity process $\Lambda$ established in Lemma \ref{MartRep} along with the stochastic Fubini theorem and - in the second step - the fact that 
 $V=R*d\nu + R*W$ we obtain that
 \beqlb\label{eqn.3.03}
 \Lambda_T(t)
 \ar=\ar  \mathbf{E}[\Lambda_T] + \int_0^T W(T-r) \int_0^{t\wedge r}  R(r-s)\widetilde{N}(ds)dr\cr
 \ar\ar + \int_0^T \nu(d(T-r)) \int_0^{t\wedge r}  R(r-s)\widetilde{N}(ds)\cr
 \ar=\ar  \mathbf{E}[\Lambda_T] + \int_0^t  V(T-s) \widetilde{N}(ds).
 \eeqlb
 Next, we introduce a semi-martingale $ \{Z_T(t):t\in[0,T]\}$ by
 \beqlb\label{eqn.3002}
 Z_T(t)
 \ar:=\ar \Lambda_T(t) - \int_0^t W(T-s) \Lambda(s)ds + \int_0^t  f(T-s) N(ds)\cr
 \ar=\ar \mathbf{E}[\Lambda_T]+   \int_0^t (f- W)(T-s) \Lambda(s)ds + \int_0^t   (V +f)(T-s) \widetilde{N}(ds).
 \eeqlb
 Applying It\^{o}'s formula to $\exp\{Z_T(t) \}$ and then using the fact that $W=  e^{V+f}-1-V$, we obtain that
 \beqlb\label{eqn.3008}
 e^{Z_T(t)}= e^{\mathbf{E}[\Lambda_T]} + \int_0^t e^{Z_T(s-)} \big[ e^{(V+f)(T-s)}-1 \big] \widetilde{N}(ds),\quad t\in[0,T] , 
 \eeqlb
 which is a $(\mathscr{F}_t)$-local martingale.  If the local martingale were a true martingale, then the desired result follows from 
 \beqnn
 Z_T(T)= 
 \Lambda*d\nu(T) + f*N(T) 
 \eeqnn 
 by taking expectations on both sides of (\ref{eqn.3008}) with $t=T$. In fact, in this case,  
 \beqnn
 \mathbf{E}\big[ \exp\big\{ \Lambda*d\nu(T) + f*N(T) \big\} \big] = \mathbf{E}\big[\exp\{  Z_T(T)\} \big] =   \exp\big\{\mathbf{E}[\Lambda_T]\big\} = \exp\big\{  H_\mu *d\nu(T) + W * H_\mu (T) \big\}.
 \eeqnn
 
 It hence remains to prove that $\{ \exp\{Z_T(t) \}: t\in[0,T]  \}$ is a true $(\mathscr{F}_t)$-martingale. The proof is similar to those of Lemma 6.8 in \cite{Xu2021b} and Lemma~6.3 in \cite{Jaber2021}. We first define the process
 \beqnn
 U_T(t):=  \int_0^t  \big[ e^{(V+f)(T-s)}-1\big] \widetilde{N}(ds ), \quad t\geq 0.
 \eeqnn
 In view the local integrability of $H_\mu$ and the local boundedness of $(V,f)$, this is an $(\mathscr{F}_t)$-local martingale with quadratic variation
 \beqnn
 \big[ U_T \big]_{t} = \int_0^t  \big| e^{(V+f)(T-s)}-1\big|^2 N(ds ), \quad t\geq 0.
 \eeqnn
 It follows from  the Burkholder-Davis-Gundy inequality (e.g. see Theorem 26.12 in \cite[p.524]{Kallenberg2002}) that 
 \beqnn
 \mathbf{E}\Big[\sup_{t\in[0,T]}|U_T(t)|^2\Big] 
  \leq  C\cdot \mathbf{E}\big[ \big[ U_T \big]_T\big]
 \ar=\ar  C\cdot\int_0^T\mathbf{E} \big[\Lambda(s)\big] \big| e^{(V+f)(T-s)}-1\big|^2 ds  \cr
 \ar\leq\ar  C\cdot \sup_{t\in[0,T]} \big| e^{(V+f)(t)}-1\big|^2 \cdot \mathcal{I}_{H_\mu}(T) <\infty
 \eeqnn
 and hence that $ \{U_T(t):t\in[0,T]\}$ is a locally uniformly square integrable $(\mathscr{F}_t)$-martingale. Let
\[
	\mathcal{E}_{U_T}:= \{\mathcal{E}_{U_T}(t): t\geq 0 \} 
\]
	be the Dol\'ean-Dade exponential of $U_T$. By It\^o's formula,
 \beqnn
 \mathcal{E}_{U_T}(t)
 =\exp\Big\{ \int_0^t \big( f- W \big) (T-s)\Lambda(s) ds + \int_0^t  (V+f)(T-s) \widetilde{N}(ds)  \Big\}, \quad t\geq 0.
 \eeqnn
  From this and (\ref{eqn.3002}), we see that  
  \[
  	e^{Z_T(t)}=e^{\mathbf{E}[\Lambda_T]}  \mathcal{E}_{U_T}(t), \quad t\in[0,T].
\] 
 Since $ \mathcal{E}_{U_T}$ is a non-negative local martingale it is a super-martingale and so $\mathbf{E}[\mathcal{E}_{U_T}(t)] \leq 1$. It hence suffices to prove that 
 \[
 	\mathbf{E}[\mathcal{E}_{U_T}(t)]=1, \quad t\geq 0.
\]

To this end, we introduce, for each $t_0\geq 0$ and $n\geq 1$ the quantities
 \beqnn
 \tau_n:=\inf\{ s\geq 0: \mathcal{I}_{\Lambda}(s)\geq n \} \wedge t_0
 \quad\mbox{and}\quad
 \mathcal{E}_{U_T}^{(n)}(t):=  \mathcal{E}_{U_T}(\tau_n \wedge t), 
 \quad t\geq 0. 
 \eeqnn 
 Since  $(V,f)$ is locally bounded, there exists a constant $C>0$ such that uniformly in $t \in [0,T]$ the following holds:
 \beqnn
 \int_0^t \mathbf{1}_{\{s\leq \tau_n\}} \Lambda(s) \big| 1-\big( 1-(V+f)(T-s) \big)     \exp\big\{(V+f)(T-s)\big\}  \big|  ds 
 \leq C\cdot \int_0^{\tau_n}\Lambda(s) ds \leq C\cdot n. 
 \eeqnn
 Hence, the process $\mathcal{E}_{U_T}^{(n)}$ is a martingale for each $n\geq 1$, due to Theorem~IV.3 in \cite{LepingleMemin1978}. 
 Thus, 
 \beqnn
 \begin{split}
 1 = \mathbf{E}\big[\mathcal{E}_{U_T}^{(n)}(t_0)\big]
 & = \mathbf{E}\big[\mathcal{E}_{U_T}^{(n)}(t_0); \tau_n = t_0\big] + \mathbf{E}\big[\mathcal{E}_{U_T}^{(n)}(t_0);\tau_n<t_0\big]  \\
 & = \mathbf{E}\big[\mathcal{E}_{U_T}(t_0);\tau_n= t_0\big] + \mathbf{E}\big[\mathcal{E}_{U_T}^{(n)}(\tau_n);\tau_n<t_0\big] .
 \end{split}
 \eeqnn
 By the monotone convergence theorem and the fact that $\tau_n \to t_0$ a.s.~as $n\to\infty$, we have that
 $$ \mathbf{E}\big[\mathcal{E}_{U_T}(t_0); \tau_n=t_0 \big]\to  \mathbf{E}\big[ \mathcal{E}_{U_T}(t_0) \big]$$
and so it suffices to prove that $\mathbf{E}[\mathcal{E}_{U_T}^{(n)}(\tau_n);\tau_n<t_0] \to 0$ as $n\to\infty$. 

To this end, we define a probability law  $\mathbf{Q}_T^{(n)}$ on $(\Omega,\mathscr{F},\mathscr{F}_t)$ by 
 $$  \frac{d  \mathbf{Q}_T^{(n)} }{ d \mathbf{P} } = \mathcal{E}_{U_T}^{(n)} (\tau_n). $$
 By the definition of $\tau_n$ and Chebyshev's inequality it follows that
 \beqnn
 \mathbf{E}\big[\mathcal{E}_{U_T}^{(n)}(\tau_n); \tau_n<t_0 \big]  =\mathbf{Q}_T^{(n)} \big(\tau_n<t_0 \big) = \mathbf{Q}_T^{(n)} \big( \mathcal{I}_{\Lambda}(t_0) \geq n\big) \leq \frac{1}{n}\mathbf{E}^{\mathbf{Q}_T^{(n)}}\big[ \mathcal{I}_{\Lambda}(t_0)\big] 
 \eeqnn
and hence the desired result holds if we can establish a uniform upper bound on $\mathbf{E}^{\mathbf{Q}_T^{(n)}}\big[ \mathcal{I}_{\Lambda}(t_0)\big]$. In what follows we prove that there exists $\beta > 0$ such that
\[
	\sup_{n \geq 1}\mathbf{E}^{\mathbf{Q}_T^{(n)}}\big[ \mathcal{I}_{\Lambda}(t_0)\big] \leq 2 e^{\beta t_0} \mathcal{I}_{H_\mu}(t_0). 
\]

As  in~\cite[Section~2]{HorstXu2021}, on an extension of the original probability space we can define a time-homogeneous Poisson random measure (PRM) $N_{0}(ds,dz)$  on $(0,\infty)^2$ with intensity $ds dz$ such that 
 \beqlb
 N(t)\ar=\ar \int_0^t \int_0^{\Lambda(s-)} N_{0}(ds,dz)
  \quad{and} \quad
 \Lambda(t) =  H_\mu (t)+ \int_0^t  \int_0^{\Lambda(s-)}R(t-s)  \widetilde{N}_{0}(ds,dz ), \label{eqn.3001}
 \eeqlb
 for any $t\geq0$, where $\widetilde{N}_{0}(ds,dz ):= N_{0}(ds,dz )-dsdz$ denotes the compensated random measure.  By Girsanov's theorem for random measures as stated in, e.g.~Theorem~3.17 in \cite[p.170]{JacodShiryaev2003}, the PRM $N_{0}(ds,dz)$ is a random point measure under $\mathbf{Q}_T^{(n)}$ with intensity $$\mathbf{1}_{\{ s\leq \tau_n \}}\cdot \exp \{ (V+f)(T-s) \}  ds dz.$$ Moreover, the second equation in (\ref{eqn.3001}) - which holds a.s.~under the measure $\mathbf{P}$ - induces the following equality in distribution under the measure $\mathbf{Q}_T^{(n)}$: 
 \beqnn
 \Lambda (t)\ar \stackrel{d}{=}\ar   H_\mu (t) +  \int_0^t   \mathbf{1}_{\{ s\leq \tau_n \}} \Lambda (s) R(t-s)  \big(e^{(V+f)(T-s)}-1 \big)ds + \int_0^t  \int_0^{\Lambda(s-)}R(t-s)  \widetilde{N}_{0}(ds,dz ), \quad t \geq 0.  
 \eeqnn 
 Taking expectations on both sides of the above equation and then integrating them over $[0,t]$ yields that
 \beqnn
 \int_0^t \mathbf{E}^{\mathbf{Q}^{(n)}_T}\big[ \Lambda (s) \big]ds  
 \ar=\ar  \mathcal{I}_{H_\mu}(t) +  \int_0^t dr \int_0^r   \mathbf{E}^{\mathbf{Q}^{(n)}_T}\big[ \mathbf{1}_{\{ s\leq \tau_n \}}\Lambda (s) \big] R(r-s)  \big(e^{(V+f)(T-s)}-1 \big)ds \cr
 \ar\leq\ar \mathcal{I}_{H_\mu}(t) + C_0  \int_0^t dr \int_0^r   \mathbf{E}^{\mathbf{Q}^{(n)}_T}\big[  \Lambda (s) \big] R(t-s) ds, 
 \eeqnn
 for some constant $C_0>0$ that is independent of $t$ and $n$. 
 Applying Fubini's theorem to the term on the left-hand side of the equality and to the double integral on the right-hand side of the inequality and using the fact that $ \mathcal{I}_{\Lambda }$ is non-decreasing we see that
 \beqnn
  \mathbf{E}^{\mathbf{Q}^{(n)}_T}\big[  \mathcal{I}_{\Lambda }(t) \big] 
 \ar\leq\ar  \mathcal{I}_{H_\mu} (t)+ C_0  \int_0^t  \mathbf{E}^{\mathbf{Q}^{(n)}_T}\Big[  \int_0^s    \Lambda (r)dr \Big] R(t-s) ds\cr
 \ar = \ar  \mathcal{I}_{H_\mu} (t) + C_0 \int_0^t   \mathbf{E}^{\mathbf{Q}^{(n)}_T}\big[  \mathcal{I}_{\Lambda }( s) \big] R(t-s) ds. 
 \eeqnn
 From (\ref{Resolvent}),  for any $\beta>0$ such that $\int_0^\infty e^{-\beta t}\phi(t)dt \leq \frac{1/2}{1+C_0}$ we have 
 \beqnn 
  \int_0^\infty e^{-\beta t}R(t)dt = \frac{\int_0^\infty e^{-\beta t}\phi(t)dt}{1-\int_0^\infty e^{-\beta t}\phi(t)dt} \leq \frac{1}{1+2C_0} , 
 \eeqnn
 and hence 
 \beqnn
 \sup_{t\in[0,t_0]} e^{-\beta t} \mathbf{E}^{\mathbf{Q}^{(n)}_T}\big[  \mathcal{I}_{\Lambda }(t) \big] 
 \ar\leq\ar e^{-\beta t} \mathcal{I}_{H_\mu} (t) + C_0 \int_0^t  e^{-\beta s}  \mathbf{E}^{\mathbf{Q}^{(n)}_T}\big[  \mathcal{I}_{\Lambda }( s) \big]  e^{-\beta (t-s)} R(t-s) ds\cr 
 \ar\leq\ar \sup_{t\in[0,t_0]} e^{-\beta t}\mathcal{I}_{H_\mu} (t) + \sup_{t\in[0,t_0]} e^{-\beta t}  \mathbf{E}^{\mathbf{Q}^{(n)}_T}\big[  \mathcal{I}_{\Lambda }(t) \big] \cdot  C_0    \int_0^\infty  e^{-\beta s} R(s) ds\cr
 \ar\leq\ar \sup_{t\in[0,t_0]} e^{-\beta t} \mathcal{I}_{H_\mu} (t) + \frac{1}{2}\cdot \sup_{t\in[0,t_0]} e^{-\beta t}  \mathbf{E}^{\mathbf{Q}^{(n)}_T}\big[  \mathcal{I}_{\Lambda }(t) \big] ,
 \eeqnn
 which yields that $\sup_{n\geq 1} \sup_{t\in[0,t_0]} e^{-\beta t} \mathbf{E}^{\mathbf{Q}^{(n)}_T}\big[  \mathcal{I}_{\Lambda }(t) \big]  \leq 2\cdot \mathcal{I}_{H_\mu}(t_0)$.
 Consequently, 
 \beqnn
 \sup_{n\geq 1}  \mathbf{E}^{\mathbf{Q}^{(n)}_T}\big[  \mathcal{I}_{\Lambda}(t_0) \big] \leq e^{\beta t_0}  \sup_{n\geq 1}  \sup_{t\in[0,t_0]} e^{-\beta t} \mathbf{E}^{\mathbf{Q}^{(n)}_T}\big[  \mathcal{I}_{\Lambda}(t) \big]  \leq  2 e^{\beta t_0}  \mathcal{I}_{H_\mu} (t_0).
 \eeqnn
 \qed

  \medskip
 \textit{Proof of Theorem~\ref{Thm201}.}  
 The existence of a non-continuable solution $(V,T_\infty)$ follows from standard arguments. In view of Lemma~\ref{Lemma.FourLapGHP} it remains to prove that  $T_\infty=\infty$ and  that $ {\tt Re}V \leq 0$. Choosing $\mu=\phi$, we see that
 \beqnn
 H_\phi= \phi+\phi*R=R
 \quad\mbox{and}\quad
 H_\phi*\phi +H_\phi *W=R*\phi +R*W= V.
 \eeqnn
 Taking these back into equation (\ref{FourLapFunHP}) shows that  for $\nu \in M(\mathbb{R}_+;\mathbb{C}_-)$ and $f\in L^\infty_{\rm loc}(\mathbb{R}_+;\mathbb{C}_-) $, 
 \beqnn
 \mathbf{E}\Big[\exp \Big\{ \Lambda*d\nu(T)+\int_0^Tf(T-t)N(dt) \Big\}\Big]= e^{  V(T) },
 \quad  T\in[0,T_\infty). 
 \eeqnn 
 It hence follows that 
 \beqnn
 \exp\big\{{\tt Re}V(T)\big\} = \big|\exp\big\{  V (T) \big\}\big|
 \ar\leq\ar \mathbf{E}\Big[\Big|\exp\Big\{ \Lambda*d\nu(T)+\int_0^Tf(T-t)N(dt)\Big\}\Big|\Big]\cr
 \ar=\ar\mathbf{E}\Big[ \exp\Big\{ {\tt Re}\Big(\Lambda*d\nu(T)+\int_0^Tf(T-t)N(dt)\Big) \Big\} \Big]\leq 1 , 
 \eeqnn
 and hence that ${\tt Re}V(T) \leq 0$. Furthermore, from (\ref{VolRiccati}) we can see that 
 \beqnn
 {\tt Re}V= {\tt Re}(\phi*d\nu)+ \big(\exp\{{\tt Re}(V+f)\}\cos ({\tt Im}(V+f)) -1\big)*\phi.
 \eeqnn
 Since ${\tt Re}(V+f)\leq 0$, we have that $  {\tt Re}(\phi*d\nu)\geq {\tt Re}V\geq {\tt Re}(\phi*d\nu) -2\mathcal{I}_\phi$. Similarly,
 \beqnn
 {\tt Im}V = {\tt Im}(\phi*d\nu)+ \big(\exp\{{\tt Re}(V+f)\}\sin ({\tt Im}(V+f)) \big)*\phi
 \eeqnn
 and hence $ |{\tt Im}V| \leq |{\tt Im}(\phi*d\nu)|+ \mathcal{I}_\phi$.
 The preceding estimates yield that  $T_\infty=\infty$ and $V\in L^\infty_{\rm loc}(\mathbb{R}_+;\mathbb{C}_-)$. 
 \qed 
 
 \medskip
  {\it Proof of Proposition~\ref{LapTheta}.}
 We first prove the equivalence between claim (1) and claim (2).
 Using integration by parts, 
 \beqlb\label{eqn.LapTheta01}
 {\it \Psi}_k(t)=k\int_0^t s^{k-1}{\it \Phi}(s)ds-t^k{\it \Phi}(t),\quad t>0.
 \eeqlb
 If ${\it \Phi}\in{\rm RV}_{-\alpha}^\infty$, then $t^{k-1}{\it \Phi}(t)\in{\rm RV}_{k-\alpha-1}^\infty$ and hence claim (2) follows directly from Proposition~\ref{Thm.Karamata}.
 For the converse, using integration by parts again,
 \beqnn
 \int_0^\infty s^{-k-1}{\it \Psi}_k(s)ds = \int_0^\infty t^k\phi(t)dt \int_t^\infty s^{-k-1}ds =\frac{1}{k}  \int_0^\infty  \phi(t)dt  <\infty.
 \eeqnn
 Moreover, for any  $t,\epsilon>0$ we also have
 \beqnn
 t^{-k}{\it \Psi}_k(t) = \int_0^{\epsilon t} (s/t)^k\phi(s)ds + \int_{\epsilon t}^t (s/t)^k\phi(s)ds \leq \epsilon^k\int_0^{\epsilon t} \phi(s)ds + \int_{\epsilon t}^\infty \phi(s)ds,
 \eeqnn
 which vanishes as $t\to\infty$ and then $\epsilon \to 0+$. From these and integration by parts, we have
 \beqlb\label{eqn.3.11}
 {\it\Phi}(t)= k\int_t^\infty s^{-k-1}{\it\Psi}_k(s)ds- t^{-k}{\it\Psi}_k(t), 
 \eeqlb
 and so claim (1) follows from Proposition~\ref{Thm.Karamata}.  
 The equivalence between claims (1) and (3) and claims (2) and (4) along with the final statement of the proposition 
 follows from  Proposition~\ref{Thm.KaramataTauberian} . 
 \qed
 
  \medskip
 {\it Proof of Proposition~\ref{Thm.AsymR}.} Taking Laplace transforms of both sides of (\ref{Resolvent}) and then using (\ref{LapPhi}),
 \beqlb\label{LapResolvent}
 \widehat{\mathcal{I}}_R (1/\lambda)= \int_0^\infty e^{-t/\lambda } R(t)dt = \frac{\int_0^\infty e^{-t/\lambda}\phi(t)dt}{1- \int_0^\infty e^{-t/\lambda}\phi(t)dt } =\frac{m-{\it\hat\Phi}(1/\lambda)}{1-m+{\it\hat\Phi}(1/\lambda)},  \quad \lambda >0. 
 \eeqlb
 When $m<1$, then ${\it\hat\Phi}(1/\lambda)\to0$ and $ \widehat{\mathcal{I}}_R (1/\lambda) \to m/(1-m)$ as $\lambda \to \infty$. 
 Thus, claim (1) follows from Proposition~\ref{Thm.Karamata} and \ref{Thm.KaramataTauberian}.
 When $m=1$ and $\sigma<\infty$, then ${\it\hat\Phi}(1/\lambda) \sim \lambda^{-1}{\it \Psi}_1(\infty)= \sigma/\lambda $ and $\widehat{\mathcal{I}}_R (1/\lambda) \sim \lambda/\sigma \in {\rm RV}^\infty_1$ as $\lambda \to\infty$ and claim (2) follows from Proposition~\ref{Thm.Karamata} and \ref{Thm.KaramataTauberian}. 
 To obtain claim (3), we apply (\ref{LapResolvent}) with $m=1$ to conclude that as $\lambda \to\infty$,
 \beqnn
 \widehat{\mathcal{I}}_R (1/\lambda)=\frac{1-{\it\hat\Phi}(1/\lambda)}{{\it\hat\Phi}(1/\lambda)} \sim \frac{1}{ {\it\hat\Phi}(1/\lambda) } .
 \eeqnn
An applications of Proposition~\ref{LapTheta} shows that $\widehat{\mathcal{I}}_R (1/\lambda)   \in {\rm RV}^\infty_{\alpha}$  if and only if ${\it\hat\Phi}(1/\cdot)\in{\rm RV}_{-\alpha}^\infty$ - equivalently, if and only if  ${\it \Phi}\in{\rm RV}_{-\alpha}^\infty$. Thus,  claim (3) follows  from Proposition~\ref{Thm.Karamata} and \ref{Thm.KaramataTauberian}.
 \qed

   \section{Subcritical Hawkes processes}\label{Sec.SubHP}
 \setcounter{equation}{0}

 In this section, we prove our FLTs for subcritical Hawkes processes $(m<1)$ given in Theorem~\ref{MainThm.02}. By Remark~\ref{rem:mu},  we know that 
 \[
 	 \mathbf{E}\Big[\frac{N(nt)}{n}\Big] = \mu_0\left( t + \frac 1 n \int_0^{nt} \int_0^s R(r) dr ds \right).
 \]
 Moreover, since $m = \|\phi\|_{L^1}$, the resolvent equation $R = \phi + R \ast \phi$ yields that
 \[
 	\int_0^\infty R(t)dt = \frac{1}{1-m}-1.
 \]
 This shows that 
 \beqnn
 	\mathbf{E}\Big[\frac{N(nt)}{n}\Big] 
 	\ar=\ar \mu_0 \left( t + \frac{1}{n} \int_0^{nt} \Big( \frac{m}{1-m} - \int_s^\infty R(r)dr \Big)ds \right)  \cr
	\ar=\ar \mu_0\left( \frac{1}{1-m} t - \frac 1 n \int_0^{nt}  \int_s^\infty R(r)dr ds \right). 
 \eeqnn
 Hence
  \beqnn
  n^{1/2}\Big( \frac{N(nt)}{n} -\frac{\mu_0}{1-m} \cdot t \Big) \ar=\ar n^{1/2}\Big(\frac{N(nt)}{n} - \mathbf{E}\Big[\frac{N(nt)}{n}\Big]\Big)  -  \frac{\mu_0}{n^{1/2}} \int_0^{nt}ds\int_s^\infty R(r)dr
  \eeqnn
and it follows from the FCLT for subcritical Hawkes processes recalled in Theorem~\ref{MainThm.01}(2) that the weak limit of the sequence of normalized processes depends strongly on the tail behavior of the integral function $$\left\{ \int_0^{t}ds\int_s^\infty R(r)dr : t \geq 0 \right\}.$$ The following two results clarify the tail behavior.  
 
%

  \begin{lemma}\label{MainThm02.Prop.01}
  	If $m<1$, ${\it\Psi}_*=\infty$ and ${\it\Phi}\in {\rm RV}_{-\alpha}^\infty $ with $\alpha\in[0,1/2]$, then as $t\to\infty$,
  	\beqnn
  	\int_0^{t}ds\int_s^\infty R(r)dr \sim  \frac{ t{\it\Phi}(t)}{ (1-\alpha)(1-m)^2}\in {\rm RV}_{1-\alpha}^\infty.
  	\eeqnn
  \end{lemma}
  \proof Using integration by parts and then  (\ref{LapResolvent}), we have for each $\lambda >0$,
  \beqlb\label{LapR}
  \int_0^\infty\frac{e^{-t/\lambda}}{\lambda} \int_0^tds  \int_s^\infty  R(r)dr dt 
  \ar=\ar \lambda  \int_0^\infty \big(1-e^{-t/\lambda}\big) R(t)dt
  = \frac{\lambda}{1-m}\cdot  \frac{{\it\hat\Phi}(1/\lambda)}{1-m+{\it\hat\Phi}(1/\lambda)}. \qquad
  \eeqlb
  By Proposition~\ref{LapTheta}, we have that ${\it\hat\Phi}(1/\lambda) \sim \Gamma(1-\alpha){\it\Phi}(\lambda)$ and
  \beqlb \label{eqn4.23}
  \int_0^\infty\frac{1}{\lambda}e^{-t/\lambda} \int_0^tds  \int_s^\infty  R(r)dr dt
  \sim \frac{\Gamma(1-\alpha)}{(1-m)^2}
  \cdot \lambda  {\it\Phi}(\lambda) \in {\rm RV}_{1-\alpha}^{\infty},
  \eeqlb
  as $\lambda \to \infty$. 
 The desired result now follows from Proposition~\ref{Thm.KaramataTauberian}(2).
  \qed

  \begin{proposition}\label{MainThm02.Prop.02}
  	If $m<1$ and  ${\it\Psi}_* \in [0, \infty)$, then as $t\to\infty$,
  	\beqnn
  	\int_0^{t}ds\int_s^\infty R(r)dr \sim  \frac{{\it\Psi}_*\cdot t^{1/2}}{(1-m)^{2}} \in {\rm RV}_{1/2}^{\infty}.
  	\eeqnn
  \end{proposition}
  \proof Let us first consider the case ${\it\Psi}_*>0$. In this case, we have that 
  \[
  	{\it\Psi}_1(t)\sim {\it\Psi}_* \cdot t^{1/2} \in  {\rm RV}_{1/2}^\infty.
\]
  From (\ref{eqn.3.11}) with $k=1$ and Proposition~\ref{Thm.Karamata}, it hence follows that $${\it\Phi}(t) \sim {\it\Psi}_*\cdot t^{-1/2}\in  {\rm RV}_{-1/2}^\infty.$$
  Similarly as in the proof of Proposition~\ref{MainThm02.Prop.01} with $\alpha=1/2$, the Laplace-Stieltjes transform of our integral function at $1/\lambda$ satisfies 
  \beqnn
  \int_0^\infty\frac{1}{\lambda}e^{-t/\lambda} \int_0^tds  \int_s^\infty  R(r)dr dt
  \sim \frac{\Gamma(1/2)}{(1-m)^2}
  \cdot \lambda  {\it\Phi}(\lambda)
  \sim   \frac{{\it\Psi}_*\cdot\Gamma(1/2)}{(1-m)^2}
  \cdot \lambda^{1/2}.
  \eeqnn
 Since our integral function is monotone, the desired result follows from  Proposition~\ref{Thm.KaramataTauberian}(2).
  
  Let us now consider the case of ${\it\Psi}_*=0$. In this case, we introduce,
for any $\delta>0$ and all $t>0$, the functions 
  \beqnn
  \phi_\delta(t):= \phi(t)+ \delta \cdot (t^2+t^3)^{-1/2} \quad \mbox{and} \quad
  {\it\Psi}(\delta ,t) :=\int_0^t s \phi_\delta(s)ds  = {\it\Psi}_1(t) + 2\delta \cdot (1+t)^{1/2}.
  \eeqnn
  Then ${\it\Psi}(\delta ,t) \sim 2\delta\cdot t^{1/2}$. 
  Let $R_\delta$ be the resolvent of $\phi_\delta$. Then $R(t)\leq  R_\delta(t)$ for any $t\geq 0$.
The desired result now follows by applying the previous result with ${\it\Psi}_*=2\delta$ as 
  \beqnn
  \int_0^{t}ds\int_s^\infty R(r)dr \leq \int_0^{t}ds\int_s^\infty R_\delta (r)dr  \sim \frac{2\delta}{(1-m)^2}\cdot t^{1/2}.
  \eeqnn
  \qed

   \medskip
   {\it Proof of Theorem~\ref{MainThm.02}.}
  If ${\it\Psi}_*\in [0,\infty)$, then it follows from Theorem~\ref{MainThm.01}(2) and Proposition~\ref{MainThm02.Prop.02} that the processes
  \beqnn
  n^{1/2}\Big( \frac{N(nt)}{n} -\frac{\mu_0}{1-\rho} \cdot t \Big) \ar=\ar n^{1/2}\Big(\frac{N(nt)}{n} - \mathbf{E}\Big[\frac{N(nt)}{n}\Big]\Big)  -  \frac{\mu_0}{n^{1/2}} \int_0^{nt}ds\int_s^\infty R(r)dr ,
  \eeqnn
converge weakly to $$\sqrt{\mu_0(1-m)^{-3} } B(t) - \mu_0 {\it\Psi}_* (1-m)^{-2}\cdot \sqrt{t}$$ in $\mathbf{D}([0,\infty),\mathbb{R})$. If ${\it\Psi}_*=\infty$,  then it follows from (\ref{eqn.3.11}) with $k=1$ that $n^{1/2}{\it\Phi}(n)\geq  n^{-1/2} {\it\Psi}_1(n)\to \infty$  as $n\to\infty$ and so it follows from
Theorem~\ref{MainThm.01}(2) that
  \beqnn
  \lim_{n\to\infty} \sup_{t\in[0,T]} \frac{1}{{\it\Phi}(n)} \Big|\frac{N(nt)}{n} - \mathbf{E}\Big[\frac{N(nt)}{n}\Big]\Big| = 0\quad \mbox{in} \quad \mathbf{P},
  \eeqnn
  for any $T\geq 0$. 
  Moreover, from Proposition~\ref{MainThm02.Prop.01} and~\ref{Thm.UniConver} we also have for any $T>0$ that
  \beqnn
  \lim_{n\to\infty} \sup_{t\in[0,T]}\Big|\frac{\mu_0}{n{\it\Phi}(n)} \int_0^{nt}ds\int_s^\infty R(r)dr   -\frac{\mu_0\cdot t^{1-\alpha}}{(1-\alpha)(1-m)^{2}} \Big| =0. 
  \eeqnn
Combining the last two results proves the desired result. 
  \qed

   \section{Weakly critical Hawkes processes}
 \label{Sec.WeaklyHP}
 \setcounter{equation}{0}

 In this section we prove our FLTs for weakly critical Hawkes processes $(m=1, \sigma < \infty)$ along with the representation of the Fourier-Laplace and the speed of convergence of the rescaled processes to the limiting one. 
 For convenience, we first introduce some notation that will be used throughout this section. For each $\beta> 0$ and $n\geq 1$,  we define
 \beqnn
 \phi^{(n)}_\beta(t):=e^{-\beta t}\phi(nt),\quad  R^{(n)}_\beta(t):=e^{-\beta t} R^{(n)}(t) = e^{-\beta t} R(nt),\quad R^*_\beta(t):=\frac{1}{\sigma}e^{-\beta t},\quad t\geq 0.
 \eeqnn
 By (\ref{Resolvent}), it is easy to see that 
 \beqlb\label{eqn.30050}
 R^{(n)}_\beta=\phi^{(n)}_\beta+ n\cdot R^{(n)}_\beta*\phi^{(n)}_\beta. 
 \eeqlb
 Integrating both sides of this equation over $\mathbb{R}_+$, yields that
 \beqlb\label{eqn.3.61}
 \| R^{(n)}_\beta\|_{L^1}=\|\phi^{(n)}_\beta\|_{L^1}+ n\|\phi^{(n)}_\beta\|_{L^1}\cdot \|R^{(n)}_\beta\|_{L^1}
 \quad \mbox{and}\quad
 \| R^{(n)}_\beta\|_{L^1} = \frac{n\|\phi^{(n)}_\beta\|_{L^1}}{n\big(1-n\|\phi^{(n)}_\beta\|_{L^1} \big)}.
 \eeqlb
Since $m=\|\phi\|_{L^1}=1$ and $\sigma={\it \Psi}_1(\infty)<\infty$, we have
$n\|\phi_\beta^{(n)}\|_{L^1} <1$ for any $n\geq 1$. 
Moreover, 
\beqlb\label{eqn.3.24}
 n\big(1-n\|\phi_\beta^{(n)}\|_{L^1}\big)=\int_0^\infty n(1-e^{-\beta t/n})\phi(t)dt \to \beta \sigma>0
 \quad \mbox{and}\quad \| R^{(n)}_\beta\|_{L^1}  \to \frac{1}{\beta\sigma},
 \eeqlb
 as $n\to\infty$. Hence there exists a constant $n_0\geq 1$ such that for any $n\geq n_0$,
 \beqlb\label{eqn.3.23}
 n\big(1-n\|\phi_\beta^{(n)}\|_{L^1}\big)\geq\frac{\beta\sigma }{2}
 \quad \mbox{and}\quad
 \| R^{(n)}_\beta\|_{L^1} \leq \frac{2}{\beta\sigma}.
 \eeqlb

 \subsection{Proofs of Theorem~\ref{MainThm.04} and \ref{Thm.FourLapCIR}}\label{sec:MainThm04}

 
 To obtain the non-positivity and the upper estimates on the Fourier-Laplace exponents of $\Lambda^{(n)}$ and $\Lambda^*$, it will be convenient to first analyze the convergence of families rescaled Hawkes processes with general immigration densities. Specifically, for each $n\geq 1$, we denote by $N_n$ a Hawkes process whose density process $\Lambda_n$ is of the form (\ref{HawkesDensity}) with immigration intensity
 \beqnn
 \mu= \mu_n \in L^1_{\rm loc}(\mathbb{R}_+;\mathbb{R}_+)
 \eeqnn
 and consider the rescaled processes 
  \[
  	X^{(n)}_H:= \left\{ \left( \Lambda^{(n)}_H(t), \mathcal{I}_{\Lambda^{(n)}_H}(t), N_H^{(n)}(t),\widetilde{N}_H^{(n)}(t) \right):t\geq 0 \right\}  
\]	
where
 \begin{equation} \label{tildeNH}
 \Lambda^{(n)}_H(t):= \frac{\Lambda_{n}(nt)}{n}, \quad  
 \mathcal{I}_{\Lambda^{(n)}_H}(t) : = \frac{\mathcal{I}_{\Lambda_{n}}(nt)}{n^2} ,\quad 
 N_H^{(n)}(t):= \frac{N_{n}(nt)}{n^2},
 \quad  
 \widetilde{N}_H^{(n)}(t):= \frac{\widetilde{N}_{n}(nt)}{n}. 
  \end{equation}
 If we let 
  \beqnn
   H_{\mu_n}:= \mu_n+  \mu_n*R  \in L^1_{\rm loc}(\mathbb{R}_+;\mathbb{R}_+)  \quad \mbox{and} \quad
   H^{(n)}(t):= \frac{H_{\mu_n} (nt)}{n}, ~~ t \geq 0,
 \eeqnn
then $H^{(n)}$ is the expected value of the rescaled process $N_n$. In view of the martingale representation result (\ref{MartRep}), the process $\Lambda^{(n)}_H$ can be written as
 \beqlb\label{eqn.3006}
 \Lambda^{(n)}_H(t)=  H^{(n)}(t) + R^{(n)}*d\widetilde{N}_H^{(n)}(t)=   H^{(n)}(t) + \int_0^t R^{(n)}(t-s)\widetilde{N}_H^{(n)}(ds),\quad t\geq 0.
 \eeqlb
 
 Let us now assume that the process $H^{(n)}$ converges to a non-negative and non-decreasing c\'adl\'ag function $H^*$ on $\mathbb{R}_+$. 
 Since $\Lambda^{(n)}_H$ is the quadratic variation process of the martingale  $\widetilde{N}_H^{(n)}$, Proposition~\ref{Thm.AsymR}(2) suggests that in this case the intensity process $\Lambda^{(n)}_H$ converges to the unique solution $\Lambda^{*}_H$  of the SDE\,\footnote{ the existence and uniqueness of strong non-negative solutions to (\ref{CIR}) is standard; see Theorem~1 in \cite{YamadaWatanabe1971}.}
 \beqlb\label{CIR}
 \Lambda^*_H(t)= H^*(t) + \int_0^t \frac{1}{\sigma} \sqrt{\Lambda^*_H(s)}dB(s),\quad t\geq 0.
 \eeqlb
 
 We do not expect the convergence of the rescaled densities to guarantee the convergence of the rescaled Hawkes processes. Instead, we are going to prove that if $\mathcal{I}_{H^{(n)}}\to \mathcal{I}_{H^*}$ locally uniformly, then
\[
 	\Big(\Lambda^{(n)}_H, \mathcal{I}_{\Lambda^{(n)}_H},N_H^{(n)},\widetilde{N}_H^{(n)}\Big)\to \Big(\Lambda^*_H, \mathcal{I}_{ \Lambda^*_H}, \mathcal{I}_{ \Lambda^*_H},\sigma(\Lambda^*_H-H^*)\Big)
\]
 weakly in $M(\mathbb{R}_+;\mathbb{R}_+)\times\mathbf{D}(\mathbb{R}_+;\mathbb{R}^3) $  as $n\to\infty$. From this, we will then deduce that Theorem~\ref{MainThm.04} holds.

\subsubsection{Proof of Theorem~\ref{MainThm.04}}

 We start with the following two lemmas that provide an explicit exponential-affine representation of the Fourier-Laplace functionals of $\Lambda^*_H$ and $\Lambda^{(n)}_H$.

 \begin{lemma}\label{LaplaceCIR}
 	For each $w\in M(\mathbb{R}_+;\mathbb{C}_-)$ and $g\in L^\infty_{\rm loc}(\mathbb{R}_+;\mathtt{i}\mathbb{R})$, 
    we have that $V^*\in L^\infty_{\rm loc}(\mathbb{R}_+;\mathbb{C}_-)$ and 
 	\beqlb\label{eqn.3010}
 	\mathbf{E}\Big[\exp\Big\{  \Lambda^*_H* dw(T) + \int_0^T g(T-t)\sqrt{\Lambda^*_H(t)}dB(t)\Big\}\Big]= \exp\Big\{  \sigma \int_0^T H^*(T-t) dV^*(t) \Big\},\quad T\geq 0. 
 	\eeqlb
 	Moreover, if $w\in M(\mathbb{R}_+;\mathbb{R}_-)$ and $g\equiv0$, then $0\geq V^* \geq \sigma^{-1}\cdot \mathcal{I}_w$.
 \end{lemma}
 \proof 
 Since $V^*$ is locally of bounded variation,  for each $T\geq 0$ the following random variable is well defined: 
 \beqnn
 \Lambda^*_T:=  \sigma \int_0^T\Lambda^*_H (T-t) dV^*(t)= \sigma \int_0^T \Lambda^*_H(t) dV^*(T-t).
 \eeqnn
 Taking expectations on both sides of (\ref{CIR}), we see that $\mathbf{E}[\Lambda^*_H(t)]= H^*(t)$ for $t\geq 0$ and hence
 \beqnn
 \mathbf{E}[\Lambda^*_T]=  \sigma \int_0^T H^*(T-t) dV^*(t).
 \eeqnn
 
Let us now consider the Doob martingale $\big\{\Lambda^*_T(t):=\mathbf{E}\big[\Lambda^*_T\big| \mathscr{F}_t\big]:t\in[0,T] \big\}$.
 By the stochastic Fubini theorem,
 \beqnn
 \Lambda^*_T(t)
 \ar=\ar \mathbf{E}[\Lambda^*_T] + \int_0^T  \int_0^{r\wedge t} \sqrt{\Lambda^*_H(s)}dB(s)   dV^*(T-r)
 =\mathbf{E}[\Lambda^*_T] + \int_0^t V^{*}(T-s)  \sqrt{\Lambda^*_H(s)}dB(s).
 \eeqnn
 Moreover, we define the semi-martingale $\{\xi_T(t):t\in[0,T]\}$ by
 \beqnn
 \xi_T(t)\ar:=\ar \Lambda^*_T(t) + \int_0^t g(T-s)\sqrt{\Lambda^*_H(s)}dB(s) + \int_0^t \Lambda^*_H(s)w(d(T-s)) - \int_0^t \sigma \Lambda^*_H(s)dV^*(T-s)\cr
 \ar=\ar  \mathbf{E}[\Lambda^*_T] + \int_0^t  (V^{*}+g)(T-s)  \sqrt{\Lambda^*_H(s)}dB(s) + \int_0^t \Lambda^*_H(s)w(d(T-s))  - \int_0^t \sigma \Lambda^*_H(s)dV^*(T-s).
 \eeqnn
 Applying It\^o's formula to $\exp\{ \xi_T(t) \}$ and then using the equality (\ref{CBI.Riccati}),  we see that
 \beqlb\label{eqn.3009}
 \exp\big\{\xi_T(t)\big\}
 \ar=\ar \exp\big\{ \mathbf{E}[\Lambda^*_T] \big\} +\int_0^t  e^{\xi_T(s)} V^*(T-s)  \sqrt{\Lambda^*_H(s)}dB(s),\quad t\in[0,T],
 \eeqlb
 which is a local martingale. Using the same arguments given in the proof of Lemma~\ref{Lemma.FourLapGHP} one can actually verify that the local martingale is indeed a true martingale.
 Since 
 $$\xi_T(T)= \Lambda^*_H* dw(T) + \int_0^T g(T-t)\sqrt{\Lambda^*_H(t)}dB(t) \in \mathbb{C}_-,$$  
  the equation (\ref{eqn.3010}) follows by taking expectations on both sides of (\ref{eqn.3009}) with $t=T$.   In particular, if $H^*(t):= \mathbf{1}_{\{t\geq 0\}}$, then 
  \[
  	\exp\big\{ \sigma V^*(T) \big\}=\mathbf{E} \left[\exp \left\{   \Lambda^*_H* dw(T)+ \int_0^T g(T-t)\sqrt{\Lambda^*_H(t)}dB(t) \right\} \right].
\]
 Since $\Lambda^*_H * dw(T) \in\mathbb{C}_- $ and $  \int_0^T g(T-t)\sqrt{\Lambda^*_H(t)}dB(t) \in {\tt i}\mathbb{R}$, we see that
 \beqnn
 \exp\big\{\sigma\cdot {\tt Re}V^*(T) \big\} = |\exp\big\{ \sigma V^*(T) \big\}|
 \ar\leq\ar  \mathbf{E}\big[\exp\big\{ {\tt Re} (\Lambda^*_H* dw(T))\big\}\big] \leq 1 .
 \eeqnn
 Thus,  ${\tt Re}V^*(T)\leq 0$ and hence $V^*\in L^\infty_{\rm loc}(\mathbb{R}_+;\mathbb{C}_-)$.
 In addition, if $w \in M(\mathbb{R}_+;\mathbb{R}_-)$ and $g\equiv0$, then $$V^* \in L^\infty_{\rm loc}(\mathbb{R}_+;\mathbb{R}_-)$$ and by Jensen's inequality, 
 \beqnn
 \exp\big\{\sigma \cdot  V^*(T) \big\}= \mathbf{E}\big[\exp\big\{  \Lambda^*_H* dw(T)\big\}\big] \geq \exp\big\{   \mathbf{E}[\Lambda^*_H]* dw(T)\big\}  = \exp\big\{ H^** dw(T) \big\} = \exp \big\{\mathcal{I}_w(T) \big\}.
 \eeqnn
This shows that  $0\geq V^*(T)\geq \sigma^{-1}\cdot \mathcal{I}_w(T)$.
 \qed

  \begin{lemma}\label{Lemma.FourLapGHP01}
 	For each $n\geq 1$, $w\in M(\mathbb{R}_+;\mathbb{C}_-)$ and $g \in L^\infty_{\rm loc}(\mathbb{R}_+;\mathtt{i}\mathbb{R})$, we have  that
 	\beqlb\label{eqn.FouLapGHP01}
 	\mathbf{E}\big[\exp\big\{  \Lambda^{(n)}_H* dw(T) +g* d\widetilde{N}^{(n)}_H(T) \big\}\big]
 		= \exp \big\{ H^{(n)}* dw(T)+  H^{(n)} * W^{(n)}(T) \big\}, \quad T\geq 0,
 	\eeqlb
 	where $$W^{(n)}:=n^2\big(e^{(V^{(n)} +g)/n}-1-(V^{(n)} +g)/n \big)$$ and $V^{(n)} \in L^\infty_{\rm loc}(\mathbb{R}_+;\mathbb{C}_-)$ is  the unique solution to the following two equivalent Volterra integral equations:
 	\beqlb
 	V^{(n)} \ar=\ar R^{(n)}* dw + n^2 \big(e^{(V^{(n)}+g)/n} -1-  (V^{(n)}+g)/n \big)* R^{(n)},\label{Laplace03} \\
 	V^{(n)}\ar=\ar \phi^{(n)}* dw + n^2 \big(e^{ (V^{(n)}+g)/n } -1 -g/n\big)* \phi^{(n)}. \label{Laplace04}
 	\eeqlb
  In addition,  the following hold:
  \begin{enumerate}
  	\item[(1)] If $w \in M(\mathbb{R}_+;\mathbb{R}_-)$ and $g\equiv 0$, then $V^{(n)} \in L^\infty_{\rm loc}(\mathbb{R}_+;\mathbb{R}_-)$ and $0\geq V^{(n)} \geq R^{(n)}* dw$.
  	
  	\item[(2)] For each $T\geq 0$, there exist an integer $n_0\geq 1$ and a constant $C>0$ such that for any $g \in L^\infty_{\rm loc}(\mathbb{R}_+;\mathtt{i}\mathbb{R})$ and any absolutely continuous measure $w(dt)=f(t)dt$ with density $f\in L^\infty_{\rm loc}(\mathbb{R}_+;\mathbb{C}_-)$, 
  	\beqlb \label{eqn.3011}
    \sup_{n\geq n_0}\|V^{(n)}\|_{L^\infty_T}\leq C  \big( \|f\|_{L^\infty_T}+\|g\|_{L^\infty_T}^2\big) 
  \quad\mbox{and}\quad 
 	  \sup_{n\geq n_0} \|W^{(n)}\|_{L^\infty_T}\leq C  \big( \|f\|_{L^\infty_T}+\|g\|_{L^\infty_T}+\|g\|_{L^\infty_T}^2\big)^2 .  
 	  \ \ 
  	\eeqlb
  	\end{enumerate}
 \end{lemma}
 \proof For each $n\geq 1$ and $t\geq 0$, let 
 \[
 	g_n(t):=\frac{g(t/n )}{n} \quad \mbox{and} \quad 
	w_n(dt):= \frac 1 n w(n^{-1}\cdot dt).
\] 
 Let $V_n \in L^\infty_{\rm loc}(\mathbb{R}_+;\mathbb{C}_-)$ be the unique solution of the Volterra integral equation (\ref{VolRiccati}) with 
 \[
 	\nu(dt)=w_n(dt)-g_n(t)dt \quad \mbox{and} \quad f=g_n.
\] 
 For $t\geq 0$, let
 \beqnn
 W_n(t):= e^{(V_n+g_n)(t)}-1-(V_n+g_n)(t)
 \quad \mbox{and}\quad 
 V^{(n)}(t):=n\cdot V_n(nt).
 \eeqnn 
 By a change of variables, it follows that $V^{(n)}\in L^\infty_{\rm loc}(\mathbb{R}_+;\mathbb{C}_-)$ is the unique solution to the equivalent equations (\ref{Laplace03})-(\ref{Laplace04}), and 
 \beqnn
 W^{(n)}(t)=n^2 \cdot W_n(nt),\quad  \Lambda^{(n)}_H* dw(t)= \Lambda_{H,n}* dw_n(nt),\quad g* d\widetilde{N}^{(n)}_H(t)= g_n* dN_n(nt) - \Lambda_{H,n}*g_n(nt). 
 \eeqnn
 Claim (1) follows directly from Theorem~\ref{Thm201}, and the equation (\ref{eqn.FouLapGHP01}) follows from Proposition~\ref{Lemma.FourLapGHP}: 
 \beqnn
 \mathbf{E}\big[\exp\big\{  \Lambda^{(n)}_H* dw(T) +g* d\widetilde{N}^{(n)}_H(T) \big\}\big]
 \ar=\ar \exp\big\{ H_{\mu_n}*(dw_n-g_n)(nT)+H_{\mu_n}* W_n (nT) \big\}  \\ \vspace{10pt}
 \ar=\ar \exp\big\{ H^{(n)}* dw(T)+H^{(n)}* W^{(n)} (T) \big\}.
 \eeqnn
 
To prove Claim (2) we rewrite the equation  (\ref{Laplace04}) as
 \beqnn
 V^{(n)}\ar=\ar \phi^{(n)}*\big( f+ e^{g/n}\cdot n^2 (e^{V^{(n)}/n} -1) + n^2  (e^{g/n} -1 -g/n ) \big).
 \eeqnn
 Using the inequalities $|e^{x}-1|\leq |x|$ and $|e^x-1-x|\leq |x|^2$ for any $x\in \mathbb{C}_-$, we see that
 \beqnn
 \begin{split}
 |V^{(n)}| & \leq \phi^{(n)}*\Big(|f|+ n^2  \big|e^{V^{(n)}/n} -1\big| + n^2 \big|e^{g/n} -1 -g/n\big|\Big) 
 \leq   \phi^{(n)}  * \big( n |V^{(n)}|  + |f| +|g|^2 \big).
 \end{split}
 \eeqnn
 For $\beta>0$, let  $V^{(n)}_{\beta}(t):= e^{-\beta t}V^{(n)}(t)$. The preceding inequalities imply that
 \beqlb\label{eqn.3.22}
 |V^{(n)}_{\beta}|  \leq    \phi^{(n)}_{\beta}  * \big( n |V^{(n)}_{\beta}|  +  |f|+|g|^2 \big) . 
 \eeqlb
 It now follows from Young's convolution inequality that
 \beqnn
 \|V^{(n)}_{\beta}\|_{L^\infty_T} \leq \big( n \|V^{(n)}_{\beta}\|_{L^\infty_T} + \|f \|_{L^\infty_T}+\|g \|^2_{L^\infty_T}\big) \cdot \|\phi_\beta^{(n)}\|_{L^1}.
 \eeqnn
 In view of the inequalities (\ref{eqn.3.23}), there exists a constant $n_0\geq 1$ such that for any $n\geq n_0$,
 \beqlb\label{eqn.3.23.01}
 \|V^{(n)}_{\beta}\|_{L^\infty_T} \leq \frac{ n\|\phi_\beta^{(n)}\|_{L^1}}{ n\big(1-n\|\phi_\beta^{(n)}\|_{L^1}\big)}\cdot \big(\|f\|_{L^\infty_T}+ \|g \|_{L^\infty_T}^2 \big)  \leq \frac{2}{\beta\sigma} \cdot \big(\|f\|_{L^\infty_T}+ \|g \|_{L^\infty_T}^2 \big), 
 \eeqlb
from which we deduce that there exists a constant $C>0$ that does not dependent on $f$ and $g$ such that
\[
	\sup_{n\geq 1} \|V^{(n)}\|_{L^\infty_T} \leq \frac{C}{\beta\sigma} \cdot e^{\beta T}  \cdot  \big(\|f\|_{L^\infty_T}+ \|g \|_{L^\infty_T}^2 \big).  
\]	
To  prove the second inequality in (\ref{eqn.3011}) we first use the inequality $|e^x-1-x|\leq x^2$ for any $x\in\mathbb{C}_-$ and then the Minkowski inequality, to obtain a constant $C>0$ independent of $f$ and $g$ such that
 \beqnn
 \|W^{(n)}\|_{L^\infty_T} 
 \leq   \|V^{(n)} +g\|_{L^\infty_T}^2
 \leq \big(\|V^{(n)} \|_{L^\infty_T} + \|g\|_{L^\infty_T}\big)^2 
 \leq C\big(\|f\|_{L^\infty_T}+ \|g\|_{L^\infty_T} +\|g\|^2_{L^\infty_T}\big)^2.
 \eeqnn
 \qed

The preceding results allow us to establish the convergence of the solution to our Volterra equations  (\ref{Laplace03}) to the uinque solution of the equation (\ref{CBI.Riccati}) along with the rate of convergence assuming that $g \equiv 0$ and that the measure $w\in M(\mathbb{R}_+;\mathbb{R}_-)$ has a density. 

  \begin{proposition}\label{Prop.ConV}
 	Let $V^{(n)}$ and $V^*$  be  solutions to  (\ref{Laplace03}) and (\ref{CBI.Riccati}) respectively with 
 \beqnn
	g \equiv 0 \quad \mbox{and} \quad w(dt)=f(t)dt.
 \eeqnn
 	For each $T\geq0$ and $f\in C^1(\mathbb{R}_+;\mathbb{R}_-)$, we have $	\| V^{(n)} - V^*  \|_{L^\infty_T} \to 0 $ as $n\to\infty$.
 
 \end{proposition}
 \proof Recalling that $R^*(t) \equiv 1/\sigma$ and using (\ref{Laplace03}) and (\ref{CBI.Riccati}), we see that 
 \beqlb\label{eqn.50001}
V^{(n)} - V^*=\sum_{i=1}^3 \varepsilon_i^{(n)} 
\quad \mbox{with}\quad  
 \begin{array}{l}
 \displaystyle{\varepsilon_1^{(n)} := \Big(f+ \frac{ (V^*)^2}{2}\Big)*(R^{(n)}-R^*),}\vspace{7pt}\cr
  \displaystyle{ \varepsilon_2^{(n)}:=\frac{1}{2}\big( (V^{(n)})^2 - (V^*)^2\big)*R^{(n)}, }  \vspace{7pt} \cr
  \displaystyle{\varepsilon_3^{(n)} :=  n^2\Big(e^{V^{(n)}/n}-1- \frac{V^{(n)}}{n}-\frac{1}{2}\Big(\frac{V^{(n)}}{n}\Big)^2\Big) *R^{(n)} .}
 \end{array}
 \eeqlb
 For  a constant $\beta>0$ to be specified later, let
 \beqnn
 V^{(n)}_\beta(t):= e^{-\beta t}V^{(n)} (t),\quad
 V^*_\beta(t):= e^{-\beta t}V^*(t),\quad
 \varepsilon_{i,\beta}^{(n)}(t):= e^{-\beta t}\varepsilon_i^{(n)}(t),\quad t\geq 0,~~i=1,2,3.
 \eeqnn
Then,
 \beqnn
 	\|V^{(n)} - V^*\|_{L^\infty_T} \leq e^{\beta T} \|V^{(n)}_\beta - V^*_\beta\|_{L^\infty_T} \quad 
	\mbox{and} \quad 
	\|V^{(n)}_\beta - V^*_\beta\|_{L^\infty_T}\leq  \sum_{i=1}^3\|\varepsilon_{i,\beta}^{(n)}\|_{L^\infty_T}.
  \eeqnn
  Hence it suffice to prove $	\|V^{(n)}_\beta - V^*_\beta\|_{L^\infty_T} \to0 $. 
We consider the three quantities  $\varepsilon_{i.\beta}^{(n)}$ (i=1,2,3) separately, starting with the second one. Since 
 \beqnn
	\varepsilon_{2,\beta}^{(n)} =  \frac{1}{2} \big(( V^{(n)}+V^* )(V^{(n)}_\beta - V^*_\beta)\big)*R^{(n)}_\beta,
 \eeqnn
 it follows from Lemma~\ref{LaplaceCIR} and Lemma \ref{Lemma.FourLapGHP01}(2) that there exits a constant $C_0>0$ depending only on $T$ such that for any $n\geq 1$ and $f\in C^1(\mathbb{R}_+;\mathbb{R}_-)$, 
 \beqlb\label{eqn3.23}
 \|\varepsilon_{2,\beta}^{(n)}\|_{L^\infty_T}
 \ar\leq\ar  \frac{1}{2} \big(\| V^{(n)}\|_{L^\infty_T} + \| V^* \|_{L^\infty_T}\big)\cdot \|R^{(n)}_\beta\|_{L^1} \cdot \|V^{(n)}_\beta - V^*_\beta\|_{L^\infty_T}\cr
 \ar\leq\ar  \frac{C_0}{2} \|f\|_{L^\infty_T}\cdot  \|R^{(n)}_\beta\|_{L^1} \cdot \|V^{(n)}_\beta - V^*_\beta\|_{L^\infty_T}.
 \eeqlb
  From  (\ref{eqn.3.23}), we have that
 $ \| R^{(n)}_\beta\|_{L^1} \leq 2/(\beta\sigma)$ for $n\geq n_0$. Taking this back into (\ref{eqn3.23}) and then choosing 
 \beqlb\label{eqn.beta}
 \beta = \frac{2C_0}{\sigma}\cdot \|f\|_{L^\infty_T}, 
 \eeqlb
 we see that $\|\varepsilon_{2,\beta}^{(n)}\|_{L^\infty_T} \leq \|V^{(n)}_\beta - V^*_\beta\|_{L^\infty_T}/2$ and hence that
 \beqnn
 \|V^{(n)}_\beta - V^*_\beta\|_{L^\infty_T}\leq   \|\varepsilon_{1,\beta}^{(n)}\|_{L^\infty_T} + \|\varepsilon_{3,\beta}^{(n)}\|_{L^\infty_T} + \frac{1}{2}\|V^{(n)}_\beta - V^*_\beta\|_{L^\infty_T}. 
 \eeqnn
As a result,
\beqnn
\|V^{(n)}_\beta - V^*_\beta\|_{L^\infty_T} \leq 2 \big(  \|\varepsilon_{1}^{(n)}\|_{L^\infty_T} + \|\varepsilon_{3}^{(n)}\|_{L^\infty_T}  \big) 
\eeqnn
 and so it remains to  prove that both $\|\varepsilon_{1}^{(n)}\|_{L^\infty_T} $ and $ \|\varepsilon_{3}^{(n)}\|_{L^\infty_T}$ vanish as $n\to\infty$. 

 Using integration by parts on $\varepsilon_{1}^{(n)}$ yields that
 \beqnn
 \varepsilon_{1}^{(n)}= f(0)\cdot (\mathcal{I}_{R^{(n)}} -\mathcal{I}_{R^*})+ (f'+V^*V^{*'})* (\mathcal{I}_{R^{(n)}} -\mathcal{I}_{R^*} ).
 \eeqnn
 By the second result in Lemma~\ref{LaplaceCIR} and (\ref{CBI.Riccati}), we have $\|V^*\|_{L^\infty_T} + \|V^{*'}\|_{L^\infty_T} <\infty$ and then 
 \beqnn
  \| \varepsilon_{1}^{(n)} \|_{L^\infty_T}
 \ar \leq\ar   |f(0)|\cdot \big\|\mathcal{I}_{R^{(n)}} -\mathcal{I}_{R^*} \big\|_{L^\infty_T}+ \|f'+V^*V^{*'}\|_{L^1_T}\cdot \big \|\mathcal{I}_{R^{(n)}} -\mathcal{I}_{R^*} \big\|_{L^\infty_T} \cr
 \ar\leq\ar \big(   \|f\|_{L^\infty_T} + \|f'\|_{L^\infty_T}\cdot T + \|V^*\|_{L^\infty_T} \cdot \|V^{*'}\|_{L^\infty_T}\cdot T  \big)\cdot \big \|\mathcal{I}_{R^{(n)}} -\mathcal{I}_{R^*} \big\|_{L^\infty_T} ,
 \eeqnn
 which goes to $0$ as $n\to\infty$, due to Corollary~\ref{ConvergenceR}. 
 
 Finally, applying the inequality $|e^{x}-1-x-x^2/2|\leq |x|^3$ for any $x\leq 0$ to $\varepsilon_3^{(n)}$, we have  that
 \beqnn
 |\varepsilon_3^{(n)}| \leq n^{-1} \cdot |V^{(n)}|^3*R^{(n)}
 \quad\mbox{and hence}\quad 
 \|\varepsilon_3^{(n)}\|_{L^\infty_T} \leq n^{-1}\cdot\|V^{(n)}\|_{L_T^\infty}^3 \cdot \mathcal{I}_{R^{(n)}}(T).
 \eeqnn
 An applications of  Lemma~\ref{Lemma.FourLapGHP01}(2) and Corollary~\ref{ConvergenceR} yields that $ \|\varepsilon_3^{(n)}\|_{L^\infty_T} \leq C/n \to 0 $ as $n\to\infty$. 
 \qed

The next lemma provides the desired sufficient condition for the weak convergence of the rescaled intensities and (compensated) Hawkes processes in terms of their integrated expected values.   

  \begin{lemma}\label{Lemma.ConverGHP}
 If $\mathcal{I}_{H^{(n)}}\to \mathcal{I}_{H^*}$ locally uniformly, then
\[
 	\big(\Lambda^{(n)}_H, \mathcal{I}_{\Lambda^{(n)}_H},N_H^{(n)},\widetilde{N}_H^{(n)}\big)\to \big(\Lambda^*_H, \mathcal{I}_{ \Lambda^*_H}, \mathcal{I}_{ \Lambda^*_H},\sigma(\Lambda^*_H-H^*)\big)
\]
 	weakly in $M(\mathbb{R}_+;\mathbb{R}_+)\times\mathbf{D}(\mathbb{R}_+;\mathbb{R}^3) $  as $n\to\infty$.

 \end{lemma}
 \proof   For each $f\in C^1(\mathbb{R}_+;\mathbb{R}_-)$, let $V^{(n)}$ and $V^*$  be  solutions to  (\ref{Laplace03}) and (\ref{CBI.Riccati}) respectively with $g \equiv 0$ and $w(dt)=f(t)dt$.
 From  Lemma~\ref{LaplaceCIR} and  \ref{Lemma.FourLapGHP01}, we have 
 \beqnn
  \mathbf{E}\big[\exp\big\{f*\Lambda^{(n)}_H(T)\big\}\big] 
  \ar=\ar \exp \big\{ H^{(n)}* f(T)+  H^{(n)} * \big( n^2 (e^{V^{(n)}/n}-1-V^{(n)}/n)\big)(T) \big\} , \cr
  \mathbf{E}\big[\exp\big\{f*\Lambda^*_H(T)\big\}\big] \ar=\ar \exp\Big\{ H^* *f(T) + \frac{1}{2} \cdot H^* * |V^*|^2(T) \Big\}. 
 \eeqnn
 By (\ref{eqn.3011}) and the inequality $|e^{x}-1-x-x^2/2|\leq |x|^3$ for any $x\leq 0$, we have as $n\to\infty$, 
 \beqnn
 n^2 \big(e^{V^{(n)}(t)/n}-1-V^{(n)}(t)/n\big) \sim   |V^{(n)}(t)|^2 /2, 
 \eeqnn
 which converges to $ |V^*(t)|^2/2$ uniformly in $t\in[0,T]$, due to Proposition~\ref{Prop.ConV}. 
From this and the assumption that $\mathcal{I}_{H^{(n)}}\to \mathcal{I}_{H^*}$ locally uniformly, for each $T\geq 0$ we have as $n\to\infty$, that  
 \beqnn
 \mathbf{E}\big[\exp\big\{f*\Lambda^{(n)}_H(T)\big\}\big]  \to \mathbf{E}\big[\exp\big\{f*\Lambda^*_H(T)\big\}\big],
 \eeqnn 
 and hence  $f*\Lambda^{(n)}_H(T)\to f*\Lambda^*_H(T)$ weakly. 
 Since $C^1(\mathbb{R}_+;\mathbb{R}_-)$ is dense in $C(\mathbb{R}_+;\mathbb{R}_-)$, it follows that $\Lambda^{(n)}_H \to \Lambda^*_H$ weakly in $M(\mathbb{R}_+;\mathbb{R}_+)$.
 By footnote~\ref{Footnote.03}, this is equivalent to $\mathcal{I}_{\Lambda^{(n)}_H}\to\mathcal{I}_{\Lambda_H^*}$ weakly in $\mathbf{C}(\mathbb{R}_+;\mathbb{R}_+) $.  

 Let us now consider the rescaled process $N^{(n)}_H$ and the corresponding compensated  process $\widetilde N^{(n)}_H$ as introduced in \eqref{tildeNH}. By definition, 
 \beqnn
 N^{(n)}_H= 
 \mathcal{I}_{\Lambda^{(n)}_H}+ n^{-1} \cdot \widetilde{N}^{(n)}_H. 
 \eeqnn
 Moreover, the compensated Hawkes process $\widetilde{N}^{(n)}_H$ is a locally square-integrable martingale  with (predictable) quadratic variation
 \beqnn
 [\widetilde{N}^{(n)}_H]_t= N^{(n)}_H(t),\quad  \langle \widetilde{N}^{(n)}_H\rangle_t =  
 \mathcal{I}_{\Lambda^{(n)}_H} (t) \quad \mbox{and} \quad
 \mathbf{E} \big[  [\widetilde{N}^{(n)}_H]_t \big]= \mathbf{E} \big[  \langle \widetilde{N}^{(n)}_H\rangle_t \big] = {\cal I}_{H^{(n)}}(t),
 \quad t\geq 0. 
 \eeqnn
 Hence, the Burkholder-Davis-Gundy inequality (e.g. see Theorem 26.12 in \cite[p.524]{Kallenberg2002}) yields a constant  $C>0$ that is independent of $n$, such that 
 \beqnn
 \mathbf{E}\Big[\sup_{t\in[0,T]}\big| \widetilde{N}_H^{(n)}(t)/n\big|^2\Big] \leq C \cdot  \mathbf{E} \big[  [\widetilde{N}^{(n)}_H]_T \big]/n^2 
 =   C\cdot {\cal I}_{H^{(n)}}(T) /n^2\to 0 \quad \mbox{as} \quad n \to \infty, 
 \eeqnn
 from which we conclude that $\sup_{t\in[0,T]}|\widetilde{N}_H^{(n)}(t)/n|\to0$ in $\mathbf{P}$. 
 In particular,  it follows that   
 \[
 	( \Lambda^{(n)}_H, \mathcal{I}_{\Lambda^{(n)}_H},N_H^{(n)} )\to  ( \Lambda^*_H,\mathcal{I}_{ \Lambda^*_H}, \mathcal{I}_{ \Lambda^*_H})
 \] 
 weakly in $M(\mathbb{R}_+;\mathbb{R}_+)\times\mathbf{D}(\mathbb{R}_+;\mathbb{R}^2)$. 
 Furthermore, an application of Theorem~3.11 in  \cite[p.473]{JacodShiryaev2003} yields that $\widetilde{N}^{(n)}_H$ converges weakly in $\mathbf{D}(\mathbb{R}_+;\mathbb{R})$ to a continuous Gaussian martingale $G^*$ with predictable quadratic variation $\mathcal{I}_{\Lambda^*_H}$.
 By Corollary~3.33(b) in \cite[p.353]{JacodShiryaev2003} it follows that 
 \[
 	( \Lambda^{(n)}_H, \mathcal{I}_{\Lambda^{(n)}_H},N_H^{(n)},\widetilde{N}_H^{(n)})\to  ( \Lambda^*_H,\mathcal{I}_{ \Lambda^*_H}, \mathcal{I}_{ \Lambda^*_H},G^*) 
 \]	
 weakly in $M(\mathbb{R}_+;\mathbb{R}_+)\times\mathbf{D}(\mathbb{R}_+;\mathbb{R}^3)$. 
 By Skorokhod's representation theorem and the continuity of limit process, we may actually assume that this convergence holds almost surely and locally uniformly.    

 It remains to characterize the limit process $G^*$. 
 Integrating both sides of (\ref{eqn.3006})  and then using Fubini's theorem, we have  for $t\geq 0$ that
 \beqlb\label{eqn.3007}
 \begin{split}
 \mathcal{I}_{\Lambda^{(n)}_H}(t) & = \mathcal{I}_{H^{(n)}}(t) + \int_0^t R^{(n)}(t-s) \widetilde{N}_H^{(n)} (s)ds \\
& =  \mathcal{I}_{H^{(n)}}(t) + \int_0^t R^{(n)}(t-s) G^*(s)ds + \int_0^t R^{(n)}(t-s) \big(\widetilde{N}_H^{(n)}-G^*\big)(s)ds.
\end{split}
 \eeqlb
 Since $\sup_{n\geq 1}\mathcal{I}_{R^{(n)}}(T)<\infty$ and  $\sup_{t\in[0,T]}|\widetilde{N}_H^{(n)}(t)-G^*(t) | \to 0$ a.s. as $n\to\infty$ for any $T\geq 0$, the last term on the right side of the second equality tends to $0$ as $n\to\infty$. 
 By the stochastic Fubini theorem, 
 \beqnn
  \int_0^t R^{(n)}(t-s) G^*(s)ds =  \int_0^t \mathcal{I}_{R^{(n)}}(t-s) dG^*(s),
 \eeqnn
 which converges almost surely as $n\to\infty$ to 
 \beqnn
 \int_0^t \frac{t-s}{\sigma}  dG^*(s) =  \frac{1}{\sigma} \int_0^tG^*(s)ds.  
 \eeqnn
 Letting $n\to\infty$ in (\ref{eqn.3007}) and then using the preceding results yields
 \beqnn
  \mathcal{I}_{\Lambda^*_H}(t) = \mathcal{I}_{H^*}(t) +\frac{1}{\sigma}  \int_0^t G^*(s)ds ,\quad t\geq 0. 
 \eeqnn
  Integrating both sides of (\ref{CIR})  and then  using Fubini's theorem, we also have 
  \beqnn
   \mathcal{I}_{\Lambda^*_H}(t) = \mathcal{I}_{H^*}(t) + \frac{1}{\sigma}  \int_0^t\int_0^s \sqrt{\Lambda_{H}^*(r)}dB(r) ds ,\quad t\geq 0.
  \eeqnn
 Comparing the preceding two equations and then using the continuity of $G^*$, we have for any $t\geq 0$ that
 $$G^*(t)= \int_0^t \sqrt{\Lambda_{H}^*(s)}dB(s)= \sigma \big( \Lambda^*_H(t)-H^*(t) \big),\quad a.s. $$ 
 \qed

We are now ready to prove Theorem~\ref{MainThm.04}. 

 \medskip
 \textit{Proof of Theorem~\ref{MainThm.04}.}
 The rescaled process $\Big(\Lambda^{(n)},\mathcal{I}^{(n)}_\Lambda, N^{(n)}, \widetilde{N}^{(n)} \Big)$ is equal in distribution to the process 
 $\Big( \Lambda^{(n)}_H,\mathcal{I}_{\Lambda_H^{(n)}}, N^{(n)}_H, \widetilde{N}^{(n)}_H \Big)$ with $H_{\mu_n}=\mu_0+ \mu_0 \cdot \mathcal{I}_R$. Corollary~\ref{ConvergenceR} yields that 
 \[
 	H^{(n)}(t)=\mu_0/n + \mu_0 \cdot  \mathcal{I}_R(nt)/n \to \mu_0\cdot t/\sigma 
\]	
locally uniformly as $n\to\infty$. The desired convergence hence follows from Lemma~\ref{Lemma.ConverGHP}.
 \qed

\subsubsection{Proof of Theorem ~\ref{Thm.FourLapCIR}.}

 The representation (\ref{FourLapCIR}) and  claim (1) follow from Lemma~\ref{LaplaceCIR} with $H^*(t)=\mu_0\cdot t/\sigma$ for $t\geq 0$.
 To prove claim (2) we consider the exponential kernel
 \beqnn
 \phi(t)=\frac{1}{\sigma } \cdot e^{-t/\sigma},\quad t\geq 0.
 \eeqnn
 In this case, $R(t) \equiv  1/\sigma$ on $\mathbb{R}_+$. For the sequence of Hawkes process $\{(N_{n},\Lambda_{n}) : n \in \mathbb N \}$  with immigration rates
 $ \mu_n=n\cdot \phi$ this yields
 \[
 	H_{\mu_n}=  n\big(\phi+ R*\phi \big) \equiv n/\sigma \quad \mbox{and} \quad  H^{(n)}\equiv 1/\sigma
\]
and the  equation (\ref{eqn.3006}) reduces to
  \beqlb\label{eqn.300601}
 \Lambda^{(n)}_H(t)=   \frac{1}{\sigma}+ \frac{1}{\sigma}\widetilde{N}_H^{(n)}(t)=   \frac{1}{\sigma}+ \int_0^t \frac{1}{\sigma}  \widetilde{N}_H^{(n)}(ds),\quad t\geq 0.
 \eeqlb
 
Let  $f\in L^\infty_{\rm loc}(\mathbb{R}_+;\mathbb{C}_-)$ and $g\in C^1 (\mathbb{R}_+;\mathtt{i}\mathbb{R})$. For the particular choice of $H^{(n)}$ it follows form Lemma~\ref{Lemma.FourLapGHP01} that
 	\beqlb\label{eqn.FouLapGHP}
 \mathbf{E}\big[\exp\big\{  \Lambda^{(n)}_H* f(t) +g* d\widetilde{N}^{(n)}_H(t) \big\}\big]
 = \exp \big\{ V^{(n)}(t)   \big\}, \quad t\geq 0.
 \eeqlb
%
%
 Analogously, let $\Lambda^*_H$ be the unique solution of (\ref{CIR}) with $H^* \equiv 1/\sigma$. 
 Then it follows from Lemma~\ref{LaplaceCIR} that
 \beqlb\label{eqn.300602}
 \mathbf{E}\Big[\exp\Big\{  \Lambda^*_H*f(t) + \int_0^t g(t-s)\sqrt{\Lambda^*_H(s)}dB(s)\Big\}\Big]= \exp\big\{    V^*(t) \big\},\quad t\geq 0.
 \eeqlb
 
If $g\in C^1 (\mathbb{R}_+;\mathtt{i}\mathbb{R})$, then  the preceding representations of the functions allow us to establish the uniform convergence of $V^{(n)}$ to $V^*$ using probabilistic methods. In fact, 
 by the stochastic Fubini theorem,  
 \beqnn
 g* d\widetilde{N}^{(n)}_H(t) \ar=\ar \int_0^t g(t-s) \widetilde{N}_H^{(n)}(ds) =  g(0) \cdot \widetilde{N}^{(n)}_H(t) + \int_0^t g'(t-s)\widetilde{N}^{(n)}_H(s)ds 
 \eeqnn
 and
 \beqnn
 \int_0^t g(t-s)  \sqrt{\Lambda^*_H(s)}dB(s) \ar=\ar \sigma \int_0^t g(t-s) d(\Lambda^*_H(s)-1) = g(0) \cdot ( \sigma \Lambda^*_H(t)- \sigma )   +    g'*( \sigma \Lambda^*_H- \sigma )(t).
 \eeqnn
 Applications of Theorem 3.12 in \cite{HorstXu2022} and Lemma~\ref{Lemma.ConverGHP} show that   
 \beqnn
\Big ( \Lambda^{(n)}_H(t),\widetilde{N}^{(n)}_H(t) \Big)\to \Big(\Lambda^*_H(t), \sigma\Lambda^*_H(t)-\sigma \Big)= \Big(\Lambda^*_H(t),\int_0^t \sqrt{\Lambda^*_H(s)}dB(s) \Big),\quad t\geq 0,
 \eeqnn
  weakly in $\mathbf{D}(\mathbb{R}_+;\mathbb{R}_+\times \mathbb{R})$ as $n\to\infty$. 
  By Skorokhod's representation theorem and the continuity of limit process, we may assume that this convergence holds almost surely and locally uniformly, which yields that 
  \beqlb\label{eqn.300603}
  \sup_{t\in[0,T]} \big|\Lambda^{(n)}_H* f(t) -\Lambda^*_H* f(t) \big| + \sup_{t\in[0,T]} \Big| g* d\widetilde{N}^{(n)}_H(t) - \int_0^t g(t-s)  \sqrt{\Lambda^*_H(s)}dB(s)  \Big|  \to 0 \quad a.s. , 
  \eeqlb
  as $n\to\infty$. 
  From this and (\ref{eqn.FouLapGHP})-(\ref{eqn.300602}), we see that $ \|e^{V^{(n)}}-e^{V^*}\|_{L^\infty_T} $ is equal to
  \beqnn
  \sup_{t\in[0,T]}  \Big| \mathbf{E}\Big[\exp\big\{  \Lambda^{(n)}_H* f(t) +g* d\widetilde{N}^{(n)}_H(t) \big\} -\exp\Big\{  \Lambda^*_H*f(t) + \int_0^t g(t-s)\sqrt{\Lambda^*_H(s)}dB(s)\Big\}\Big] \Big| . 
  \eeqnn  
  By using Jensen's inequality and the fact that the function $e^x$ is uniformly Lipschitz continuous and bounded by $1$ on $\mathbb{C}_-$, the preceding quantity can be bounded by 
  \beqnn
   \mathbf{E}\Big[    \sup_{t\in[0,T]}  \Big|   \Lambda^{(n)}_H* f(t) - \Lambda^*_H*f(t)+g* d\widetilde{N}^{(n)}_H(t) -\int_0^t g(t-s)\sqrt{\Lambda^*_H(s)}dB(s) \Big| \wedge 2\Big]  .
  \eeqnn 
  Using the dominated convergence theorem and (\ref{eqn.300603}), yields that 
  \beqnn
   	\big\|\exp\{V^{(n)} \}-\exp\{V^*\} \big\|_{L^\infty_T} \to 0 \quad \mbox{as} \quad n\to\infty.
 \eeqnn
By Lemma~\ref{Lemma.FourLapGHP01} we know that $V^{(n)}, V^* \in  L^\infty_{\rm loc} (\mathbb{R}_+; \mathbb{C}_-)$. Hence, $\|V^{(n)}- V^*\|_{L^\infty_T}\to 0$ as $n\to\infty$.  
Thus, by Lemma~\ref{Lemma.FourLapGHP01}(2),  
  \beqnn
  \|V^* \|_{L^\infty_T} \leq \sup_{n\geq 1}\| V^{(n)}\|_{L^\infty_T}+ \lim_{n\to\infty}\|V^*-V^{(n)}\|_{L^\infty_T}  \leq C\cdot(\|f\|_{L^\infty_T} +\|g\|_{L^\infty_T}^2 ),  
  \eeqnn
   for some constant $C>0$ that is independent of  $(n,f,g)$. 
 Hence (\ref{UpperBound}) holds uniformly for $f\in L^\infty_{\rm loc}(\mathbb{R}_+;\mathbb{C}_-)$ and $g\in C^1(\mathbb{R}_+;\mathtt{i}\mathbb{R})$. 
   \medskip
   
   For the general functions locally bounded function $g\in L^\infty_{\rm loc}(\mathbb{R}_+;\mathtt{i}\mathbb{R})$, one can always find a sequence $\{g_\epsilon\}_{\epsilon>0}\subseteq C^1(\mathbb{R}_+;\mathtt{i}\mathbb{R})$ satisfying that  for any $T\geq 0$,
   \beqnn
   \|g_\epsilon\|_{L^\infty_T}\leq 2 \|g\|_{L^\infty_T}
   \quad \mbox{and}\quad 
   \lim_{\epsilon\to 0+} \|g_\epsilon-g\|_{L^2_T} =0.
   \eeqnn
   Let $V_\epsilon^*  $ be the unique differentiable solution of (\ref{CBI.Riccati}) with $g$ replaced by $g_\epsilon$.  
   The preceding result shows that 
   \beqlb\label{eqn.30014}
  \|V_\epsilon^* \|_{L^\infty_T}   \leq \frac{C}{2}\cdot \big( \|f\|_{L^\infty_T} +\|g_\epsilon\|_{L^\infty_T}^2 \big) \leq C\cdot  \big(\|f\|_{L^\infty_T} +\|g\|_{L^\infty_T}^2 \big),
   \eeqlb
    for some constant $C>0$ independent of $\epsilon$ and $(f,g)$. 
    By (\ref{CBI.Riccati}),
    \beqnn
    V^*(t)-V_\epsilon^*(t)
    \ar=\ar  \frac{1}{2\sigma}\int_0^t \big[(  V^*(s)+g(s))^2-(V_\epsilon^*(t)+g_\epsilon(s))^2 \big]ds \cr
    \ar=\ar  \frac{1}{2\sigma}\int_0^t  (  V^*(s)-V_\epsilon^*(s)+g(s)-g_\epsilon(s))    (  V^*(s)+V_\epsilon^*(s)+g(s)+g_\epsilon(s)) ds,\quad t\geq 0.
    \eeqnn
 By (\ref{eqn.30014}) and the fact that $\|V^*\|_{L^\infty_T}<\infty$, there exists a constant $C_1>0$ independent of $\epsilon$ such that 
    \beqnn
    \|V^*-V_\epsilon^*\|_{L^\infty_t} \leq C_1 \|g-g_\epsilon\|_{L^1_T} + C_1\int_0^t  | V^*(s)-V_\epsilon^*(s)|  ds , 
    \eeqnn
    for any $t\in[0,T]$.  
    By Gr\"onwall's inequality,  
    \[
    	\|V^*-V_\epsilon^*\|_{L^\infty_T}  \leq C_1 \|g-g_\epsilon\|_{L^1_T}\cdot e^{C_1\cdot T} \to 0 \quad \mbox{as} \quad \epsilon \to 0+.
\] 
  Along with (\ref{eqn.30014}) this yields a constant $C>0$ that is independent of $(f,g)$ such that,
    \beqnn 
    \|V^*\|_{L^\infty_T}\leq  \lim_{\epsilon \to 0+}\|V^*-V_\epsilon^*\|_{L^\infty_T}+  \sup_{\epsilon>0}\|V_\epsilon^*\|_{L^\infty_T} \leq C\cdot(\|f\|_{L^\infty_T} +\|g\|_{L^\infty_T}^2 ).
    \eeqnn
  \qed

 \subsection{Proof of Theorem~\ref{MainThm.07}}\label{sec:MainThm.07}
 
 In this section we establish the rate of convergence of the rescaled processes $X^{(n)}$ to the limiting process $X^*$. As a preparation, the next proposition provides uniform moment estimates for the sequence $$\left\{(\mathcal{I}_{\Lambda^{(n)}},N^{(n)}, \widetilde{N}^{(n)})\right\}_{n\geq 1}$$ that will be used in the proofs that follow.
 
 \begin{proposition}\label{MonementEsti}
 	For each $T\geq 0$ and $p\geq 1$, there exists a constant $C>0$ such that for any $n\geq 1$,
 	\beqlb\label{eqn.3.56}
 	\mathbf{E} \big[ \big|\mathcal{I}_{\Lambda^{(n)}}(T) \big|^{p} \big] 
 	+\mathbf{E}\big[\big|N^{(n)}(T) \big|^{p} \big] 
 	+ \mathbf{E} \Big[\sup_{t\in[0,T]}|\widetilde{N}^{(n)}(t)|^{p}\Big] \leq C .
 	\eeqlb
 	
 \end{proposition}
 \proof From Lemma~\ref{MartRep} and Proposition~\ref{Thm.AsymR}(2), we first have  
 \beqlb\label{eqn.30021}
 \sup_{n\geq 1}\mathbf{E}[N^{(n)}(T)]=\sup_{n\geq 1}\mathbf{E}[|\mathcal{I}_{\Lambda^{(n)}}(T)|] = \mu_0 \cdot \sup_{n\geq 1}\big(T/n+ \mathcal{I}_{R}^2(nT)/n^2\big) <\infty .
 \eeqlb
 Since the martingale $\widetilde{N}$ has quadratic variation $[\widetilde{N}]_t= N(t)$, by using the Burkholder-Davis-Gundy inequality; see Theorem 26.12 in \cite[p.524]{Kallenberg2002},  we also have
 \beqnn
 \sup_{n\geq 1}\mathbf{E} \Big[\sup_{t\in[0,T]}|\widetilde{N}^{(n)}(t)|^2 \Big] 
 \ar=\ar \sup_{n\geq 1}\mathbf{E} \Big[\sup_{t\in[0,T]}|\widetilde{N}(nt)/n|^2 \Big] \cr
 \ar\leq\ar C \cdot \sup_{n\geq 1}\mathbf{E} \Big[N(nT)/n^2 \Big] \cr
 \ar=\ar C \cdot \sup_{n\geq 1}\mathbf{E} \big[ N^{(n)}(T)\big]<\infty.
 \eeqnn 
 By H\"older's inequality, $ \sup_{n\geq 1}\mathbf{E} \big[\sup_{t\in[0,T]}|\widetilde{N}^{(n)}(t)| \big]<\infty$ and hence (\ref{eqn.3.56}) holds with $p=1$. 
 
 It is thus enough to prove that  (\ref{eqn.3.56}) holds for any $p\geq 2$ under the assumption that it holds for $p-1$. 
  Using the Burkholder-Davis-Gundy inequality again, and using the induction hypothesis,  
 \beqlb\label{eqn.3.57}
 \sup_{n\geq 1}\mathbf{E} \Big[\sup_{t\in[0,T]}|\widetilde{N}^{(n)}(t)|^{p}\Big]
 \leq C \cdot \sup_{n\geq 1}\mathbf{E} \big[ |N^{(n)}(T)|^{p/2}\big] <\infty.
 \eeqlb
 Using the change of variables to (\ref{IntegralLambda}) with $\mu \equiv \mu_0$ yields that 
 \beqlb\label{eqn.3.55}
 \mathcal{I}_{\Lambda^{(n)}}(T)= \mu_0 \cdot \big(T/n+ \mathcal{I}_{R}^2(nT)/n^2\big) + \int_0^{T} R^{(n)}(T-s) \widetilde{N}^{(n)}(s)ds.
 \eeqlb
 Raising both sides to the $p$-th power and then using the power inequality,  there exists a constant $C>0$ such that for any  $n\geq 1$,
 \beqnn
 |\mathcal{I}_{\Lambda^{(n)}}(T)|^{p}\leq C \cdot \big|T+ \mathcal{I}_{R}^2(nT)/n\big|^{p} + C \cdot \Big|\int_0^{T} R^{(n)}(T-s) \widetilde{N}^{(n)}(s)ds\Big|^{p}.
 \eeqnn
 The first term on the right side of this inequality is uniformly bounded; see (\ref{eqn.30021}). 
 The second term  can be bounded by $C\cdot \big|\mathcal{I}_{ R^{(n)}}(T)\big|^{p}\cdot \sup_{t\in[0,T]}  |\widetilde{N}^{(n)}(t)|^{p}$. 
 Taking expectations on the both sides of the preceding inequality and then using (\ref{eqn.3.57}), we have that
 \beqnn
 \sup_{n\geq 1} \mathbf{E} \big[ |\mathcal{I}_{\Lambda^{(n)}}(T)|^{p} \big]\leq C + C \cdot \sup_{n\geq 1} \mathbf{E}\Big[ \sup_{t\in[0,T]}  |\widetilde{N}^{(n)}(t)|^{p} \Big]  <\infty. 
 \eeqnn 
 By the power inequality, there exists a constant $C>0$ such that for any $n\geq 1$,
 \beqnn
 \mathbf{E}\big[|N^{(n)}(T)|^{p}\big]\leq C\cdot \mathbf{E}\big[ |\mathcal{I}_{\Lambda^{(n)}}(T)|^{p}\big] +C\cdot n^{-1}\cdot \mathbf{E}\big[|\widetilde{N}^{(n)}(T)|^{p}\big].
 \eeqnn
 The preceding results yield that $\sup_{n\geq 1} \mathbf{E}[|N^{(n)}(T)|^{p}]<\infty$. 
 \qed
 
 As another preparation, we are now going to prove the convergence of time-scaled resolvent $R^{(n)}$ to the limit $R^*$ in $L^2_{\rm loc}(\mathbb{R}_+;\mathbb{R})$ under the assumption  that $\phi$ has bounded variation. 
 For $f \in L^1(\mathbb{R}^d;\mathbb{R})$ with $d\in\mathbb{Z}_+$, let $\mathcal{F}f$ be the Fourier transform of $f$ and $\overline{\mathcal{F}f}$ the corresponding conjugate, i.e., 
 \beqnn
 \mathcal{F}f(z):= \int_{\mathbb{R}^d} e^{\mathtt{i}\cdot 2\pi \langle z, x\rangle}f(x)dx
 \quad\mbox{and}
 \quad 
 \overline{\mathcal{F}f}(z):=\mathcal{F}f(-z),\quad z\in \mathbb{R}^d.
 \eeqnn 
 
 	
 	\begin{proposition}\label{UpperphiR}
 		If $\phi$ has bounded  variation,  for any $\beta>0$ and $T\geq 0$  we have the following.
 		\begin{enumerate}
 			\item[(1)] There exist constants $ C >0$ and $n_0\geq 1$ such that for any $z \in\mathbb{R}$ and $n\geq n_0$,
 			\beqlb \label{eqn.30022}
 			\big|\mathcal{F}R^{(n)}_\beta(z) \big|\leq C\cdot \Big( \frac{1}{|z|} \wedge 1\Big). 
 			\eeqlb
 			
 			\item[(2)] $\displaystyle\big\|R^{(n)}_\beta-R^*_\beta \big\|_{L^2}\to 0$ 
 			and 
 			$\displaystyle	\big\|R^{(n)} -R^* \big\|_{L^2_T}\to 0$ as $n\to\infty$.
 		 
 			\end{enumerate}

 	\end{proposition}
 	\proof 
  We first prove claim (1). 
  By a change of variables, 
 	\beqnn
 	\mathcal{F} \phi^{(n)}_\beta(z) = \frac{1}{n}\int_0^\infty e^{\mathtt{i}\cdot 2\pi\frac{z}{n}\cdot t} e^{-\frac{\beta}{n}\cdot t} \phi(t)dt ,\quad z \in\mathbb{R}.
 	\eeqnn
 	Hence, taking Fourier transforms on both sides of (\ref{eqn.30050}) yields that 
 	\beqlb\label{eqn.50021}
 	\mathcal{F}R^{(n)}_\beta(z) = \frac{1}{n} \int_0^\infty e^{\mathtt{i}\cdot 2\pi \frac{z}{n}\cdot t} e^{-\frac{\beta}{n}\cdot t} R(t)dt = \frac{1}{n} \frac{\int_0^\infty e^{\mathtt{i}\cdot 2\pi\frac{z}{n}\cdot t} e^{-\frac{\beta}{n}\cdot t} \phi(t)dt }{1-\int_0^\infty e^{\mathtt{i}\cdot 2\pi\frac{z}{n}\cdot t} e^{-\frac{\beta}{n}\cdot t} \phi(t)dt } = \frac{\mathcal{F} \phi^{(n)}_\beta(z)}{1-n \mathcal{F} \phi^{(n)}_\beta(z)}.
 	\eeqlb
 	Since the total variation of $\phi^{(n)}_\beta$ is uniformly bounded, by Proposition~5.1 in \cite{Xu2021}, there exist some constants $n_0\geq 1$ and $C_1,C_2>0$ such that for any $n\geq n_0$ and $z\in\mathbb{R}$,
 	\beqnn
 	\big|\mathcal{F} \phi^{(n)}_\beta(z)  \big| \leq \frac{C_1}{n} \Big( \frac{1}{|z|/n} \wedge 1 \Big) = C_1 \cdot \Big(\frac{1}{|z|} \wedge \frac{1}{n} \Big)
 	\quad \mbox{and}\quad  	\big|1-n\mathcal{F} \phi^{(n)}_\beta(z) \big| \geq C_2\cdot  \Big(  \frac{|z|}{n}\wedge 1\Big). 
 	\eeqnn
 	As a result, 
 	\[
 	\big|\mathcal{F}R^{(n)}_\beta(z) \big| \leq C/|z| \quad \mbox{ uniformly in } \quad n\geq n_0, ~~ z\in\mathbb{R}.
 	\] 
 	Since we laso have that $ \sup_{z\in\mathbb{R}} \big|\mathcal{F}R^{(n)}_\beta(z) \big| \leq \|R^{(n)}_\beta\|_{L^1}$, 
 	claim (1) follows by putting the two estimates together and then using  (\ref{eqn.3.23}).
 	
 	We now start to prove claim (2). 
 	 For each $z\in\mathbb{R}$, it follows from Corollary~\ref{ConvergenceR} that $ \mathcal{F}R^{(n)}_\beta(z)- \mathcal{F}R^*_\beta(z)\to0$  as $n\to\infty$. 
 	Applications of (\ref{eqn.30022}) and the Plancherel theorem give that 
 	\[
 	\mathcal{F}R^{(n)}_\beta \in L^2(\mathbb{R};\mathbb{C}) \quad \mbox{and} \quad R^{(n)}_\beta \in L^2(\mathbb{R}_+;\mathbb{R}_+).
 	\] 
 	Using  (\ref{eqn.30022}) again and then the dominated convergence theorem, we see that
 	\beqnn
 	\lim_{n\to\infty} \big\|R^{(n)}_\beta-R^*_\beta \big\|_{L^2}= \lim_{n\to\infty} \big\|\mathcal{F}R^{(n)}_\beta - \mathcal{F}R^*_\beta \big\|_{L^2}=
 	\Big\| \lim_{n\to\infty}\big( \mathcal{F}R^{(n)}_\beta - \mathcal{F}R^*_\beta \big)\Big\|_{L^2}=0.
 	\eeqnn
 	For the second convergence, we notice that $|R^{(n)}(t) -R^*(t)|= e^{\beta t} |R^{(n)}_\beta(t) -R^*_\beta(t)|$ and hence
 	\beqnn
 	\big\|R^{(n)} -R^* \big\|_{L^2_T} \leq  e^{\beta T}   \big\|R^{(n)}_\beta -R^*_\beta  \big\|_{L^2_T} \leq  e^{\beta T}   \big\|R^{(n)}_\beta -R^*_\beta  \big\|_{L^2},
 	\eeqnn
 	which also vanishes as $n\to\infty$. 
 	\qed

 In view of (\ref{FourLapCIR}) and (\ref{eqn.FouLapGHP}), the difference between the Laplace transforms of $X^{(n)}$ and $X^*$ can be controlled by a functional of the difference of $V^{(n)}$ and $V^*$. Extending  the proof of Proposition~\ref{Prop.ConV}, the next lemma provides an upper estimate of $\| V^{(n)} - V^*  \|_{L^\infty_T}$ and $	\| V^{(n)} - V^*  \|_{L^1_T}$ for general $\mathbb{C}_-$-valued density functions $f$.

  \begin{lemma}\label{Lemma.ConV01}
  	Let $V^{(n)}$ and $V^*$  be  solutions of  (\ref{Laplace03}) and (\ref{CBI.Riccati}), respectively with $w(dt)=f(t)dt$. 
 	For any $T\geq0$, the following hold.
 	\begin{enumerate}
 		\item[(1)] 	There exits a constant $C>0$ such that for any $n\geq 1$,  $f\in C^1(\mathbb{R}_+;\mathbb{C}_-)$ and $g\in C^1(\mathbb{R}_+;\mathtt{i}\mathbb{R})$, 
 		\beqlb\label{eqn.3.25.01}
 		\| V^{(n)} - V^*  \|_{L^\infty_T}
 		\ar\leq\ar C\cdot \big(1+ \|f\|_{L^\infty_T}  + \|f'\|_{L^\infty_T}   +\|g\|_{L^\infty_T}  + \|g'\|_{L^\infty_T}  \big)^6 \cr
 		\ar\ar\cr
 		\ar\ar \times \exp\big\{C\cdot (\|f\|_{L^\infty_T}+\|g\|_{L^\infty_T}+\|g\|_{L^\infty_T}^2)\big\}\cdot \big( n^{-1} +\|\mathcal{I}_{R^{(n)}}-\mathcal{I}_{R^*}\|_{L^\infty_T} \big).
 		\eeqlb

 		\item[(2)]  there exists a constant $C>0$ such that for any $n\geq 1$, $f\in L^\infty_{\rm loc}(\mathbb{R}_+;\mathbb{C}_-)$ and $g\in L^\infty_{\rm loc}(\mathbb{R}_+;\mathtt{i}\mathbb{R})$,  
 		\beqlb\label{eqn.3.25.02}
 		\| V^{(n)} - V^*  \|_{L^1_T}
 		\ar\leq\ar C\cdot \big(1+ \|f\|_{L^\infty_T}   +\|g\|_{L^\infty_T}   \big)^6\cr
 		\ar\ar\cr
 		\ar\ar \times  \exp\big\{C\cdot (\|f\|_{L^\infty_T}+\|g\|_{L^\infty_T}+\|g\|_{L^\infty_T}^2)\big\}  \cdot \big( n^{-1}+ \|R^{(n)} -R^* \|_{L^2_T}\big).
 		\eeqlb
 		\end{enumerate}

 \end{lemma}
 \proof Let $V^*_g:=V^*+g$ and $V^{(n)}_g:=V^{(n)}+g$. 
 We have $V^{(n)} - V^*=V^{(n)}_g -V^*_g $ and by (\ref{UpperBound}) and the first inequality in (\ref{eqn.3011}),
 \beqlb\label{eqn.5002}
 \sup_{n\geq 1} \|V^{(n)}_g\|_{L^\infty_T} + \|V^*_g \|_{L^\infty_T} \leq  C\cdot \big(\|f\|_{L^\infty_T} + \|g\|_{L^\infty_T} + \|g\|_{L^\infty_T}^2\big).
 \eeqlb
 
 \begin{enumerate}
 	\item[(1)] Notice that (\ref{eqn.50001}) also holds with $(V^{(n)},V^*)$ replaced by $(V^{(n)}_g,V^*_g)$.  Similarly as in the proof of Proposition~\ref{Prop.ConV}, there exists a constant $C>0$ such that for any $n\geq 1$,  $f\in C^1(\mathbb{R}_+;\mathbb{C}_-)$ and $g\in C^1(\mathbb{R}_+;\mathtt{i}\mathbb{R})$, 
 	\beqlb\label{eqn.5003}
 	\| V^{(n)} - V^*  \|_{L^\infty_T} 
 	\ar=\ar \| V^{(n)}_g - V^*_g  \|_{L^\infty_T} \cr
 	\ar\leq\ar C\cdot \exp\big\{C\cdot \big(\|V^{(n)}_g\|_{L^\infty_T} + \|V^*_g\|_{L^\infty_T} \big) \big\} \cdot \Big( \|V^{(n)}_g\|_{L_T^\infty}^3 \cdot \mathcal{I}_{R^{(n)}}(T)\cdot n^{-1} \cr
 	\ar\ar + \big(   \|f\|_{L^\infty_T} + \|f'\|_{L^\infty_T} + \|V^*_g\|_{L^\infty_T} \cdot \|V^{*'}_g\|_{L^\infty_T} \big)\cdot \big \|\mathcal{I}_{R^{(n)}} -\mathcal{I}_{R^*} \big\|_{L^\infty_T}  \Big). 
 	\eeqlb
  with the constant $\beta$ is defined as in (\ref{eqn.beta}). 
  By (\ref{CBI.Riccati}),  we have $V^{*'}_g= f/\sigma +g' +  |V^{*}_g|^2/(2\sigma)$. 
  Using  the Minkowski inequality and (\ref{eqn.5002}) yield that 
 	\beqnn
 	 \|V^{*'}_g \|_{L^\infty_T}  \leq C \big( \|f\|_{L^\infty_T}+ \|f\|_{L^\infty_T}^2 + \|g\|_{L^\infty_T}^2 + \|g\|_{L^\infty_T}^4 + \|g'\|_{L^\infty_T} \big).
 	\eeqnn
    Plugging this and  (\ref{eqn.5002}) back into (\ref{eqn.5003}) implies the desired upper bound  (\ref{eqn.3.25.01}) immediately.

  \item[(2)] For  a constant  $\beta>0$ to be specified later, let
  \beqnn
  f_\beta(t):=e^{-\beta t}f(t),\quad
  V^{(n)}_{g,\beta}(t):=e^{-\beta t}V^{(n)}_g(t)
  \quad \mbox{and}\quad 
  V^*_{g,\beta}(t):=e^{-\beta t} V^*_g(t),\quad t\geq 0. 
  \eeqnn
 It is easy to identify that $V^{(n)}_{g,\beta}-V^*_{g,\beta} = \sum_{i=1}^3  \varepsilon_{i,\beta}^{(n)} $ with 
  \beqlb 
 \varepsilon_{1,\beta}^{(n)}(t) \ar :=\ar \int_0^t\Big( f_\beta+  \frac{1 }{2}V^*_gV^*_{g,\beta}\Big)(t-s)\cdot \big(R^{(n)}_\beta-R^*_\beta\big)(s)ds, \\
 \varepsilon_{2,\beta}^{(n)}(t) \ar :=\ar \frac{1}{2}\int_0^t \big((V^{(n)}_g + V^*_g ) (V^{(n)}_{g,\beta}- V^*_{g,\beta})\big)(t-s)R^{(n)}_\beta(s)ds, \\
 \varepsilon_{3,\beta}^{(n)}(t) \ar  :=\ar  n^2 \int_0^t \Big(e^{V^{(n)}_g/n}-1- \frac{V^{(n)}_g}{n}-\frac{1}{2}\Big(\frac{V^{(n)}_g}{n}\Big)^2\Big)(t-s) e^{-\beta(t-s)} R^{(n)}_\beta (s)ds.
 \eeqlb
 The Minkowski inequality induces that $\|V^{(n)}_{g,\beta} -V^*_{g,\beta}\|_{L^1_T} \leq \sum_{i=1}^3 \|\varepsilon_{i,\beta}^{(n)}\|_{L^1_T}$. 
 Applying Young's convolution inequality
 to $\|\varepsilon_{1,\beta}^{(n)}\|_{L^1_T}$ and then using H\"older's inequality as well as (\ref{eqn.5002}), shows that there exists a constant $C>0$ that  depends only on $T$ such that 
 \beqlb\label{eqn.3.58}
 \big\|\varepsilon_{1,\beta}^{(n)}\big\|_{L^1_T}
 \ar\leq\ar \big\|f_\beta+   V^*_gV^*_{g,\beta}/2\big\|_{L^1_T}\cdot \big\|R^{(n)}_\beta-R^*_\beta \big\|_{L^1_T} \cr
 \ar\leq\ar T\cdot \big(  \big\|f \big\|_{L^\infty_T}  +  \big\|  V^*_g \big\|_{L^\infty_T}^2 \big)\cdot \big\|R^{(n)}_\beta-R^*_\beta \big\|_{L^1_T} \cr
 \ar\leq\ar \frac{T}{|2\beta|^{1/2}}\cdot \big(  \big\|f \big\|_{L^\infty_T}  +  \big\|  V^*_g \big\|_{L^\infty_T}^2 \big)\cdot \big\|R^{(n)} -R^* \big\|_{L^2_T} \cr
 \ar\leq\ar C\cdot \big( \|f\|_{L^\infty_T} + \|f\|_{L^\infty_T}^2+ \|g\|_{L^\infty_T}^2+\|g\|_{L^\infty_T}^4 \big)\cdot \|R^{(n)} -R^* \|_{L^2_T}. 
 \eeqlb
   Applying the inequality $|e^{x}-1-x-x^2/2|\leq |x|^3$ for any $x\in\mathbb{C}_-$ to $\varepsilon_{3,\beta}^{(n)}$, we have 
 \beqnn
 |\varepsilon_{3,\beta}^{(n)}| \leq n^{-1}\cdot |V_g^{(n)}|^3 * R^{(n)}. 
 \eeqnn
  By using Young's convolution inequality
  to $\|\varepsilon_{3,\beta}^{(n)}\|_{L^1_T}$ and then (\ref{eqn.5002}) as well as Corollary~\ref{ConvergenceR}, 
 \beqlb\label{eqn.3.59}
   \|\varepsilon_{3,\beta}^{(n)}\|_{L^1_T} 
   \ar\leq\ar \|V_g^{(n)}\|_{L^3_T}^3 \cdot  \mathcal{I}_{R^{(n)} }(T)   \cdot n^{-1} \cr
   \ar\leq\ar   \|V_g^{(n)}\|_{L^\infty_T}^3 \cdot T \cdot  \mathcal{I}_{R^{(n)} }(T)   \cdot n^{-1} 
   \leq  C\cdot (\|f\|_{L^\infty_T} + \|g\|_{L^\infty_T} +\|g\|_{L^\infty_T}^2)^3\cdot n^{-1},
 \eeqlb
  for some constant $C$ independent of $n$, $f$ and $g$. 
   By applying Young's convolution inequality to $\|\varepsilon_{2,\beta}^{(n)}\|_{L^1_T}$,  we also get that 
 \beqnn
 \|\varepsilon_{2,\beta}^{(n)}\|_{L^1_T} \ar\leq \ar  \|(V^{(n)}_g + V^*_g ) (V^{(n)}_{g,\beta}- V^*_{g,\beta})\|_{L^1_T}\cdot \|R^{(n)}_\beta\|_{L^1_T}.
 \eeqnn
 Applications of Lemma~\ref{Lemma.FourLapGHP01}(2) and Theorem~\ref{Thm.FourLapCIR}(2) imply that 
 \beqnn
 \big\|(V^{(n)}_g + V^*_g ) (V^{(n)}_{g,\beta}- V^*_{g,\beta})\big\|_{L^1_T}
 \ar\leq\ar \big(\|V^{(n)}_g\|_{L^\infty_T} + \|V^*_g \|_{L^\infty_T}  \big) \cdot T\cdot  \big\|V^{(n)}_{g,\beta}- V^*_{g,\beta}\big\|_{L^1_T} \cr
 \ar\leq\ar C_0 \cdot\big( \|f\|_{L^\infty_T}+\|g\|_{L^\infty_T} +\|g\|_{L^\infty_T}^2\big)\cdot \big\|V^{(n)}_{g,\beta}- V^*_{g,\beta}\big\|_{L^1_T} , 
 \eeqnn
 for some constant $C_0>0$ depending only on $T$. Hence 
 \beqnn
 \|\varepsilon_{2,\beta}^{(n)}\|_{L^1_T} \leq C_0 \cdot\big( \|f\|_{L^\infty_T}+\|g\|_{L^\infty_T}\big)\cdot \|R^{(n)}_\beta\|_{L^1_T} \cdot \|V^{(n)}_{g,\beta}- V^*_{g,\beta}\|_{L^1_T}.
 \eeqnn
 Choosing $\beta= 4C_0 (\|f\|_{L^\infty_T}+\|g\|_{L^\infty_T}+\|g\|_{L^\infty_T}^2)/\sigma$ and then using the second the inequality in (\ref{eqn.3.23}), yields that 
 \beqnn
 \|\varepsilon_{2,\beta}^{(n)}\|_{L^1_T}\leq \frac{1}{2} \|V^{(n)}_{g,\beta}- V^*_{g,\beta}\|_{L^1_T}.
 \eeqnn
 Taking this, (\ref{eqn.3.58}) and (\ref{eqn.3.59}) back into $\|V^{(n)}_{g,\beta} -V^*_{g,\beta}\|_{L^1_T}$, we have
 \beqnn
 \big\|V^{(n)}_{g,\beta} -V^*_{g,\beta}\big\|_{L^1_T} \leq C\cdot (1+\|f\|_{L^\infty_T}+\|g\|_{L^\infty_T})^6 \cdot \big( n^{-1} + \|R^{(n)} -R^* \|_{L^2_T} \big)
 \eeqnn
 and  the desired result follows from the fact that 
 \beqnn\|V^{(n)} -V^*\|_{L^1_T}=\|V^{(n)}_{g} -V^*_{g}\|_{L^1_T}\leq e^{\beta T}\|V^{(n)}_{g,\beta} -V^*_{g,\beta}\|_{L^1_T}.
 \eeqnn
   	\qed
  	\end{enumerate}

 \medskip
 \textit{Proof of Theorem~\ref{MainThm.07}.}
 In what follows we prove the two inequalities (\ref{UpperBLapTran}) and (\ref{UpperBWass}). The corresponding result for kernels of bounded variation  can be proved in the same way. 

{\bf Part I.}
 We first prove the inequality (\ref{UpperBLapTran}). To this end, we recall that 
 \beqnn
 	\mathcal{I}_\Lambda^{(n)}= \mathcal{I}_{\Lambda^{(n)}} 
 	\quad \mbox{and} \quad 
 	N^{(n)}(dt)= \Lambda^{(n)}(t)dt+ n^{-1}\cdot \widetilde{N}^{(n)}(dt).
\eeqnn
Hence, for each $f = (f_1, f_2, f_3, f_4) \in \mathcal{A}_K^-$, it holds that $ f_2*\mathcal{I}_\Lambda^{(n)} = \mathcal{I}_{f_2}* \Lambda^{(n)} $
 and  $ f_3 *dN^{(n)} = f_3 *  \Lambda^{(n)} + (f_3/n)*d\widetilde{N}^{(n)}$, which allow us to express the convolutions $f*X^{(n)}(t)$ and $ f*X^*(t)$ as 
 \beqlb
 f*X^{(n)}(t)
 \ar=\ar  (f_1+\mathcal{I}_{f_2}+f_3)*\Lambda^{(n)}(t) +  \int_0^t(f_3/n+f_4)(t-s)\widetilde{N}^{(n)}(dt), \label{eqn.3.49}\\
 f*X^*(t)
 \ar=\ar (f_1+\mathcal{I}_{f_2}+f_3)*\Lambda^*(t) +  \int_0^t  f_4 (t-s)\sqrt{\Lambda^*(s)}dB(s). \label{eqn.3.50}
 \eeqlb
 Let $W^{(n)}= n^2\big( e^{(V^{(n)}+g)/n}-1- (V^{(n)}+g)/n\big)$ and $V^{(n)}$ be the unique solution of (\ref{Laplace03}) with $w(dt)= (f_1+\mathcal{I}_{f_2} +f_3)(t)dt$ and $g = f_3/n+f_4 $.
 By the inequality $|e^x-1-x|\leq |x|^2/2$ for any $x\in \mathbb{C}_-$, the power inequality and Lemma~\ref{Lemma.FourLapGHP01}(2), there exists a constant $C>0$ such that for any  $f\in\mathcal{A}_K^-$,
 \beqlb\label{eqn.30017}
 \sup_{n\geq1}\| V^{(n)} \|_{L^\infty_T}  \leq  C\cdot (K+K^2)
 \quad \mbox{and}\quad
 \sup_{n\geq1} \| W^{(n)}  \|_{L^\infty_T} \leq  \sup_{n\geq1} \| V^{(n)}\|_{L^\infty_T}^2+\|g \|_{L^\infty_T}^2 \leq C\cdot (K^2+K^4).
 \eeqlb
 Applying (\ref{eqn.FouLapGHP01}) with $H^{(n)}= \mu_0/n+ \mu_0\cdot \mathcal{I}_{R^{(n)}}$, we have
 \beqnn
 \mathbf{E}\big[e^{f*X^{(n)}(t)}]= \exp \Big\{ \frac{\mu_0}{n} \big(\mathcal{I}_{f_1+\mathcal{I}_{f_2}+f_3 }(t)+ \mathcal{I}_{W^{(n)}}(t)\big)+ \mu_0\cdot \mathcal{I}_{V^{(n)}}(t) \Big\}.
 \eeqnn
 Similarly, let $V^*$ be the unique solution of (\ref{CBI.Riccati}) with $w(dt)= (f_1 +\mathcal{I}_{f_2 } +f_3 )(t)dt$ and $g= f_4$.
 By (\ref{FourLapCIR}),  
 $$\mathbf{E}\big[e^{f*X^*(t)}\big] =  \exp\big\{ \mu_0 \cdot \mathcal{I}_{V^*}(t)  \big\}. $$
 
 By Theorem~\ref{Thm.FourLapCIR}(2), the estimate (\ref{eqn.30017}) and the local Lipschitz continuity of $e^x$ on $\mathbb{C}$, there exists a constant $C>0$ such that for any $n\geq 1$ and $t\in[0,T]$, 
 \beqnn
 \big| \mathbf{E}\big[e^{f*X^{(n)}(t)}] -\mathbf{E}[e^{f*X^*(t)}\big] \big|
 \leq \frac{C}{n} \big|\mathcal{I}_{f_1 +\mathcal{I}_{f_2}  +f_3}(t)+ \mathcal{I}_{W^{(n)}}(t)\big| + C \big| \mathcal{I}_{V^{(n)}}(t) - \mathcal{I}_{V^*}(t) \big|.
 \eeqnn
Moreover, since 
$\|\mathcal{I}_{f_1 +\mathcal{I}_{f_2} +f_3 }\|_{L^\infty_T}\leq C\cdot K $ 
and 
$$\sup_{n\geq 1}\|\mathcal{I}_{W^{(n)}}\|_{L^\infty_T} \leq T\cdot \sup_{n\geq 1}\| W^{(n)} \|_{L^\infty_T} \leq C\cdot( K^2 +K^4)$$
for some constant $C>0$ independent of  $f$, we have uniformly in $n\geq 1$ that
 \beqnn
 \sup_{t\in[0,T]}\big| \mathbf{E}\big[e^{f*X^{(n)}(t)}] -\mathbf{E}[e^{f*X^*(t)}\big] \big|
 \leq C\cdot \big(n^{-1}+\|  V^{(n)}- V^* \|_{L^\infty_T}\big).
 \eeqnn
 The desired inequality  (\ref{UpperBLapTran}) now follows  from Lemma~\ref{Lemma.ConV01}(1).  
 \medskip
 
 {\bf Part II.}
 We now prove the inequality (\ref{UpperBWass}). 
 Let $\mathcal{A}_K^\circ$ be the set of functions $f\in \mathcal{A}_K$ with ${\tt Im}f =0$.  
 For each $n\geq 1$ and $f^{1},f^{2} \in \mathcal{A}_K^\circ$, let 
 \beqlb
 X^{(n)}_f:= (f^{1}*X^{(n)},f^{2}*X^{(n)}),\quad
 X^*_f :=  (f^{1}*X^*,f^{2}*X^*) , 
 \quad  
  \theta_n:= \big| \log  \big( n^{-1} +\big\|\mathcal{I}_{R^{(n)}}- \mathcal{I}_{R^*}\big\|_{L^\infty_T} \big)\big|. \label{eqn.3.48}
 \eeqlb
 By the definition of Wasserstein distance of measures on $\mathbb{C}$ as recalled in Appendix~\ref{Appendix.DM}, it suffices to prove that there exists constant $C>0$  such that  for any $n\geq 1$,  $F\in \mathcal{H}_{\rm Lip}(\mathbb{R}^2;\mathbb{R}) $ and $f^{1},f^{2} \in \mathcal{A}_K^\circ$, 
 \beqlb\label{eqn.3.45}
\sup_{t\in[0,T]}  \left|\mathbf{E}\big[F(X^{(n)}_{f}(t) )\big]- \mathbf{E} \right[F(X^*_{f}(t) ) \big] \big|
 \leq C\cdot|\theta_n|^{-\kappa}.
 \eeqlb

 To this end, we need the following auxiliary functions. 
 For any $n\geq 1$ and a constant $\varrho \in (0, \frac{1}{2}-\kappa)$, let $\psi_n: \mathbb{R}^2 \mapsto [0,1]$ be a smooth function that satisfies  
 \beqnn
 \psi_n(z) =  
 \begin{cases}
	1, & \mbox{if } |z|\leq |\theta_n|^{\varrho }; \vspace{5pt}\\
	0, & \mbox{if } |z|\geq |\theta_n|^{\varrho }+4
 \end{cases} 
 \quad \mbox{and}\quad
  \|\psi'_n\|_{L^\infty}\leq \frac{1}{2} .
 \eeqnn 
 
 For each $F\in \mathcal{H}_{\rm Lip}(\mathbb{R}^2;\mathbb{R})$,   the function $F_n:= F\cdot \psi_n$ is  absolutely continuous, whose Radon–Nikodym derivative, denoted by $F'_n$, is uniformly bounded by $1$ and has support in $\{z\in\mathbb{R}^2:|z|\leq |\theta_n|^{\varrho }+4\}$. 
 By the Plancherel theorem and the equality $z\cdot\mathcal{F}F_n(z)= \mathcal{F}F'_n(z)$ for any $z\in \mathbb{R}$,
 \beqnn
  \|\mathcal{F}F_n\|_{L^2}=\|F_n\|_{L^2} \leq 2(|\theta_n|^{\varrho }+4)
 \quad\mbox{and}\quad
 \int_\mathbb{R} |z\cdot\mathcal{F}F_n(z)|^2dz= \|\mathcal{F}F_n'\|_{L^2}^2=\|F'_n\|_{L^2}^2 \leq  4(|\theta_n|^{\varrho }+4)^2. 
 \eeqnn 
 By H\"older's inequality, there exists a constant $C>0$ such that for any $z_0\geq 0$ and $n\geq1$,
 \begin{equation}\label{eqn.3.46}
 \begin{split}
 \Big|\int_{|z|\geq z_0} |\mathcal{F}F_n(z)|dz\Big|^2
 & \leq  \int_{|z|\geq z_0} (1+|z|)^2\cdot |\mathcal{F}F_n(z)|^2dz\cdot \int_{|z|\geq z_0} (1+|z|)^{-2}dz    \\
 & \leq  \frac{C}{ 1+z_0}  \int_{\mathbb{R}} \big( |\mathcal{F}F_n(z)|^2 +|z\cdot\mathcal{F}F_n(z)|^2 \big)dz \\
 & \leq  C\cdot \frac{ (1+ |\theta_n|^{\varrho } )^2}{ 1+ z_0} .
\end{split}
 \end{equation}
 
 To prove (\ref{eqn.3.45}), we first notice that 
%
 \beqlb\label{eqn.30019}
 \begin{split}
&  \sup_{t\in[0,T]} \big|\mathbf{E}[F(X^{(n)}_{f}(t) )]- \mathbf{E}[F(X^*_{f}(t)  )] \big| \\
 \leq & \sup_{t\in[0,T]} \big|\mathbf{E}[F_n(X^{(n)}_{f}(t) )]- \mathbf{E}[F_n(X^*_{f}(t)  )] \big|
 + \sup_{t\in[0,T]} \mathbf{P}( |X^{(n)}_{f}(t)| \geq |\theta_n|^{\varrho } )  +\sup_{t\in[0,T]}  \mathbf{P}(  |X^*_{f}(t)| \geq |\theta_n|^{\varrho } ).
 \end{split}
 \eeqlb
 By Chebyshev's inequality and the power inequality,  there exists a constant  $C>0$ such that for any $n\geq 1$, 
 \beqlb \label{eqn.30018}
  \sup_{t\in[0,T]}  \mathbf{P}( |X^{(n)}_{f}(t)| \geq |\theta_n|^{\varrho } ) 
 \ar\leq\ar  \frac{1}{\theta_n} \sup_{t\in[0,T]}  \mathbf{E}[|X^{(n)}_{f}(t)|^{1/\varrho}]  \cr
 \ar\leq\ar \frac{C}{\theta_n}\sum_{i=1}^2   \sup_{t\in[0,T]}  \mathbf{E}\big[|f^{i}*X^{(n)}(t)|^{1/\varrho}\big] \cr
 \ar\leq\ar  \frac{C}{\theta_n}\sum_{i=1}^2   \sup_{t\in[0,T]} \mathbf{E}\big[|(f_1^{i}+\mathcal{I}_{f_2^{i}}+f_3^{i})*\Lambda^{(n)}(t)|^{1/\varrho}\big] \cr
 \ar\ar + \frac{C}{\theta_n}\sum_{i=1}^2  \sup_{t\in[0,T]}  \mathbf{E}\big[ |(f_3^{i}/n + f_4^{i})* d\widetilde{N}^{(n)}(t)|^{1/\varrho} \big] .
 \eeqlb
 The boundedness of $f^{i}$ yields a constant $C>0$ that is independent of $n\geq 1$ such that
 \beqnn
  \sup_{t\in[0,T]} \mathbf{E}\big[|(f_1^{i}+\mathcal{I}_{f_2^{i}}+f_3^{i})*\Lambda^{(n)}(t)|^{1/\varrho}\big] \leq C \mathbf{E}\big[|\mathcal{I}_{\Lambda^{(n)}}(T)|^{1/\varrho}\big] . 
 \eeqnn
Moreover, by Theorem D.1 in \cite{Xu2021b} and the fact that $\|f^{i}\|_{L^\infty_T}\leq K$, 
 \beqnn
\sup_{t\in[0,T]}  \mathbf{E}\big[ |(f_3^{i}/n + f_4^{i})*d\widetilde{N}^{(n)}(t)|^{\frac{1}{\varrho} } \big]  \leq  C \sup_{t\in[0,T]}\mathbf{E}\Big[ \Big|\int_0^t (f_3^{i}/n + f_4^{i})^2(t-s)N^{(n)}(ds)\Big|^{\frac{1}{2\varrho}} \big] 
   \ar\leq\ar C  \mathbf{E}\big[| N^{(n)}(T)|^{\frac{1}{2\varrho}}  \big],
 \eeqnn
 for some constant $C>0$ independent of $n\geq 1$.  
 Taking these estimates back into (\ref{eqn.30018}) and the using Proposition~\ref{MonementEsti},  we have as $n\to\infty$,
 \beqnn
 \sup_{t\in[0,T]}\mathbf{P}( |X^{(n)}_{f}(t)| \geq |\theta_n|^{\varrho } )\leq C/\theta_n\to 0. 
 \eeqnn 
 Similarly, one also can prove that $\sup_{t\in[0,T]} \mathbf{P}( |X^*_{f}(t)| \geq |\theta_n|^{\varrho } ) \to 0$ as $n\to\infty$. 
To obtain (\ref{eqn.3.45}) it thus remains to prove that
 \beqlb\label{eqn.30031}
 \sup_{t\in[0,T]} \big|\mathbf{E}[F_n(X^{(n)}_{f}(t) )]- \mathbf{E}[F_n(X^*_{f}(t)  )] \big|
  \leq C\cdot|\theta_n|^{-\kappa},
 \eeqlb
 uniformly in $n\geq 1$,  $F\in \mathcal{H}_{\rm Lip}(\mathbb{R}^2;\mathbb{R}) $ and $f^{1},f^{2} \in \mathcal{A}_K^\circ$. This can be achieved by employing Fourier analysis. 
Using the Fourier inversion theorem and then  Fubini's theorem, we see that   $F_n=\mathcal{F}\overline{\mathcal{F}F_n}$ and that
 \beqnn
   \mathbf{E}[F_n(X_f^{(n)}(t))]
   \ar = \ar \mathbf{E}\Big[ \int_{\mathbb{R}^2} e^{\mathtt{i}\langle z, X_f^{(n)}(t)\rangle}\overline{\mathcal{F}F_n}(z)dz\Big]= \int_{\mathbb{R}^2}  \mathbf{E}[e^{\mathtt{i}\langle z, X_f^{(n)}(t)\rangle}]\cdot \overline{\mathcal{F}F_n}(z)dz ,\cr
  \mathbf{E}[F_n(X_f^*(t))]
   \ar=\ar \mathbf{E}\Big[  \int_{\mathbb{R}^2} e^{\mathtt{i}\langle z, X_f^*(t)\rangle}\overline{\mathcal{F}F_n}(z)dz \Big] = \int_{\mathbb{R}^2}   \mathbf{E}[e^{\mathtt{i}\langle z, X_f^*(t)\rangle}]\cdot \overline{\mathcal{F}F_n}(z)dz.
 \eeqnn
 For a constant $\delta\in(0,1)$ to be specified later, it is straightforward to verify that  
 \beqlb\label{eqn.30034}
 \sup_{t\in[0,T]} \big|\mathbf{E}[F_n(X^{(n)}_{f}(t) )]- \mathbf{E}[F_n(X^*_{f}(t)  )] \big| \leq \| \varepsilon^{(n)}_{f,1}\|_{L^\infty_T}+ \|\varepsilon^{(n)}_{f,2}\|_{L^\infty_T},
 \eeqlb
 where 
 \beqnn
 \varepsilon^{(n)}_{f,1}(t)\ar:=\ar \int_{|z|\leq \delta\theta_n} \Big(  \mathbf{E}[e^{\mathtt{i}\langle z, X_f^{(n)}(t)\rangle}]-\mathbf{E}[e^{\mathtt{i}\langle z, X_f^*(t)\rangle}]\Big) \overline{\mathcal{F}F_n}(z)dz,\cr
  \varepsilon^{(n)}_{f,2}(t)\ar:=\ar \int_{|z|> \delta\theta_n} \Big(  \mathbf{E}[e^{\mathtt{i}\langle z, X_f^{(n)}(t)\rangle}]-\mathbf{E}[e^{\mathtt{i}\langle z, X_f^*(t)\rangle}]\Big) \overline{\mathcal{F}F_n}(z)dz.
 \eeqnn
 
 We estimate the above terms separately, starting with the second one.  Applying the estimate (\ref{eqn.3.46}) with $z_0=\theta_n$ to $ \varepsilon^{(n)}_{f,2}(t)$, yields a constant $C>0$ such that for any $f\in \mathcal{A}_K$ and $F\in \mathcal{H}_{\rm Lip}$,
 \beqlb\label{eqn.30032}
  \|\varepsilon^{(n)}_{f,2}\|_{L^\infty_T}  \leq 2 \int_{|z|> \delta\theta_n} \big| \overline{\mathcal{F}F_n}(z) \big| dz  \leq C\cdot \frac{ 1+ |\theta_n|^{\varrho } }{ (1+ \theta_n)^{1/2} }
 \leq C\cdot ( 1+ \theta_n)^{\varrho -1/2}\leq C\cdot \theta_n^{-\kappa}.
 \eeqlb 
 
Estimating the first term is more involved. We start with the following representation of the exponents. 
For any $z:=(z_1,z_2) \in\mathbb{R}^2$, let $f_z:=\mathtt{i}z_1f^{1}+\mathtt{i}z_2f^{2} \in C^1(\mathbb{R}_+;\mathtt{i}\mathbb{R}^4)$.
 It follows from (\ref{eqn.3.48}) that 
 \beqnn
	\mathtt{i}\langle z, X_f^{(n)}\rangle = f_z*X^{(n)} \quad \mbox{and} \quad 
	\mathtt{i} \langle z, X_f^*\rangle = f_z*X^*. 
 \eeqnn
Moreover, it is easy to verify that uniformly in $f^{1},f^{2} \in \mathcal{A}^\circ_K$ and $z\in\mathbb{R}^2$,
 \beqlb\label{UpperBfz}
  \|f_z\|_{L^\infty_T} + \|f'_{z}\|_{L^\infty_T}    \leq CK\cdot |z|. 
 \eeqlb 

Next, let $V^{(n)}_z \in L^\infty_{\rm loc}(\mathbb{R}_+;\mathbb{C}_-)$, respectively, $V^*_z \in L^\infty_{\rm loc}(\mathbb{R}_+;\mathbb{C}_-)$ be the unique solution to the Volterra equation (\ref{Laplace03}), respectively, and (\ref{CBI.Riccati}) associated with the measure 
\[
	w(dt)= (f_{z,1}+f_{z,3}+\mathcal{I}_{f_{z,2}} )(t)dt 
\]	
and the respectively functions 
\[
	 g = f_{z,3}/n+f_{z,4} \quad \mbox{and} \quad g=f_{z,4}.
\]
From (\ref{eqn.3.49}) and (\ref{eqn.FouLapGHP01}), respectively, (\ref{eqn.3.50})  and (\ref{FourLapCIR}), it follows that 
\beqnn
 \mathbf{E}\big[e^{\mathtt{i}\langle z, X_f^{(n)}(t)\rangle} \big]
 \ar=\ar  \exp \Big\{ \frac{\mu_0}{n} \big(\mathcal{I}_{f_{z,1} +\mathcal{I}_{f_{z,2}} +f_{z,3}}(t)+ \mathcal{I}_{W^{(n)}_{z}}(t)\big)+ \mu_0\mathcal{I}_{V^{(n)}_z}(t) \Big\}, \quad t\geq 0. 
 \eeqnn
where $W_{z}^{(n)}=n^2(e^{(V^{(n)}_z+g)/n} -1- (V^{(n)}_z+g)/n) \in L^\infty_{\rm loc}(\mathbb{R}_+;\mathbb{C}_-)$ and that
 \beqnn
 \mathbf{E}\big[ e^{\mathtt{i}\langle z, X_f^*(t)\rangle} \big]
 =  \exp\big\{ \mu_0 \cdot \mathcal{I}_{V^*_z}(t)  \big\},\quad t\geq 0. 
 \eeqnn
 Taking  these two representations back into $\varepsilon^{(n)}_{f,1}(t)$, we see that $ \varepsilon^{(n)}_{f,1}(t)=\varepsilon^{(n)}_{f,1,1}(t)+\varepsilon^{(n)}_{f,1,2}(t)$ with
 \beqnn
 \varepsilon^{(n)}_{f,1,1}(t)\ar :=\ar \int_{|z|\leq \delta\theta_n} \Big( \exp \Big\{ \frac{\mu_0}{n} \big(\mathcal{I}_{f_{z,1}+ f_{z,3} +\mathcal{I}_{f_{z,2}} }(t)+ \mathcal{I}_{W^{(n)}_{z}}(t)\big)\Big\}-1\Big) \cdot \exp\big\{ \mu_0\mathcal{I}_{V^{(n)}_z}(t) \big\} \overline{\mathcal{F}F_n}(z)dz,\cr
 \varepsilon^{(n)}_{f,1,2}(t)\ar :=\ar \int_{|z|\leq \delta\theta_n} \Big( \exp \big\{ \mu_0\mathcal{I}_{V^{(n)}_z}(t) \big\} -\exp\big\{ \mu_0 \cdot \mathcal{I}_{V^*_z}(t)  \big\}\Big) \overline{\mathcal{F}F_n}(z)dz . 
 \eeqnn
 
 It remains to prove that the above terms of both of order $o(1/\theta_n)$ as $n\to\infty$. 
 Since ${\tt Re}f_z =0$ and ${\tt Re}V^{(n)}_{z}, {\tt Re}W^{(n)}_{z}\leq 0$, and since $|e^x-1|\leq |x|$ for any $x\in\mathbb{C}_-$ the first terms admits the pointwise estimate
 \beqnn
 |\varepsilon^{(n)}_{f,1,1}(t)| \leq \frac{\mu_0}{n}\int_{|z|\leq \delta\theta_n} \big| \mathcal{I}_{f_{z,1} +f_{z,3}+\mathcal{I}_{f_{z,2}} }(t)+\mathcal{I}_{W^{(n)}_{z}}(t)\big|  \cdot\big| \overline{\mathcal{F}F_n}(z)\big| dz.
 \eeqnn
 Moreover, by Lemma~\ref{Lemma.FourLapGHP01}(2) and (\ref{UpperBfz}), there exists a constant $C>0$ such that $\| W^{(n)}_{z} \|_{L^\infty_T} \leq   CK \cdot |z|$ for any $n\geq 1$, $z\in\mathbb{R}^2$ and $f^{1},f^{2} \in \mathcal{A}_K$.  
 Since $\| \mathcal{F}F_n \|_{L^\infty}\leq \|F_n\|_{L^1}  \leq 4 (|\theta_n|^{\varrho }+4)^2$, we conclude that as $n\to\infty$, 
 \beqlb \label{eqn.30033}
 \|\varepsilon^{(n)}_{f,1,1}\|_{L^\infty_T}  \leq  \frac{C}{n}\int_{|z|\leq \delta\theta_n} (1+|z|)\cdot\big|  \mathcal{F}F_n (z)\big| dz
 \leq C\cdot \frac{|\theta_n|^{2+2\varrho}+1}{n} = o(1/\theta_n).
 \eeqlb
 
We proceed with the second term.  Since $\mathcal{I}_{V^{(n)}_z}(t), \mathcal{I}_{V^*_z}(t)\in\mathbb{C}_- $ and $|e^{x}-e^{y}|\leq |x-y|$ uniformly in $x,y\in\mathbb{C}_-$,  
 \beqnn
\|\varepsilon^{(n)}_{f,1,2}\|_{L^\infty_T}
 \ar\leq\ar \mu_0\int_{|z|\leq \delta\theta_n} \big\|\mathcal{I}_{V^{(n)}_z} -  \mathcal{I}_{V^*_z} \big\|_{L^\infty_T} \cdot \big| \overline{\mathcal{F}F_n}(z) \big| dz
 \leq \mu_0\cdot T\cdot \int_{|z|\leq \delta\theta_n} \big\| V^{(n)}_z  -  V^*_z  \big\|_{L^\infty_T} \cdot \big| \overline{\mathcal{F}F_n}(z) \big| dz . 
 \eeqnn
 By Lemma~\ref{Lemma.ConV01}(1) and (\ref{UpperBfz}),  there exist a constant $C_0>0$ independent of $\delta$ such that  for any $n\geq 1$, $z\in\mathbb{R}^2$, $f^1,f^2\in\mathcal{A}_K^\circ$ and $F\in \mathcal{H}_{\rm Lip}$,
 \beqnn
 \|V^{(n)}_z-V^*_z\|_{L^\infty_T}\leq C_0e^{C_0K\cdot |z|} e^{-\theta_n}
 \quad \mbox{and}\quad
 \|\varepsilon^{(n)}_{f,1,2}\|_{L^\infty_T} \leq C_0 e^{-\theta_n} \cdot \int_{|z|\leq \delta\theta_n} e^{C_0K\cdot |z|}\cdot \big| \mathcal{F}F_n (z)\big| dz.
 \eeqnn
It is not difficult to show that 
 \beqnn
	\|\varepsilon^{(n)}_{f,1,2}\|_{L^\infty_T} \leq C_0 e^{(C_0K\delta-1)\theta_n} \cdot 4|\delta\theta_n|^2\cdot \| \mathcal{F}F_n\|_{L^\infty}.
 \eeqnn
 Since the support of $F_n$ belongs to $\{z\in\mathbb{R}^2:|z|\leq |\theta_n|^{\varrho }+4\}$, we also have that 
 \[
	 \| \mathcal{F}F_n\|_{L^\infty} \leq  \| F_n\|_{L^1} \leq 4(|\theta_n|^{\varrho }+4)^2.
	\] 
 Combining all estimates and then choosing $\delta\in(0,\frac{1}{C_0K})$, we see that as $n\to\infty$, 
 \beqnn
 \|\varepsilon^{(n)}_{f,1,2}\|_{L^\infty_T} \leq C_0  \cdot 4|\delta\theta_n|^2\cdot 4(|\theta_n|^{\varrho }+4)^2 \cdot e^{(C_0K\delta-1)\theta_n}= o(1/\theta_n).
 \eeqnn
 Taking this together with (\ref{eqn.30032}) and (\ref{eqn.30033}) backing into (\ref{eqn.30034}), the desired inequality (\ref{eqn.30031}) follows. 
 \qed

  \subsection{Proof of Proposition~\ref{Thm.ConvRateRLinfty}} \label{sec:MainThm.08} 

 Before giving the detailed proof of Proposition~\ref{Thm.ConvRateRLinfty}, we first provide asymptotic results and upper estimates for  the Laplace-Stieltjes transforms of the functions ${\it\Phi}$ and ${\it\Psi}_k$ introduced in \eqref{LapPhi} and \eqref{LapPsi}, respectively. They will play a key role in analyzing the convergence rate of rescaled resolvent.  
 Let $ {\it\bar\Psi}_1$ be the tail-function of ${\it\Psi }_1$, i.e., 
 \beqnn
 {\it\bar\Psi}_1(t):=\sigma-{\it\Psi }_1(t),\quad t\geq 0. 
 \eeqnn

 \begin{proposition}\label{TauberPhi01}
 	If $ {\it\Phi }\in {\rm RV}^\infty_{-\alpha-1} $ with $\alpha\in(0,1)$, we have as $\lambda \to 0+$,
 	$$
 	\int_0^\infty (e^{-\lambda t}-1+\lambda t) d{\it\Phi }(t)  \sim  \Gamma(-\alpha) {\it\Phi}(1/\lambda )
 	\quad \mbox{and} \quad 
 	\int_0^\infty (1-e^{-\lambda t})d{\it\Psi }_1(t)\sim - (\alpha+1)\Gamma(-\alpha)\cdot \frac{{\it\Phi}(1/\lambda)}{\lambda}.
 	$$
 \end{proposition}
 \proof The first result follows from Proposition~\ref{Prop.A.6} with $n=1$.
 Using integration by parts and then Proposition~\ref{Thm.Karamata},  we have 
 \beqlb\label{eqn.30037}
 {\it \bar{\Psi}}_1(t)= t{\it\Phi}(t) + \int_t^\infty {\it\Phi}(s)ds \sim (\alpha+1)/\alpha\cdot t {\it\Phi}(t) \in {\rm RV}^\infty_{-\alpha}
 \eeqlb
 and the second asymptotic result follows from Proposition~\ref{Prop.A.6}  with $n=0$ .
 \qed

  \begin{proposition}\label{TauberPhi02}
 	If $ {\it\Phi }\in {\rm RV}^\infty_{-\alpha-1} $ with $\alpha\in(0,1)$, we have as $z \to 0$,
 	\beqlb\label{eqn.30035}
 	\int_0^\infty (e^{\mathtt{i}zt}-1- \mathtt{i}zt) d{\it\Phi }(t) \sim  -\exp\big\{\mathtt{i}\cdot{\rm sgn}(z)\frac{\pi}{2}(\alpha-1)\big\}\cdot \Gamma(-\alpha) \cdot {\it\Phi}(1/|z|).
 	\eeqlb
 \end{proposition}
 \proof We can express the left-hand side of (\ref{eqn.30035}) as $  G_1(z)+ \mathtt{i} G_2(z)$ with
 \beqnn
 G_1(z):=  \int_0^\infty \big( \cos(zt) -1\big)d{\it\Phi }(t)
 \quad\mbox{and}\quad
 G_2(z):=  \int_0^\infty \big( \sin(zt)-zt \big) d{\it\Phi(t)}.
 \eeqnn
 Using integration by parts twice to $G_1$ and $G_2$ respectively, we have
 \beqnn
 G_1(z) =  z^2 \int_0^\infty \cos(zt)dt \int_t^\infty {\it\Phi(s)}ds
 \quad \mbox{and}\quad
 G_2(z) = -z^2 \int_0^\infty \sin(zt) dt \int_t^\infty {\it\Phi(s)}ds.
 \eeqnn
 It follows from Proposition~\ref{Thm.Karamata} that  $\int_t^\infty {\it\Phi}(s) ds \sim  t{\it\Phi}(t) /\alpha \in {\rm RV}^\infty_{-\alpha}$.
 Since the function $\cos(\cdot)$ is even, it follows from Proposition~\ref{Thm.AbelThm} that as $z\to 0$, 
 \beqnn
 G_1(z)=G_1(|z|) \sim \Gamma(1-\alpha) \sin(\pi \alpha/2) \cdot |z|  \int_{1/|z|}^\infty {\it\Phi(s)}ds \sim  - \Gamma(-\alpha) \sin(\pi \alpha/2) \cdot  {\it\Phi}(1/|z|) . 
 \eeqnn 
 Similarly, since the function $\sin(\cdot)$ is odd, we also have as $z\to 0$,
 \beqnn
 G_2(z)= {\rm sgn}(z) G_2(|z|) \sim \Gamma(-\alpha)  \cos(\pi \alpha/2) \cdot {\rm sgn}(z) {\it\Phi}(1/|z|).
 \eeqnn
 Putting these two results together, we have  as $z\to 0$,
 \beqnn
 \int_0^\infty (e^{\mathtt{i}zt}-1- \mathtt{i}zt) d{\it\Phi }(t)
 \ar\sim\ar \big(-\sin(\pi \alpha/2)+\mathtt{i}\cdot {\rm sgn}(z)\cos(\pi \alpha/2)\big) \cdot \Gamma(-\alpha)  {\it\Phi}(1/|z|),
 \eeqnn
 which equals to the right-hand side of (\ref{eqn.30035}). 
 \qed

  \begin{proposition}\label{TauberPhi03}
 	If $ {\it\Phi }\in {\rm RV}^\infty_{-\alpha-1} $ with $\alpha\in(0,1)$, we have as $z \to 0$,
 	\beqlb\label{eqn.30036}
 	\int_0^\infty (1-e^{\mathtt{i}zt} ) d{\it\Psi }_1(t) \sim  -(\alpha+1)\Gamma(-\alpha) \cdot {\rm sgn}(z)\exp\big\{-\mathtt{i}\cdot{\rm sgn}(z)\frac{\pi}{2}  \alpha\big\}\cdot  |z|^{-1}{\it\Phi}(1/|z|) .
 	\eeqlb
 \end{proposition}
 \proof Like the proof of Proposition~\ref{TauberPhi02},  the left-hand side of (\ref{eqn.30036}) equals to $\int_0^\infty ( e^{\mathtt{i}zt} -1) d{\it\bar\Psi }_1(t)$ and can be rewritten into $G_3(z)+\mathtt{i}G_4(z)$ with 
 \beqnn
 G_3(z):=\int_0^\infty \big(  \cos(zt) -1 \big) d{\it\bar\Psi}_1(t)
 \quad\mbox{and}\quad
 G_4(z):=  \int_0^\infty \sin(zt) d{\it\bar\Psi}_1(t).
 \eeqnn
 Applying integration by parts to these two integrals, we have
 \beqnn
 G_3(z)= z\int_0^\infty \sin(zt) {\it\bar\Psi}_1(t)dt
 \quad \mbox{and}\quad
 G_4(z)= -z\int_0^\infty \cos(zt) {\it\bar\Psi}_1(t) dt.
 \eeqnn
 From (\ref{eqn.30037}) and Proposition~\ref{Thm.AbelThm}, we have as $z\to 0$,
 \begin{align*}
 G_3(z)  & \sim -(\alpha+1)\Gamma(-\alpha) {\rm sgn}(z) \cos(\pi\alpha/2) \cdot \frac{{\it\Phi}(1/|z|)}{|z|} ,\\
 G_4(z)  & \sim (\alpha+1)\Gamma(-\alpha) \sin(\pi\alpha/2)\cdot  \frac{{\it\Phi}(1/|z|)}{|z|}.
 \end{align*}
 The desired result follows by putting these two results together.
 \qed

 \begin{proposition}\label{Corollary.TauPhi01}
 	If $ {\it\Phi }\in {\rm RV}^\infty_{-\alpha-1} $ with $\alpha\in(0,1)$, then for any $K>0$, there exists a constant $C>0$ such that for any $ |z|,\lambda\in[0, K]$,
 	\beqlb\label{eqn.3.36}
 	\Big| \int_0^\infty \big(e^{(\mathtt{i}z -\lambda )t}-1-(\mathtt{i}z -\lambda )t\big)d{\it\Phi}(t) \Big| \leq C\big[ {\it\Phi}(1/|z|)+ \big(1+  |z|/\lambda \big)\cdot {\it\Phi}(1/\lambda )\big].
 	\eeqlb
 \end{proposition}
 \proof The   integral on the left-hand side of (\ref{eqn.3.36}) can be written as the sum of the following three terms: 
 \beqnn
 f_1(\lambda)\ar:=\ar \int_0^\infty \big(e^{-\lambda t  }-1+\lambda t\big)d{\it\Phi}(t),\cr
 f_2(z)\ar:=\ar \int_0^\infty \big(e^{\mathtt{i}z t }-1-\mathtt{i}z t\big)d{\it\Phi}(t),\cr
 f_3(\lambda,z)\ar:=\ar \int_0^\infty \big( e^{\mathtt{i}zt}-1\big)\big(e^{-\lambda t}-1\big) d{\it\Phi}(t).
 \eeqnn
 By Proposition~\ref{TauberPhi01} and \ref{TauberPhi02}, there exists a constant $C>0$ such that  for any $|z|,\lambda\in[0,K]$, 
 \beqnn
 |f_1(\lambda)|\leq C\cdot  {\it\Phi}(1/\lambda)
 \quad \mbox{and} \quad 
 |f_2(z)| \leq C\cdot {\it\Phi}(1/|z|). 
 \eeqnn
 Additionally, an application of the inequality $|e^{\mathtt{i}zt}-1|\leq  |zt|$ for any $z,t\in\mathbb{R}$ and Proposition~\ref{TauberPhi01} shows that  uniformly in $z\in\mathbb{R}$ and $\lambda \in[0,K]$, 
 \beqnn
 |f_3(\lambda,z)| \ar\leq\ar |z|\int_0^\infty  \big(e^{-\lambda t}-1\big) t d{\it\Phi}(t) = |z|\int_0^\infty  \big(1-e^{-\lambda t}\big)  d{\it\Psi}_1(t) \leq C \cdot \frac{|z|}{\lambda} \cdot {\it\Phi}(1/\lambda).
 \eeqnn
 The desired result follows by putting these three estimates together.
 \qed

 \begin{proposition}\label{Corollary.TauPhi02}
 	If $ {\it\Phi }\in {\rm RV}^\infty_{-\alpha-1} $ with $\alpha\in(0,1)$, then  for any $\kappa>0$, there exist constants $C>0$, $\delta\in(0,1)$ and $n_0\geq 1$ such that for any $n\geq n_0$ and $|z|,\lambda\in[0,\delta n]$,
 	\beqnn
 	\Big| \int_0^\infty \Big(e^{\frac{\mathtt{i}z -\lambda }{n}t}-1-\frac{\mathtt{i}z -\lambda }{n}t\Big)d{\it\Phi}(t) \Big|
 	\ar\leq\ar C\big[ |z|^{\alpha+1+\kappa}\vee 1 +(1+|z|/\lambda) ( \lambda^{\alpha+1+\kappa}\vee 1)\big] \cdot {\it\Phi}(n).
 	\eeqnn
 \end{proposition}
 \proof By using Proposition~\ref{Corollary.TauPhi01} with $K=\delta$,  
 \beqlb\label{eqn.5004}
 \Big| \int_0^\infty \Big(e^{\frac{\mathtt{i}z -\lambda }{n}t}-1-\frac{\mathtt{i}z -\lambda }{n}t\Big)d{\it\Phi}(t) \Big| 
 \ar\leq\ar  C\big[ {\it\Phi}(n/|z|)+ \big(1+  |z|/\lambda \big)\cdot {\it\Phi}(n/\lambda )\big]\cr
 \ar=\ar C\Big[ \frac{{\it\Phi}(n/|z|)}{ {\it\Phi}(n)}+ \big(1+  |z|/\lambda \big)\cdot \frac{{\it\Phi}(n/\lambda )}{ {\it\Phi}(n)} \Big] \cdot {\it\Phi}(n),
 \eeqlb
 uniformly in $|z|,\lambda\in[0,\delta n]$.  
 The monotonicity of ${\it\Phi}$ induces that for any $|z|,\lambda\in[0,1]$,
 \beqnn
 \frac{{\it\Phi}(n/|z|)}{{\it\Phi}(n)} \vee \frac{{\it\Phi}(n/\lambda)}{{\it\Phi}(n)} \leq 1 . 
 \eeqnn
 Moreover, for any $\kappa >0$,  by Proposition~\ref{Thm.PotterThm} there exist two constants $C>0$ and $n_0\geq 1$ such that 
 \beqnn
 \frac{{\it\Phi}(n /|z|)}{{\it\Phi}(n)} \leq C |z|^{\alpha+1+\kappa}
 \quad \mbox{and}\quad
 \frac{{\it\Phi}(n/\lambda)}{ {\it\Phi}(n)} \leq C  \lambda^{\alpha+1+\kappa} ,
 \eeqnn
 uniformly in $n\geq n_0$ and $|z|,\lambda\in[1,\delta n ]$.
  The desired result follows immediately by plugging these three upper bound estimates into the right-hand side of the equality in (\ref{eqn.5004}).
 \qed

 We are ready to study the convergence rate of the  time-scaled resolvent  $R^{(n)}$ to the limit $R^*$.
 For each $\beta>0$,  by Corollary~\ref{ConvergenceR} we see that $R^{(n)}_\beta \to R^*_\beta$ in $M(\mathbb{R}_+;\mathbb{R}_+)$. 
 With a slight abuse of notation we write $d_{\mathcal{H}_2}(R^{(n)}_\beta,R^*_\beta)$ for the $2$-smooth Wasserstein distance between the two finite measures with density functions $R^{(n)}_\beta$ and $R^*_\beta$; see Appendix~\ref{Appendix.DM}.  

 \begin{lemma}\label{Lemma.01}
 	There exists a constant $C>0$ such that  for any $n\geq 1$ and $\beta\geq 1$,
 	\beqnn
 	d_{\mathcal{H}_2}(R^{(n)}_\beta,R^*_\beta)\leq 
 	\begin{cases}
 		\displaystyle{ C\cdot n{\it\Phi }(n),} & \mbox{if $ {\it\Phi }\in {\rm RV}^\infty_{-\alpha-1} $ with $\alpha\in(0,1)$;} \vspace{5pt}\\
 		\displaystyle {C\cdot n^{-1},} &  \mbox{if $ {\it\Psi }_2(\infty)<\infty$.}
 		\end{cases}
 	\eeqnn 
 	\end{lemma}
 	\proof By the definition of the metric $d_{\mathcal{H}_2}$ it is enough to prove that there exists a constant $C>0$ such that for any $n\geq 1$, $\beta \geq 1/2$ and $h\in\mathcal{H}_2$,
 	\beqlb\label{eqn.3.60}
 	\Big|\int_0^\infty h(t)\big( R^{(n)}_{2\beta}(t)-R^*_{2\beta}(t)\big)dt\Big|\leq 	\begin{cases}
 		\displaystyle{ C\cdot n{\it\Phi }(n),} & \mbox{if $ {\it\Phi }\in {\rm RV}^\infty_{-\alpha-1} $ with $\alpha\in(0,1)$;} \vspace{5pt}\\
 		\displaystyle {C\cdot n^{-1},} &  \mbox{if $ {\it\Psi }_2(\infty)<\infty$.}
 	\end{cases}
 	\eeqlb
 	
 	As preparation to prove (\ref{eqn.3.60}), we first introduce the function $h_{\beta } :=h \cdot \psi_{\beta }$, where $\psi_{\beta}$ is a smooth function on $\mathbb{R}$ that satisfies that $\psi_{\beta }(t) = e^{-\beta t}  $ if $ t\geq 0$ and
 	\beqnn
 	|\psi_{\beta} (t)| \leq   2 e^{-\beta |t|} ,\quad |\psi_{\beta }'(t)| \leq 4 e^{-\beta |t|} ,\quad |\psi_{\beta }'' (t)| \leq 8 e^{-\beta |t|},\quad t\in\mathbb{R}.
 	\eeqnn 
 	It is easy to see that $\|h_{\beta}\|_{L^1}\leq 8$, $\|h'_{\beta}\|_{L^1} \leq 32$ and $\|h''_{\beta}\|_{L^1} \leq 64$. 
 	Moreover, since $\sup_{n\geq 1}\|R^{(n)}_{\beta}\|_{L^1}+\|R^*_{\beta}\|_{L^1}<\infty$ (see (\ref{eqn.3.23})),
 	fundamental properties of Fourier transformation yield that
 	\beqlb\label{eqn.3.38}
 	\sup_{n\geq 1}\|\mathcal{F}R^{(n)}_{\beta}\|_{L^\infty}+ \|\mathcal{F}R^*_{\beta}\|_{L^\infty}\leq C_0
 	\quad \mbox{and}\quad
 	| \overline{\mathcal{F}h_{\beta}} (z)| \leq \frac{64}{1+|z|^2},
 	\eeqlb
 	for any $z\in\mathbb{R}$ and some constant $C_0>0$.
   Additionally, we will need the following bound for the term $\big| \mathcal{F}R^{(n)}_{\beta}(z)- \mathcal{F}R^*_{\beta}(z) \big|$. While 
   \beqlb\label{eqn.5007}
   \mathcal{F}R^*_{\beta}(z)= - \frac{1}{\sigma(\mathtt{i}z-\beta)},
   \eeqlb
   for $z\in\mathbb{R}$, 
   taking Fourier transforms on both sides of (\ref{eqn.3.61}) yields
   \beqnn
   \mathcal{F}R^{(n)}_{\beta}(z) =  \frac{\mathcal{F}\phi^{(n)}_\beta(z/n)  }{n(1-\mathcal{F}\phi^{(n)}_\beta(z/n) ) } = \frac{1}{n(1-\mathcal{F}\phi^{(n)}_\beta(z/n) ) }- \frac{1}{n},
   \eeqnn
   and hence
   \beqlb\label{eqn.5005}
   \big|\mathcal{F}R^{(n)}_{\beta}(z)-\mathcal{F}R^*_{\beta}(z) \big|\leq \frac{1}{n}+  \big|\varepsilon^{(n)}_R(\beta,z) \big|,
   \eeqlb
   with
   \beqlb\label{eqn.3.37}
   \varepsilon^{(n)}_R(\beta,z):=  \frac{n\int_0^\infty \big(  \exp\big\{\frac{\mathtt{i}z  -\beta}{n}t\big\}-1- \frac{\mathtt{i}z  -\beta}{n}t\big)d{\it\Phi}(t)}{n \big(1-\mathcal{F}\phi^{(n)}_\beta(z/n) \big) \cdot  \sigma (\mathtt{i}z-\beta)}.
   \eeqlb
 	
 We are now ready to prove (\ref{eqn.3.60}). 
 Using the Fourier inversion theorem and then Fubini's theorem, 
 we have $h_{\beta}(t)= \mathcal{F}\overline{\mathcal{F}h_{\beta}}(t)$
 for any $t\in\mathbb{R}$ and 
 	\beqlb\label{eqn.3.63}
 	\int_0^\infty h(t)\big( R^{(n)}_{2\beta}(t)-R^*_{2\beta}(t)\big)dt
 	\ar=\ar \int_0^\infty\big( R^{(n)}_{\beta}(t)-R^*_{\beta}(t)\big)h_\beta(t)dt\cr
 	\ar=\ar \int_{\mathbb{R}} \big( \mathcal{F}R^{(n)}_{\beta}(z)- \mathcal{F}R^*_{\beta}(z)\big)\cdot \overline{\mathcal{F}h_\beta}(z)dz .
 	\eeqlb
 To obtain (\ref{eqn.3.60}), we establish the desired upper bounds for the last integral in (\ref{eqn.3.63}) in various cases. 
 
 \textbf{Case I.} 
 Assume $ {\it\Phi }\in {\rm RV}^\infty_{-\alpha-1} $ with $\alpha\in(0,1)$. 
 For some constant $\theta\in(\alpha,1)$ to be determined later, the last integral in (\ref{eqn.3.63}) can be written as the sum of 
 \beqnn
 \varepsilon_{ 1}^{(n)} \ar:=\ar \int_{|z|\leq n^\theta} \big( \mathcal{F}R^{(n)}_{\beta}(z)- \mathcal{F}R^*_{\beta}(z) \big)\overline{\mathcal{F}h_\beta}(z) dz 
 \quad \mbox{and}\quad 
 \varepsilon_{ 2}^{(n)} := \int_{|z|> n^\theta} \big( \mathcal{F}R^{(n)}_{\beta}(z)- \mathcal{F}R^*_{\beta}(z) \big)\overline{\mathcal{F}h_\beta}(z) dz.
 \eeqnn
 We hence need to prove that 
 \beqnn
 \big| \varepsilon_i^{(n)} \big| \leq C\cdot n {\it\Phi}(n), \qquad i=1,2.
 \eeqnn
 
  We start with the second term.  
 Taking the two estimates in (\ref{eqn.3.38}) back into $\varepsilon_{ 2}^{(n)}$ and then using the fact that $n^{-\theta}=o(n {\it\Phi}(n))$, we have uniformly in $n\geq 1$ that
 \beqlb\label{eqn.3.62}
 | \varepsilon_{ 2}^{(n)}|\leq \int_{|z|> n^\theta} \frac{64C_0}{1+|z|^2}dz \leq \frac{C}{n^\theta} \leq C \cdot n {\it\Phi}(n).
 \eeqlb

 Next, we consider the term $ \varepsilon_{ 1}^{(n)}$.  
 Using inequalities in (\ref{eqn.3.38}) and (\ref{eqn.5005}), we have uniformly in $n\geq 1$  that
 \beqlb\label{eqn.5006}
  | \varepsilon_{ 1}^{(n)}|  \leq \frac{64}{n} \int_{|z|\leq n^\theta}\frac{dz}{1+|z|^2}
 + C_0\cdot \int_{|z|\leq n^\theta}\frac{ \big|\varepsilon^{(n)}_R(\beta,z) \big|}{1+|z|^2} dz  
 \leq \frac{C}{n} + C_0\cdot \int_{|z|\leq n^\theta}\frac{ \big|\varepsilon^{(n)}_R(\beta,z) \big|}{1+|z|^2} dz . 
 \eeqlb
 The numerator of $\varepsilon^{(n)}_R(\beta,z)$ can be estimated as follows.  For any $0<\kappa<\alpha\wedge(1-\alpha)\wedge (\alpha/\theta-\alpha)$, it follows from Proposition~\ref{Corollary.TauPhi02} that uniformly in $|z| \leq n^\theta$,
 \beqnn
 \Big|n\int_0^\infty \Big(  \exp\Big\{\frac{\mathtt{i}z  -\beta}{n}t\Big\} -1- \frac{\mathtt{i}z  -\beta}{n}t \Big)d{\it\Phi}(t)\Big|
 \ar\leq\ar   C \cdot \big( 1+ |z|^{\alpha+ \kappa +1}  \big) \cdot n {\it\Phi}(n).
 \eeqnn
 To estimate the denominator from below we use Proposition~\ref{Corollary.TauPhi02} again, to obtain constants $C_1>0$ and $n_0\geq 1$ such that for any $n\geq n_0$ and $|z| \leq n^\theta$, 
 \beqnn
 \big|n(1-\mathcal{F}\phi^{(n)}_\beta(z/n) )\big|\ar=\ar \Big|n\int_0^\infty\Big( \exp\Big\{\frac{\mathtt{i}z  -\beta}{n}t\Big\} -1 \Big)d{\it\Phi}(t) \Big|\cr
 \ar=\ar \Big| \sigma(\mathtt{i}z-\beta) + 
 n\int_0^\infty \Big( \exp\Big\{\frac{\mathtt{i}z  -\beta}{n}t\Big\}  -1-\frac{\mathtt{i}z  -\beta}{n}t\Big)d{\it\Phi}(t) \Big|\cr
 \ar\geq\ar \sigma \big|\mathtt{i}z-\beta\big| - \Big|   n\int_0^\infty \Big( \exp\Big\{\frac{\mathtt{i}z  -\beta}{n}t\Big\}  -1-\frac{\mathtt{i}z  -\beta}{n}t\Big)d{\it\Phi}(t) \Big|\cr
 \ar\geq\ar   \frac{ \sigma|z|+ \sigma\beta}{2} - C_1n {\it\Phi}(n /|z|)-C_1  (1+ |z|/\beta) \cdot n  {\it\Phi}(n/\beta) . 
 \eeqnn 
 It follows from the monotonicity of ${\it\Phi}$ that $\sup_{|z|\leq 1 } n{\it\Phi}(n /|z|) \leq n {\it\Phi}(n) \to 0$ as $n\to\infty$. Furthermore, 
 by Proposition~\ref{Thm.PotterThm},  there exist two constants $C>0$ and $n_0\geq 1$ such that for any $n\geq n_0$ and $1\leq |z|\leq n^{\theta}$,
 \beqnn
 n{\it\Phi}(n /|z|) \leq C\cdot |z|^{\alpha +1+\kappa} \cdot n{\it\Phi}(n)\leq C \cdot n^{1+\theta(\alpha+\kappa)}  {\it\Phi}(n)\cdot |z|.
 \eeqnn
 Since $n^{1+\theta(\alpha+\kappa)}  {\it\Phi}(n)\to 0$ and $\sup_{|z|\leq n^\theta}(1+ |z|/\beta) \cdot n {\it\Phi}(n/\beta) \leq C n^{1+\theta} {\it\Phi}(n)\to 0$  as $n\to\infty$,
 there exist constants $C>0$ and $n_1\geq n_0$ such that uniformly in $n\geq n_1$ and $|z| \leq n^\theta$,
 \beqlb\label{Upper.3.39}
 \big|n(1-\mathcal{F}\phi^{(n)}_\beta(z/n) )\big| \geq \frac{ \sigma|z|+ \sigma\beta}{4} 
 \quad \mbox{and hence}\quad 
 \big|\varepsilon^{(n)}_R(\beta,z)\big| \leq C \cdot \frac{ 1+ |z|^{\alpha+ \kappa +1} }{1+|z|^2}\cdot n {\it\Phi}(n). 
 \eeqlb 
 Plugging this back into the last integral in (\ref{eqn.5006}) induces  that uniformly in $n\geq 1$, 
 \beqnn
 \int_{|z|\leq n^\theta}\frac{ \big|\varepsilon^{(n)}_R(\beta,z) \big|}{1+|z|^2} dz  \leq C \cdot n {\it\Phi}(n)
 \quad \mbox{and hence}\quad 
  | \varepsilon_{ 1}^{(n)}|  \leq C \cdot n {\it\Phi}(n). 
 \eeqnn 
 
 {\bf Case~II.} Assume $ {\it\Psi }_2(\infty)<\infty$. 
 For any $\epsilon \in(0,1/4)$, we represent the last integral in (\ref{eqn.3.63}) as the sum of
 \beqnn
 \varepsilon_{3}^{(n)} \ar:=\ar \int_{|z|\leq \epsilon n} \big( \mathcal{F}R^{(n)}_{\beta}(z)- \mathcal{F}R^*_{\beta}(z) \big)\overline{\mathcal{F}h_\beta}(z) dz
 \quad \mbox{and}  \quad 
 \varepsilon_{4}^{(n)} := \int_{|z|>\epsilon n} \big( \mathcal{F}R^{(n)}_{\beta}(z)- \mathcal{F}R^*_{\beta}(z) \big)\overline{\mathcal{F}h_\beta}(z) dz.
 \eeqnn
 
  From (\ref{eqn.3.38}), we see that $ |\varepsilon_{4}^{(n)}| \leq C_0/(\epsilon n)$. To estimate the first term we recall  the definition of $\varepsilon^{(n)}_R(\beta,z)$ given in (\ref{eqn.3.37}) to get a constant $C>0$ that is independent of $n\geq 1$ and $\beta\geq 1/2$ such that
\beqlb\label{eqn.3.39}
 	\begin{split}
 		|\varepsilon_3^{(n)}|  \leq \frac{ C}{n} \int_{|z|\leq \epsilon n }\frac{ dz}{1+|z|^2}
 		+  C \int_{|z|\leq \epsilon n }\frac{|\varepsilon^{(n)}_R(\beta,z)|}{1+|z|^2} dz 
 		 \leq \frac{C}{n} +  C \int_{|z|\leq \epsilon n }\frac{|\varepsilon^{(n)}_R(\beta,z)|}{1+|z|^2} dz.
 	\end{split}
 \eeqlb
 
 Since $ {\it\Psi }_2(\infty)<\infty$ and $ | e^{z}-1-z  | \leq |z|^2$ for any $z\in \mathbb{C}_-$, there exists a constant $C>0$ such that for any $z\in\mathbb{R}$ and $n\geq 1$,
 \beqnn
 \Big|n\int_0^\infty \Big(  \exp\Big\{\frac{\mathtt{i}z  -\beta}{n}t\Big\} -1- \frac{\mathtt{i}z  -\beta}{n}t \Big)d{\it\Phi}(t)\Big|\leq \frac{|\mathtt{i}z  -\beta|^2}{n} \int_0^\infty t^2\phi(t)dt \leq \frac{C}{n}(|z|^2+\beta^2).
 \eeqnn
 This yields a constant $n_1>0$  such that for any $n\geq n_1$ and $|z|\leq \epsilon n$, 
 \beqnn
 \big| n\big(1-\mathcal{F}\phi^{(n)}_\beta(z/n) \big) \big|
 \ar=\ar \Big|  n\int_0^\infty \Big( \exp\Big\{\frac{\mathtt{i}z  -\beta}{n}t\Big\}-1 - \frac{\mathtt{i}z -\beta}{n}t\Big)d{\it\Phi}(t) - \sigma(\mathtt{i}z-\beta)\Big| \geq \frac{\sigma(|z|+\beta)}{2}. 
 \eeqnn 
 
 Taking these two upper bounds back into $\varepsilon^{(n)}_R(\beta,z)$ and using the inequality $|\mathtt{i}z-\beta|\geq (|z|+\beta)/2$ yields a constant $C>0$ such that for any $n\geq 1$ and $|z|\leq \epsilon n$,
 \beqlb\label{eqn.3.41}
 |\varepsilon^{(n)}_R(\beta,z)| \leq \frac{C}{n}
 \quad\mbox{and hence}\quad
 \int_{|z|\leq \epsilon n }\frac{|\varepsilon^{(n)}_R(\beta,z)|}{1+|z|^2} dz \leq  \frac{C}{n}\int_{|z|\leq \epsilon n }\frac{dz}{1+|z|^2}   \leq \frac{C}{n}.
 \eeqlb
 Taking this estimate back into (\ref{eqn.3.39}), we can get the desired result immediately.
 \qed

 \medskip
 {\it Proof for Proposition~\ref{Thm.ConvRateRLinfty}(i).}
 Using integration by parts,  we have for any $\beta>0$ and $t\geq 0$,
 \beqnn
 \mathcal{I}_{R^{(n)}}(t)\ar=\ar e^{\beta t} \mathcal{I}_{R^{(n)}_\beta}(t) - \int_0^t \beta e^{\beta s} \mathcal{I}_{R^{(n)}_\beta}(s)ds
 \quad\mbox{and}\quad
 \mathcal{I}_{R^*}(t)= e^{\beta t} \mathcal{I}_{R^*_\beta}(t) - \int_0^t \beta e^{\beta s} \mathcal{I}_{R^*_\beta}(s)ds,
 \eeqnn
 from which we deduce that
\[
	\|\mathcal{I}_{R^{(n)}}-\mathcal{I}_{R^*}\|_{L^\infty_T}\leq 2e^{\beta T} \cdot \|\mathcal{I}_{R^{(n)}_\beta}-\mathcal{I}_{R^*_\beta}\|_{L^\infty_T}\leq
 2e^{\beta T} \cdot \|\mathcal{I}_{R^{(n)}_\beta}-\mathcal{I}_{R^*_\beta}\|_{L^\infty}.
 \]
 
 We notice that $\|\mathcal{I}_{R^{(n)}_\beta}-\mathcal{I}_{R^*_\beta}\|_{L^\infty}$ equals the Kolmogorov distance  $d_{\rm Kol}(R^{(n)}_\beta,R^*_\beta)$ of two finite measures with density functions  $R^{(n)}_\beta$ and $R^*_\beta$; see Appendix~\ref{Appendix.DM}. 
 Using Proposition~\ref{Prop.DM} and the fact that $\|R^*_\beta\|_{L^\infty}\leq 1/\sigma$,  there exists a constant $C>0$ such that for any $n\geq 1$ and $\beta\geq 1$,
 \beqnn
 d_{\rm Kol}(R^{(n)}_\beta,R^*_\beta) \leq C \cdot \big(\mathrm{d}_{\mathcal{H}_2}(R^{(n)}_\beta,R^*_\beta) \big)^{1/3}  
 \quad \mbox{and}\quad 
 \|\mathcal{I}_{R^{(n)}}-\mathcal{I}_{R^*}\|_{L^\infty_T}\leq C\cdot e^{\beta T}\cdot \big(\mathrm{d}_{\mathcal{H}_2}(R^{(n)}_\beta,R^*_\beta)\big)^{1/3}
 \eeqnn 
 and the desired upper bounds now follow from Lemma~\ref{Lemma.01}. 
 \qed

 \medskip
 {\it Proof for Proposition~\ref{Thm.ConvRateRLinfty}(ii).}
 For each $t\geq 0$ and $\beta\geq 1$, 
 \beqnn
 |R^{(n)}(t) -R^*(t)|= e^{\beta t} |R^{(n)}_\beta(t) -R^*_\beta (t)|
 \quad \mbox{and} \quad 
 \|R^{(n)} -R^* \|_{L^2_T} \leq e^{\beta T} \cdot \|R^{(n)}_\beta-R^*_\beta\|_{L^2}. 
 \eeqnn 
 An application of  (\ref{eqn.30022}) and the Plancherel theorem gives that $\|R^{(n)}_\beta-R^*_\beta\|_{L^2}=\|\mathcal{F}R^{(n)}_\beta-\mathcal{F}R^*_\beta\|_{L^2}$. We now distinguish two cases.

 \medskip
 {\bf Case~I.} If ${\it \Phi}\in {\rm RV}_{-\alpha-1}^\infty$ with $\alpha\in(0,1)$, then for some constant $\theta\in(\alpha,1)$ to be specified later, it holds that 
 \[
 	\|\mathcal{F}R^{(n)}_\beta-\mathcal{F}R^*_\beta\|_{L^2}^2= \varepsilon_5^{(n)} +\varepsilon_6^{(n)} 
\]	
with
 \beqnn
 \varepsilon_5^{(n)} := \int_{|z|\leq n^\theta} \big|\mathcal{F}R^{(n)}_\beta(z)-\mathcal{F}R^*_\beta(z)\big|^2 dz
 \quad\mbox{and}\quad
 \varepsilon_6^{(n)} := \int_{|z|> n^\theta } \big|\mathcal{F}R^{(n)}_\beta(z)-\mathcal{F}R^*_\beta(z)\big|^2 dz.
 \eeqnn
 By  (\ref{eqn.30022}) and (\ref{eqn.5007}), there exists a constant $C>0$ such that $|\varepsilon_6^{(n)}| \leq C\cdot n^{-\theta}$ for any $n\geq 1$.
 
 Estimating the first rem is more involved. We first recall the quantity $ \varepsilon^{(n)}_R(\beta,z)$ defined in (\ref{eqn.3.37}) and the inequality (\ref{eqn.5005}). 
 By the power inequality, 
 \beqlb\label{eqn.30023}
 \varepsilon_5^{(n)} \ar\leq\ar \frac{4}{n^{2-\theta}} +2 \int_{|z|\leq n^\theta }\big| \varepsilon^{(n)}_R(\beta,z) \big|^2 dz.
 \eeqlb
 For any constant $0<\kappa<\alpha\wedge(1-\alpha)\wedge (\alpha/\theta-\alpha)$ to be specified again later, it follows from (\ref{Upper.3.39}) that
 \beqlb\label{eqn.3.40}
 \int_{|z|\leq n^\theta }\big| \varepsilon^{(n)}_R(\beta,z) \big|^2 dz
 \leq C\cdot \int_{|z|\leq n^\theta } \Big|\frac{ 1+ |z|^{\alpha+ \kappa +1} }{1+|z|^2}\Big|^2 dz\cdot \big| n {\it\Phi}(n)\big|^2.
 \eeqlb

 \smallskip
 {\it Case~I.1.} If $\alpha\in(0,1/2)$, then we choose $\theta\in(2\alpha,1)$ and $0<\kappa<\alpha\wedge(1-\alpha)\wedge (\alpha/\theta-\alpha)\wedge (1/2-\alpha)$.
 A simple calculation shows that the integral on the right side of (\ref{eqn.3.40}) can be uniformly bounded.
 Putting these estimates together, there exists a constant $C>0$ such that for any $n\geq 1$,
 \beqnn
 \|R^{(n)}_\beta-R^*_\beta\|_{L^2}^2 \leq C\cdot (n^{-\theta}+n^{ \theta-2} +  | n {\it\Phi}(n) |^2) \leq C \cdot | n {\it\Phi}(n) |^2,
 \eeqnn
 and hence $\|R^{(n)} -R^* \|_{L^2_T} \leq C\cdot e^{\beta T} \cdot n {\it\Phi}(n)  $.

 \smallskip
 {\it Case~I.2} \  If  $\alpha\in[1/2,1)$, then we choose $\theta\in(\alpha\vee(1-2\ell),1)$ and $\kappa=\ell$. In this case, the integral on the right side of (\ref{eqn.3.40}) can be bounded by $C\cdot n^{2(\alpha+\ell)-1}$ uniformly in $n\geq 1$.
 Putting these estimates together, there exists a constant $C>0$ such that for any $n\geq 1$,
 \beqnn
 \|R^{(n)}_\beta-R^*_\beta\|_{L^2}^2 \leq C\cdot (n^{-\theta}+n^{ \theta-2} +  n^{2(\alpha+\ell)+1}|{\it\Phi}(n) |^2) \leq C \cdot n^{2(\alpha+\ell)+1}|{\it\Phi}(n) |^2, 
 \eeqnn
 and hence $\|R^{(n)} -R^* \|_{L^2_T} \leq C\cdot e^{\beta T} \cdot  n^{\alpha+\ell+1/2} {\it\Phi}(n) $.

 \medskip
 {\bf Case~II.} If ${\it\Psi}_2(\infty)<\infty$,  for some $\epsilon\in(0,1/4)$ we extress the term $\|\mathcal{F}R^{(n)}_\beta-\mathcal{F}R^*_\beta\|_{L^2}^2$ as the sum of
 \beqnn
 \varepsilon_7^{(n)} := \int_{|z|\leq \epsilon n} \big|\mathcal{F}R^{(n)}_\beta(z)-\mathcal{F}R^*_\beta(z)\big|^2 dz
 \quad\mbox{and}\quad
 \varepsilon_8^{(n)} := \int_{|z|> \epsilon n} \big|\mathcal{F}R^{(n)}_\beta(z)-\mathcal{F}R^*_\beta(z)\big|^2 dz.
 \eeqnn
Using  (\ref{eqn.30022}) and (\ref{eqn.5007}) to $\varepsilon_8^{(n)} $ and  the first inequality in (\ref{eqn.3.41}) to $\varepsilon_7^{(n)} $, we have uniformly in $n\geq 1$ that
 \beqnn
 |\varepsilon_8^{(n)}| \leq \frac{C}{n}
 \quad \mbox{and}\quad 
  \varepsilon_7^{(n)} \ar\leq\ar \frac{2}{\epsilon n } + \int_{|z|\leq \epsilon n}  \big| \varepsilon^{(n)}_R(\beta,z) \big|^2 dz \leq \frac{C}{n}. 
 \eeqnn 
 Putting these estimates together, there exists a constant $C>0$ such that  $\|\mathcal{F}R^{(n)}_\beta-\mathcal{F}R^*_\beta\|_{L^2}\leq C\cdot n^{-1/2}$ and then $\|\mathcal{F}R^{(n)} -\mathcal{F}R^* \|_{L^2_T}\leq C\cdot e^{\beta T}\cdot n^{-1/2}$ for any $n\geq 1$. 
 \qed

   \section{Strongly critical Hawkes processes}
 \label{Sec.StronglyHP}
 \setcounter{equation}{0}

 In this section, we prove Theorems~\ref{MainThm.05} and \ref{MainThm.06}. The proofs are based on the asymptotics of the average event rate and its integrals given in the following lemmas. 

  \begin{lemma}\label{Lemma2.11}
 	For any $T>0$, it holds that 
	\[
		\sup_{t\in[0,T]} \left|  \frac{\mathbf{E}[N(nt)]}{\mathcal{I}^2_R(n)}-\mu_0\cdot t^{\alpha+1} \right|\to 0 \quad \mbox{as} \quad n\to\infty.
	\]
 \end{lemma}
 \proof It has been shown in Remark~\ref{rem:mu} that $\mathbf{E}[N(nt)]  = \mu_0 \big(nt+ \mathcal{I}^2_R(nt)\big)$. Hence, 
 \beqnn
 \sup_{t\in[0,T]} \Big|  \frac{\mathbf{E}[N(nt)]}{  \mathcal{I}^2_R(n)}-\mu_0\cdot t^{\alpha+1} \Big| \leq
 \mu_0\cdot \frac{nT}{\mathcal{I}^2_R(n)}  + \mu_0\cdot \sup_{t\in[0,T]}\Big|
 \frac{\mathcal{I}^2_R(nt)}{\mathcal{I}^2_R(n)}  -  t^{ \alpha +1} \Big|.
 \eeqnn
 In view of equation (\ref{eqn.3.21}), the first term on the right side of above inequality vanishes as $n\to\infty$. By Proposition~\ref{Thm.AsymR}(3) the function  $\mathcal{I}^2_R$ belongs to ${\rm RV}^{\infty}_{\alpha+1}$, and so it follows from Proposition~\ref{Thm.UniConver} that  the second term tends to $0 $ as well as $n\to\infty$.
 \qed

  \begin{lemma}\label{Lemma2.12}
 	For each $T>0$, as $n\to\infty$,  
 	\beqnn
 	\sup_{t\in[0,T]} \frac{1}{\mathcal{I}^2_R(n)}\cdot \big|\mathcal{I}_\Lambda(nt)-\mathbf{E} [\mathcal{I}_\Lambda(nt)] \big|  \to 0  \quad \mbox{in $\mathbf{P}$ and in $L^2(\mathbf{P})$} .
 	\eeqnn
 \end{lemma}
 \proof  Lemma~\ref{MartRep} along with the representation (\ref{IntegralLambda}) with $\mu\equiv \mu_0$ yields that 
 \[
 	\mathcal{I}_\Lambda(nt)-\mathbf{E} [\mathcal{I}_\Lambda(nt)]= \int_0^{t} R(n(t-s))  \widetilde{N}(ns) ds. 
\]	
	Using the relation (\ref{eqn.3.21}) in the last step, we hence get that
 \beqnn
 	\sup_{t\in[0,T]} \frac{1}{\mathcal{I}^2_R(n)}\cdot \big|\mathcal{I}_\Lambda(nt)-\mathbf{E} [\mathcal{I}_\Lambda(nt)] \big|  \ar\leq\ar \sup_{t\in[0,T]} \frac{n}{\mathcal{I}^2_R(n) }\int_0^{t} R(n(t-s)) \big|\widetilde{N}(ns)\big|ds \\
 	\ar \leq \ar \sup_{t\in[0,T]}\Big|\frac{\widetilde{N}(nt)}{n}\Big| \cdot \frac{n\mathcal{I}_R(nT) }{\mathcal{I}^2_R(n)} \\
	\ar \leq \ar C \sup_{t\in[0,T]}\Big|\frac{\widetilde{N}(nt)}{n}\Big|.
  \eeqnn
 Similarly as in the proof of Proposition~\ref{MonementEsti}, by using the Burkholder-Davis-Gundy inequality (e.g. see Theorem 26.12 in \cite[p.524]{Kallenberg2002}) together with Lemma~\ref{Lemma2.11}, we thus have that 
 \beqnn
 \mathbf{E}\Big[\sup_{t\in[0,T]}\big| \widetilde{N} (nt)/n \big|^2\Big] \ar\leq\ar \frac{C}{n^2} \cdot\mathbf{E}\big[N(nT) \big] \leq \frac{C}{n^2}\cdot\mathcal{I}^2_R(n) ,
 \eeqnn
 which tends $0$ as $n\to\infty$, from which the desired result follows.
 \qed

We proceed to analyze the asymptotics of the processes $\widetilde{N}$ and  $\mathcal{I}_\Lambda$.
The process $\widetilde{N}$ is an  $(\mathscr{F}_t)$-martingale with quadratic variation $\langle \widetilde{N}\rangle= \mathcal{I}_\Lambda$. In view of  Lemma~\ref{Lemma2.12} it seems natural to consider the following rescaled process:  
 \beqlb
 \widetilde{N}^{(n)}_R(t):= |\mathcal{I}_R^2(n)|^{-1/2}\cdot\widetilde{N}(nt)  ,\quad t\geq 0,\ n\geq 1.
 \eeqlb

 \begin{lemma}\label{Lemma2.13}
  As $n\to\infty$, we have that $\widetilde{N}^{(n)}_R \to G$ weakly in $\mathbf{D}(\mathbb{R}_+;\mathbb{R})$ where
  \beqnn
  G(t):= \sqrt{\mu_0(\alpha+1)} \int_0^t s^{\alpha/2}dB(s),\quad t\geq 0.
  \eeqnn
 \end{lemma}
 \proof The process $\widetilde{N}^{(n)}_R$ is an $(\mathscr{F}_{nt})$-martingale with quadratic variation
 \beqnn
 \langle\widetilde{N}_R^{(n)}\rangle_t= \frac{\mathcal{I}_\Lambda(nt)}{\mathcal{I}_R^2(n)},\quad t\geq 0
 \eeqnn
 that converges pointwise to $\mu_0 \cdot t^{\alpha+1}$ in probability as $n\to\infty$; see Lemma~\ref{Lemma2.12}. 
 
 By Theorem 4.13 in \cite[p.358]{JacodShiryaev2003}, the sequence $\{ \widetilde{N}^{(n)}_R\}_{n\geq 1} $ is tight;  it also $C$-tight, since the maximum jump size of $\widetilde{N}^{(n)}_R$ equals to $ |\mathcal{I}_R^2(n)|^{-1/2}$, which tends to $0$ as $n\to\infty$. 
 From the martingale representation theorem; see Theorem 7.1 in \cite[p.84]{IkedaWatanabe1989}, any limit process is equal in distribution to $G$.
 \qed

 \begin{lemma}\label{Lemma2.14}
 As $n\to\infty$, we have  weakly in $\mathbf{D}(\mathbb{R}_+,\mathbb{R})$ that
 \beqlb\label{eqn.30041}
 n|\mathcal{I}^2_R(n)|^{-3/2}\cdot \big(\mathcal{I}_\Lambda(nt)- \mathbf{E}[\mathcal{I}_\Lambda(nt)] \big)
 \to \sqrt{\mu_0(\alpha+1)}\int_0^t (t-s)^\alpha s^{\alpha/2}dB(s).
 \eeqlb
 \end{lemma}
 \proof
 By the Skorokhod representation theorem, Lemma~\ref{Lemma2.13} and the continuity of $G$, we may assume $\widetilde{N}^{(n)}_R \to G$ a.s. and uniformly on compact sets.  As argued in the proof of Lemma~\ref{Lemma2.12}, 
 \beqnn
 	\mathcal{I}_\Lambda(nt)-\mathbf{E} [\mathcal{I}_\Lambda(nt)]= \int_0^{t} R(n(t-s))  \tilde{N}(ns) ds 
 \eeqnn
and 
 \beqlb\label{eqn.3.08}
 n|\mathcal{I}^2_R(n)|^{-3/2}\cdot \big(\mathcal{I}_\Lambda(nt)- \mathbf{E}[\mathcal{I}_\Lambda(nt)] \big)
 \ar=\ar \frac{n^2}{\mathcal{I}^2_R(n)}\int_0^{t} R(n(t-s))  \widetilde{N}^{(n)}_R(s)ds \cr
   \ar=\ar  \frac{n^2}{\mathcal{I}^2_R(n)} \int_0^{t} R(n(t-s))  \big(\widetilde{N}^{(n)}_R(s)-G(s) \big) ds \cr
   \ar\ar + \frac{n^2}{\mathcal{I}^2_R(n)} \int_0^{t} R(ns)  G(t-s)ds.
 \eeqlb
 In view of (\ref{eqn.3.21}), the first term on the right-hand side of the second equality can be bounded uniformly in $t\in[0,T]$  by 
 \beqnn
	 \|\widetilde{N}^{(n)}_R-G\|_{L^\infty_T} \cdot  n\mathcal{I}_R(nT)/\mathcal{I}_R^2(n) \leq C \cdot \|\widetilde{N}^{(n)}_R-G\|_{L^\infty_T}
  \eeqnn
that tends to $0$ as $n\to\infty$. Let us then turn to the last term in (\ref{eqn.3.08}), which we denote by $M^{(n)}$ for convenience. It remains to prove the weak convergence of $M^{(n)} $ to the right-hand side of (\ref{eqn.30041}). 

To prove the tightness of the $\{M^{(n)}\}_{n\geq 1}$ we fix an $(\mathscr{F}_{nt})$-stopping time $\tau\in[0,T]$ and some $\delta\in(0,1)$. By (\ref{eqn.3.21}) we have that
 \beqnn
 \mathbf{E}\big[|M^{(n)}(\tau+\delta)-M^{(n)}(\tau)|\big] 
 \ar\leq\ar \frac{n^2}{\mathcal{I}^2_R(n)} \int_0^{T} R(ns)  \mathbf{E}\big[|G(\tau+\delta-s)-G(\tau-s)|\big]ds\cr
 \ar\leq\ar  \frac{n}{\mathcal{I}^2_R(n)} \int_0^{nT} R(s) ds \cdot \mathbf{E}[\zeta(\delta)] \\
\ar  \leq \ar  \frac{n\mathcal{I}_R(nT)}{\mathcal{I}^2_R(n)}  \cdot \mathbf{E}[\zeta(\delta)]  \\ 
\ar \leq \ar C \cdot \mathbf{E}[\zeta(\delta)] ,
 \eeqnn
 for some constant  $C>0$ that is independent of $(n,\tau)$, where 
 $$\zeta(\delta):= \sup_{s,t\in[0,T+1],|t-s|\leq \delta}|G(t)-G(s)|.$$
 Notice that $\sup_{\delta\in(0,1)}\zeta(\delta) \leq 2 \sup_{t\in[0,T]}|G(t)|$ and that the continuous martingale $G$ has predictable quadratic variation $\langle G \rangle_t =\mu_0(\alpha+1) \int_0^t s^\alpha ds $. 
 By the Burkholder-Davis-Gundy inequality; see Theorem 17.7 in \cite[p.333]{Kallenberg2002}, 
 \beqnn
 \mathbf{E}\Big[\sup_{\delta\in(0,1)}\zeta(\delta)\Big] \leq 2\cdot \mathbf{E}\Big[  \sup_{t\in[0,T]}|G(t)| \Big] \leq C \int_0^T s^\alpha ds =\frac{C}{1+\alpha}\cdot T^{\alpha+1}. 
 \eeqnn
 By the dominated convergence theorem and  the locally uniform continuity of $G$ it follows that
 \beqnn
 \lim_{\delta\to 0+} \mathbf{E}[\zeta(\delta)] = \mathbf{E}\Big[  \lim_{\delta\to 0+} \zeta(\delta)\Big] =0.
 \eeqnn
 The desired tightness of $\{M^{(n)}\}_{n\geq 1}$ thus follows from Aldous's criterion; see \cite{Aldous1978}. 
 
 To characterize the limit processes, we introduce the following two $\sigma$-finite measures on $\mathbb{R}_+$: 
 \beqnn
 \nu_n(ds):= n^2/\mathcal{I}^2_R(n)\cdot R(ns)ds
 \quad \mbox{and}\quad 
 \nu_*(ds):=\alpha s^{\alpha-1}ds. 
 \eeqnn 
 Then, $$M^{(n)}=G* \nu_n=G * \nu_* + G*(\nu_n-\nu_*).$$
 From Proposition~\ref{Thm.AsymR}(3) and the locally uniform continuity of $G$, for any $t\geq 0$, we know that $\nu_n \to \nu_*$ weakly on $(0,t]$ and  hence $G*(\nu_n-\nu_*) (t)\to 0$ a.s. as $n\to\infty$.
As a result, 
\[
 	M^{(n)}(t) \to G * \nu_*(t) \quad \mbox{a.s. for any $t\geq 0$} 
\]	
and by the stochastic Fubini lemma, 
 \beqnn
   G * \nu_* (t) =  \int_0^t (t-s)^\alpha dG(s) = \sqrt{\mu_0(\alpha+1)}\int_0^t (t-s)^\alpha s^{\alpha/2}dB(s). 
 \eeqnn 
 \qed

  \medskip
  \textit{Proof of Theorem~\ref{MainThm.05}.}
   From the martingale representation theorem (Lemma~\ref{MartRep}), we deduce that
  \beqnn
  n|\mathcal{I}^2_R(n)|^{-3/2}\cdot \big(N(nt)- \mathbf{E}[N(nt)] \big)
  \ar=\ar  \frac{n \cdot \tilde{N}^{(n)}_R(t)}{\mathcal{I}^2_R(n)} + \frac{n}{|\mathcal{I}^2_R(n)|^{3/2}} \cdot \big(\mathcal{I}_\Lambda(nt)- \mathbf{E}[\mathcal{I}_\Lambda(nt)] \big).
  \eeqnn
The second statement hence follows from Lemma~\ref{Lemma2.13} and \ref{Lemma2.14}.
 Moreover, by Lemma~\ref{Lemma2.11},
 \beqnn
 \sup_{t\in[0,T]}\big| N(nt)/\mathcal{I}^2_R(n) -  \mu_0\cdot t^{\alpha+1}    \big|
 \sim \frac{ |\mathcal{I}^2_R(n)|^{1/2}}{n}\sup_{t\in[0,T]}\Big| \frac{n}{|\mathcal{I}^2_R(n)|^{3/2}} \cdot \big(N(nt)- \mathbf{E}[N(nt)] \big)\Big|,
 \eeqnn
 as $n\to\infty$ and so the first statement follows directly from the second along with (\ref{eqn.3.21}).
 \qed

 \medskip

 \textit{Proof of Theorem~\ref{MainThm.06}.}
 By the martingale representation theorem (Lemma~\ref{MartRep}), we rewrite our normalized process as follows:
  \beqnn
 \gamma(n) \Big(\frac{N(nt)}{\mathcal{I}^2_R(n)}  -  \mu_0\cdot t^{\alpha+1} \Big)
  \ar=\ar \frac{\gamma(n)}{\mathcal{I}^2_R(n)}\big( N(nt)-\mathbf{E}[N(nt)] \big) + \gamma(n)\Big(  \frac{\mathbf{E}[N(nt)]}{\mathcal{I}^2_R(n)} - \mu_0\cdot t^{\alpha+1}\Big)\cr
  \ar=\ar \gamma(n)\frac{n }{\mathcal{I}^2_R(n)} \cdot \mu_0 t   +  \mu_0\cdot\gamma(n)A(n)\cdot \frac{1}{A(n)}\Big( \frac{\mathcal{I}^2_R(nt)}{\mathcal{I}^2_R(n)}  - t^{\alpha+1} \Big) \cr
  \ar\ar + \gamma(n)\frac{|\mathcal{I}^2_R(n)|^{1/2}}{n} \cdot  \frac{n}{|\mathcal{I}^2_R(n)|^{3/2}}\cdot \big(N(nt)-\mathbf{E}[N(nt)]\big) .
  \eeqnn
By assumption and Theorem~\ref{MainThm.05}(2), the first and third terms on the right side of the last equality converge to the corresponding limits uniformly on compact sets.
 
As for the second term, the assumptions $\mathcal{I}^2_R(n) \in {\rm 2RV}_{\alpha+1,\rho}^\infty(A)$ and $\gamma(n)A(n)\to \gamma_2^*$ imply that for any $t\geq 0$,
  \beqlb\label{eqn.3.06}
  \Big| \mu_0 \cdot \gamma(n)A(n)\cdot \frac{1}{A(n)}\Big( \frac{\mathcal{I}^2_R(nt)}{\mathcal{I}^2_R(n)}  - t^{\alpha+1} \Big) -\mu_0\cdot \gamma_2^*\cdot t^{\alpha+1} \int_1^t u^{\rho -1}du \Big|\to 0, 
  \eeqlb
  as $n\to\infty$.
  It remains to prove that this convergence also holds uniformly in $t\in[0,T]$ for any $T>0$.

 \medskip
 \textit{Case~I.} If $-\alpha-1<\rho\leq 0$, then a simple calculation shows that $\sup_{t\in[0,T]}t^{\alpha+1} \int_1^t u^{\rho -1}du <\infty$ and there exists a constant $C>0$ independent of $n$ and $t$ such that the left-hand side of (\ref{eqn.3.06}) can be bounded by
  \beqnn
  C \cdot \big|\gamma(n)A(n)  - \gamma_2^*\big|+
  C \cdot \Big|\frac{1}{A(n)}\Big( \frac{\mathcal{I}^2_R(nt)}{\mathcal{I}^2_R(n)}  - t^{\alpha+1} \Big) - t^{\alpha+1} \int_1^t u^{\rho -1}du \Big|.
  \eeqnn
By assumption (\ref{eqn.20021}) the first term vanishes as $n\to\infty$.  The second term can bounded as follows. 
  By Corollary~16.3.3 in \cite{Mao2013}, for $ \epsilon\in(0,1)$ and $\delta \in (0,\alpha+\rho+1)$,  
   there exists a constant $x_0>0$ such  that for any $n\geq 1$ and $t\in[x_0/n,T]$,
  \beqlb\label{eqn.3.09}
  \Big|  \frac{1}{A(n)}\Big( \frac{\mathcal{I}^2_R(nt)}{\mathcal{I}^2_R(n)} - t^{\alpha+1} \Big)- t^{\alpha+1} \int_1^t u^{\rho -1}du\Big|
  \leq C \cdot t^{\alpha+1}\cdot \big( t^{\rho-\delta}\vee t^{\rho+\delta}\big) \cdot \epsilon,
  \eeqlb
  which is bounded by $C (T^{ \alpha+\rho+\delta+1}\vee 1) \cdot\epsilon$ uniformly in $t\in[x_0/n,T]$. 
 Since 
 \beqnn
	A(n)\mathcal{I}^2_R(n)\to\infty \quad \mbox{and} \quad n^{\alpha+1}A(n) \to \infty, 
 \eeqnn
	we also have as $n\to\infty$ that
 \beqnn
   \sup_{t\in [0,x_0/n]} \Big|  \frac{1}{A(n)}\Big(\frac{\mathcal{I}^2_R(nt)}{\mathcal{I}^2_R(n)}   - t^{\alpha+1} \Big)- t^{\alpha+1} \int_1^t u^{\rho -1}du\Big|
   \ar\leq\ar \frac{C}{A(n)\mathcal{I}^2_R(n)} + \frac{C}{n^{2-\alpha}A(n)} + \frac{C\log n}{n^{2-\alpha}A(n)} \to 0. 
 \eeqnn
 Putting these two results together and then letting $\epsilon\to 0+$, we deduce that (\ref{eqn.3.06}) holds uniformly on any compact set.

 \medskip
 {\it Case~II.} If $\rho\leq -\alpha-1$, then using that $n|\mathcal{I}^2_R(n)|^{-1/2} \in {\rm RV}_{(1-\alpha)/2}^{\infty}$ and  $\mathcal{I}^2_R(n)/n\in {\rm RV}_{\alpha}^{\infty}$, we have as $n\to\infty$ that
 \beqnn
 n|\mathcal{I}^2_R(n)|^{-1/2} +\mathcal{I}^2_R(n)/n =o(1/A(n))
 \quad \mbox{and}\quad 
 \gamma(n)\in {\rm RV}_{\frac{1-\alpha}{2}\wedge  \alpha }^\infty,
 \eeqnn
from  which we deduce that $\gamma^*_2=0$.  
  Additionally, by (\ref{eqn.3.09}), there exits a constant $C>0$ such for any $n\geq 1$,
   \beqnn
  \sup_{t\in[x_0/n,T]}\frac{1}{A(n)}\Big|\frac{\mathcal{I}^2_R(nt)}{\mathcal{I}^2_R(n)} - t^{\alpha+1}\Big| \leq  \frac{C}{ n^{ \alpha+\rho -\delta+1}}.
   \eeqnn
  Since $\gamma(n)A(n)/n^{  \alpha+\rho -\delta+1 } \to 0$ as $n\to\infty$, we see that
   \beqnn
     \lim_{n\to\infty}\sup_{t\in[x_0/n,T]} \gamma(n) \Big| \frac{\mathcal{I}^2_R(nt)}{\mathcal{I}^2_R(n)}  - t^{\alpha+1} \Big|=0.
   \eeqnn
   On the other hand, there exists a constant $C>0$ such that for any $n\geq 1$,
   \beqnn
  \sup_{t\in[0,x_0/n]}\gamma(n) \Big| \frac{\mathcal{I}^2_R(nt)}{\mathcal{I}^2_R(n)} - t^{\alpha+1}\Big| \ar\leq\ar C \Big(\frac{\gamma(n)}{\mathcal{I}^2_R(n)} + \frac{\gamma(n)}{n^{\alpha+1}}\Big),
  \eeqnn
  which also goes to $0$ as $n\to\infty$.
  Putting  these two estimates together, one can see that (\ref{eqn.3.06}) holds uniformly on any compact set. 
  \qed

   \bigskip  \bigskip \bigskip  \bigskip
   \centerline{\Large\bf APPENDIX}

   \appendix

   \renewcommand{\theequation}{A.\arabic{equation}}
   
 \section{Distance between finite measures}\label{Appendix.DM}
 \setcounter{equation}{0}

 This section recalls several well-known distances on space of finite measures.
 A collection of Borel-measurable $\mathbb{C}$-valued functions on a topological space $(\mathbb{E},\mathscr{E})$, denoted by $\mathcal{S}(\mathbb{E})$, is said to be \textit{separating} if any two finite measures $\nu_1,\nu_2$ on $\mathbb{E}$ are equal if $\nu_1(f)=\nu_2(f)$ for any $f\in \mathcal{S}$. Examples of separating classes include the sets  $\mathcal{S}(\mathbb{E})=\{f:f \mbox{ is a bounded and measureable function on }\mathbb{E}\}$, $\mathcal{S}(\mathbb{E})=\{f:f\in C_b(\mathbb{E},\mathbb{R})\}$ and $\mathcal{S}(\mathbb{R}^d)= \{e^{\mathtt{i} \langle \lambda ,\cdot\rangle}:\lambda \in \mathbb{R}^d\}$ if $\mathbb{E}=\mathbb{R}^d$.

 \begin{definition} 
 	 For a separating class $\mathcal{S}(\mathbb{E})$ and two finite measures $\nu_1,\nu_2$ on $\mathbb{E}$ satisfying  that $|\nu_1(f)|+|\nu_2(f)|<\infty$ for any $f\in \mathcal{S}(\mathbb{E}) $, the {\rm $\mathcal{S}(\mathbb{E})$-distance} between $\nu_1$ and $\nu_2$ is given by 
 \beqlb\label{S-Distance}
 d_{\mathcal{S}(\mathbb{E})}(\nu_1,\nu_2):= \sup_{f\in \mathcal{S}(\mathbb{E})} \big|\nu_1(f)-\nu_2(f)\big|.
 \eeqlb

 \end{definition}

Due to the one-to-one corresponding between measures on $\mathbb{C}$ and $\mathbb{R}^2$, the $\mathcal{S}(\mathbb{C})$-distance between two measures on $\mathbb{C}$ is equivalent to the $\mathcal{S}(\mathbb{R}^2)$-distance between the two corresponding measures on $\mathbb{R}^2$. 
 Specially, for some $n\geq 1$, if finite the measures $\nu_1,\nu_2$ on $\mathbb{R}^n$ have density functions $\rho_1,\rho_2$, then we also write $ d_{\mathcal{S}(\mathbb{R}^n)}(\rho_1,\rho_2)$ for $ d_{\mathcal{S}(\mathbb{R}^n)}(\nu_1,\nu_2)$ for convention.
 We now recall several well-known distances on the space of finite measures on $\mathbb{R}^n$.

  \begin{enumerate}
   \item[$\bullet$]{$d_{\rm Kol}(\nu_1,\nu_2)$:} the {\bf Kolmogorov distance} of two finite measures $\nu_1$ and $\nu_2$ on $\mathbb{R}^n$ is obtained from (\ref{S-Distance}) by taking $\mathcal{S}(\mathbb{R}^n)$ to be the class of all functions $f:\mathbb{R}^n\to \mathbb{R}$ of the type $f(x_1,\cdots,x_n)= \prod_{i=1}^n\mathbf{1}_{(-,\infty,x_i]}(z_i)$ with $(z_1,\cdots,z_n)\in\mathbb{R}^n$, i.e.,
       \beqlb
       d_{\rm Kol}(\nu_1,\nu_2)
       \ar:=\ar \sup_{(z_1,\cdots,z_n)\in\mathbb{R}^n}\Big| \nu_1\Big(\prod_{i=1}^n(-\infty,x_i]\Big) - \nu_2\Big(\prod_{i=1}^n(-\infty,x_i]\Big) \Big|.
       \eeqlb

   \item[$\bullet$] {$d_{\rm W}(\nu_1,\nu_2)$:} the {\bf Wasserstein distance} of two finite measures $\nu_1$ and $\nu_2$ on $\mathbb{R}^n$ is obtained from (\ref{S-Distance}) by taking $\mathcal{S}(\mathbb{R}^n)=\mathcal{H}_{\rm Lip}(\mathbb{R}^n;\mathbb{R})$ that is the class of Lipschitz continuous functions $f:\mathbb{R}^n\to \mathbb{R}$ with Lipschitz constant $\|f\|_{\rm Lip} \leq 1$, where
       \beqnn
       \|f\|_{\rm Lip} := \sup_{x,y\in\mathbb{R}^n,x\neq y} \frac{|f(x)-f(y)|}{|z-y|}.
       \eeqnn

   \item[$\bullet$] {$d_{\mathcal{H}_k}(\nu_1,\nu_2)$:} the {\bf $k$-smooth Wasserstein distance} of two finite measures $\nu_1$ and $\nu_2$ on $\mathbb{R}^n$ for some $k\in\mathbb{Z}_+$ is obtained from (\ref{S-Distance}) by taking $\mathcal{S}(\mathbb{R}^n)=\mathcal{H}_k(\mathbb{R}^n;\mathbb{R})$ that is the class of functions $f \in C^k(\mathbb{R}^n;\mathbb{R})$  with  $\|\frac{\partial^i}{\partial x^i}f\|_{ L^\infty} \leq 1$ for $i=0,1,\cdots, k$,

  \end{enumerate}

 As a direct consequence of Proposition 1.2(2) in \cite{Ross2011} and Proposition~2.1 in \cite{GauntLi2022}, the next proposition deduces the Kolmogorov distance between two finite measures on $\mathbb{R}^n$ in terms of their Wasserstein distance and their $k$-smooth  Wasserstein distance.

  \begin{proposition}\label{Prop.DM}
  Let $\nu_1$ be a finite measure on $\mathbb{R}^n$ with bounded density function. 
  For any $K\geq 1$ and $k\in\mathbb{Z}_+$, there exists a constant $C>0$ such that for any finite measure $\nu_2$ on $\mathbb{R}^n$ with $\nu_2(\mathbb{R}^n)\leq K$,
   \beqnn
   d_{\rm Kol}(\nu_1,\nu_2) \leq C\cdot \min \Big\{ \big( d_{\rm W}(\nu_1,\nu_2) \big)^{1/2} , \big( d_{\mathcal{H}_k}(\nu_1,\nu_2) \big)^{\frac{1}{k+1}} \Big\}.
   \eeqnn

  \end{proposition}

   \renewcommand{\theequation}{B.\arabic{equation}}

 \section{Regular variation}\label{Appendix.RV}
 \setcounter{equation}{0}

 In this section we recall selected results for regularly varying functions. The reader may refer to the monographs \cite{BinghamGoldieTeugels1987,DeHaanFerreira2006,Resnick2007} for a detailed discussion of regular variation.
 
 \begin{proposition}[Uniform convergence theorem]\label{Thm.UniConver}
 If $F \in {\rm RV}_\alpha^\infty$ with $\alpha\in\mathbb{R}$, then (\ref{Def.RV}) holds uniformly in $x \in [a,b]$ if $\alpha=0$, uniformly in $x\in(0,b]$ if $\alpha >0$ and uniformly in $x\in [a,\infty)$ if $\alpha<0$ for each $0<a<b<\infty$.
 \end{proposition}

 \begin{proposition}[Potter's theorem]\label{Thm.PotterThm}
  If $F\in {\rm RV}_\alpha^\infty$ with $\alpha\in\mathbb{R}$, then for any $C,\delta>0$, there exist a constant $t_0>0$ such that 
  \[
  	\frac{F(x)}{F(y)} \leq C \max \Big\{ (x/y)^{\alpha+\delta},  (x/y)^{\alpha-\delta} \Big\} \quad \mbox{for any} \quad x,y\geq t_0.
\]
 \end{proposition}

 \begin{proposition}[Karamata's theorem]\label{Thm.Karamata}
 If $\alpha\geq -1$, then $F \in {\rm RV}_{\alpha}^\infty$ implies that
 \beqlb\label{KaraThm01}
 \int_{0}^t F(s)ds \in {\rm RV}_{\alpha+1}^\infty
 \quad\mbox{and}\quad
 \lim_{t\to\infty} \frac{tF(t)}{\int_{0}^t F(s)ds} = \alpha+1.
 \eeqlb
 If $\alpha<-1$ or if $\alpha=-1$ and $\int_{t_0}^\infty F(s)ds<\infty$ for some $t_0>0$, then $F\in {\rm RV}_{\alpha}^\infty$ implies that
 \beqlb\label{KaraThm02}
 \int_t^\infty F(s)ds \in {\rm RV}_{\alpha+1}^\infty
 \quad\mbox{and}\quad \lim_{t\to\infty} \frac{tF(t)}{\int_t^\infty  F(s)ds} = -\alpha-1.
 \eeqlb
 \end{proposition}

 \begin{proposition}[Karamata's Tauberian theorem]\label{Thm.KaramataTauberian}
If $F$ is locally integrable and if its Laplace-Stieltjes transform $\hat{F} $ is well-defined on $(0,\infty)$ 
 for some $\alpha>-1$, then the following holds: 
 \begin{enumerate}
 	\item[(1)] If $F\in {\rm RV}^\infty_\alpha$, then $\hat{F}(1/\lambda) \sim \Gamma(1+\alpha) F(\lambda)$ as $\lambda \to \infty$.
 	
 	\item[(2)] If $F$ is eventually monotone and $\hat{F}(1/\cdot) \in {\rm RV}^\infty_\alpha$, then $ \Gamma(1+\alpha) F(t)\sim \hat{F}(1/t)$ as $t \to \infty$.
 	
 	\end{enumerate}
 \end{proposition}

   \begin{proposition}[Theorem~4.10.3 in \cite{BinghamGoldieTeugels1987}]\label{Thm.AbelThm}
   If $F$  is monotone on $\mathbb{R}_+$ and belongs to $ {\rm RV}_{-\alpha}^\infty$ for some $\alpha\in (0,1)$, then the following hold: as $z\to 0+$,
   \beqnn
   \int_0^\infty \cos(zt) F(t)dt \ar\sim\ar \Gamma(1-\alpha) \sin (\alpha\pi/2)  \frac{F(1/z)}{z}
   \quad \mbox{and}\quad 
   \int_0^\infty \sin(zt) F(t)dt \sim \Gamma(1-\alpha) \cos (\alpha\pi/2)  \frac{F(1/z)}{z}.
   \eeqnn
	
   \end{proposition}

   \begin{proposition}[Theorem~8.1.6 in \cite{BinghamGoldieTeugels1987}]\label{Prop.A.6}
  If $F$ is a distribution function on $\mathbb{R}_+$ and $\ell\in {\rm RV}_0^\infty$,  for some $n\in\mathbb{Z}_+$ and $\alpha=n+\beta$ with $\beta\in[0,1]$, then the following staements are equivalent: as $t\to\infty$,
   \begin{enumerate}
   	\item[(1)] $\hat{F}(1/t)-\sum_{i=0}^n \frac{a_i}{i!} (-1/t)^{i} \sim (-1)^{n+1} t^{-\alpha} \ell(t) $ with $a_i= \int_0^\infty s^idF(s)$. 
   	
   	\item[(2)] $1-F(t) \sim (-1)^n t^{-\alpha}\ell(t)/\Gamma(1-\alpha)$ when $\beta\in(0,1)$; $\int_t^\infty s^n dF(s) \sim n!\, \ell(t)$ when $\beta=0$; $\int_0^t s^{n+1} dF(s) \sim (n+1)!\, \ell(t)$ when $\beta=1$.
   	
   	\end{enumerate}

   \end{proposition}

    \bibliographystyle{plain}

 \bibliography{Reference}

  	\end{document}